\documentclass[12pt]{article}

\usepackage{amssymb,amsmath,amsthm,amsfonts,mathrsfs}
\usepackage{subfigure,epsfig,psfrag,epstopdf}
\usepackage{fullpage,graphicx,color}
\usepackage{algorithm,algorithmic}
\usepackage{tabularx,array,booktabs,bm,multirow}
\usepackage{longtable}
\usepackage{cite,url}
\usepackage[bookmarks=false]{hyperref}
\usepackage[labelfont=bf]{caption}

\newcommand{\mJ}{\mathcal{J}}

\newcommand{\mM}{\mathcal{M}}

\newcommand{\mO}{\mathcal{O}}
\newcommand{\RR}{\mathbb{R}}
\newcommand{\ZZ}{\mathbb{Z}}

\def\bnabla{\boldsymbol{\nabla}}

\def\Dpartial#1#2{\dfrac{\partial #1}{\partial #2}}
\newcommand{\argmin}{\operatorname{argmin}}
\newcommand{\argmax}{\operatorname{argmax}}
\newcommand{\mean}{\operatorname{mean}}

\begin{document}

\title{Multiscale Optimization via Enhanced Multilevel PCA-based Control Space Reduction
  for Electrical Impedance Tomography Imaging}
\date{}
\author{{\bf \normalsize Maria M.F.M.~Chun} \\
{\it \small Department of Mathematical Sciences, Florida Institute of Technology} \and
{\bf \normalsize Briana L.~Edwards} \\
{\it \small Department of Mathematical Sciences, Florida Institute of Technology} \and
{\bf \normalsize Vladislav Bukshtynov}\footnote{Corresponding author: \url{vbukshtynov@fit.edu}} \\
{\it \small Department of Mathematics and Systems Engineering, Florida Institute of Technology,}\\
{\it \small Melbourne, FL 32901, USA}}
\maketitle

\begin{abstract}

An efficient computational approach for imaging binary-type physical properties suitable for various models
in biomedical applications is developed and validated. The proposed methodology includes gradient-based
multiscale optimization with multilevel control space reduction based on principal component analysis,
optimal switching between the fine and coarse scales, and their effective re-parameterization. The reduced
dimensional controls are used interchangeably at both scales to accumulate the optimization progress and
mitigate side effects. Computational efficiency and superior quality of obtained results are achieved through
proper communication between solutions obtained at the fine and coarse scales. Reduced size of control spaces
supplied with adjoint-based gradients facilitates the application of this algorithm to models of high complexity
and also to a broad range of problems in biomedical sciences and outside. The performance of the complete
computational framework is tested with 2D inverse problems of cancer detection by electrical impedance
tomography~(EIT) in applications to synthetic models and models based on real breast cancer images. The results
demonstrate the superior performance of the new method and its high potential for minimizing possibilities for
false positive and false negative screening and improving the overall quality of the EIT-based procedures in
medical practice.

{\bf Keywords:} PDE-constrained optimization $\circ$ gradient-based method $\circ$ control space reduction $\circ$
multiscale parameter estimation $\circ$ electrical impedance tomography $\circ$ cancer detection

\end{abstract}

\section{Introduction}
\label{sec:intro}

Electrical Impedance Tomography~(EIT) is a non-invasive medical imaging method where surface
electrodes attached to the skin above the human body parts apply small constant or alternating
currents \cite{Nishimura2020} to examine suspected tissues for abnormalities, e.g., cancer
development \cite{Zou2003Review,Brown2003,Adler2008,Holder2004,Uhlmann2009,Lionheart2004,
Abascal2011,Choi2007,Boverman2008,Liu2019}. Our current work focuses on improving the EIT-based
procedures used to obtain equivalently accurate and effective results compared to those of
computed tomography~(CT) and magnetic resonance imaging~(MRI). Despite widespread and successful
use for cancer detection and monitoring, both CT and MRI come with adverse effects caused by
ionizing radiation \cite{Power2016} and harmful chemical reactions from injected dyes
\cite{Rogosnitzky2016}, respectively. Both techniques are relatively expensive and accompanied
by risks and uncertainties, raising multiple concerns among patients and healthcare professionals.
Alternatively, EIT techniques provide a non-invasive imaging method with no damage to the body and
other benefits such as portability, low cost, and, most importantly, safety. Here, we aim to
improve this methodology by enabling the detection of defective (cancerous) tissues before they
become too aggressive and potentially lethal to increase the chances of curing and minimizing
damaging effects on the patient's body \cite{Lymperopoulos2017}.

Briefly, EIT uses an acclaimed fact that electrical properties, such as electrical conductivity or
permittivity, change if the body tissue status changes from healthy to cancer-affected
\cite{Holder2004,Brown2003,Cheney1999}. For example, studies revealed that cancerous breast tissues
have different conductivity when compared to non-malignant ones, an effect explained as a result of
the increased density of the tumor stroma \cite{Jossinet1998}. Differences in the electrical
conductivity scope between different types of tissues are the contrast source in the EIT images.
This physical phenomenon supports EIT techniques to produce images of biological tissues by interpreting
their response to applied voltages or injected currents. In turn, the inverse EIT problem reconstructs
the electrical conductivity or permittivity by measuring voltages or currents at electrodes placed
on the test volume surface. This so-called Calderon-type inverse problem \cite{Calderon1980} is highly
ill-posed; refer to the topical review paper \cite{Borcea2002}. In 1989, Cheng et al.~\cite{Cheng1989}
proposed a mathematical concept for solving EIT problems by performing both analytical and computational
analyses of such solutions. Since the 1980s, various computational techniques have suggested a range
of computational solutions to this inverse problem. We refer to the most recent papers
\cite{Adler2015,Bera2018,Wang2020} to review the current state of the art and the existing challenges
closely associated with EIT and its applications.

In this work, we develop and validate an efficient computational framework to offer the optimal
reconstruction of biomedical images based on measurements obtained with noise. We see this approach as
advantageous in various applications for medical practices dealing with models characterized mainly by
near-binary distributions, e.g., used to represent electrical conductivity. Similar to the prototype
computational algorithm proposed in \cite{KoolmanBukshtynov2021}, we supply the gradient-based multiscale
optimization with multilevel control space reduction that applies interchangeably to both fine and coarse
scales. In the upgraded model, proper intercommunication between these scales assures computational
efficiency and the superior quality of obtained results. To enhance the performance of this approach, we
developed an algorithm for multiscale optimization using various techniques for control space reduction
and assuming the identification of multiple regions as cancerous spots. The impact of the noise in
measurements and the application of regularization techniques are also systematically analyzed.

As proven computationally in many applications, fine-scale optimization performed on fine meshes allows for
obtaining high-resolution images. These fine meshes also contribute enormously toward increased sensitivity
by enforcing the accuracy of computed gradients if employed. The size of the control (optimization) space
defined over fine scales can significantly reduce after applying parameterization, e.g., by using linear
transformations based on available sample solutions (realizations) when applying principal component
analysis~(PCA). However, fine-scale optimization may still suffer from so-called over-parameterization if
the problem is under-determined, i.e., when the number of controls (optimization variables) overweighs the
available data, namely the number of measurements. On the other hand, optimization performed on coarse meshes
usually terminates much faster due to the size of the control space. However, such solutions are less accurate
due to the problems with the sensitivity naturally ``coarsened'' by the low-resolution meshes. In addition,
coarse-scale optimization may suffer from being over-determined if the sizes of the control space and data
in use are imbalanced; see, e.g., \cite{BukshtynovBook2023,Holder2004,AbdullaBukshtynovSeif2021}.

The proposed multiscale optimization framework utilizes the most advantages of using fine and coarse meshes
by switching periodically between fine and coarse scales to help mitigate their known side effects. For
example, images from fine-scale solutions may not provide clear boundaries between regions identified by
different physical properties in space. As a result, a smooth transition between healthy and cancer-affected
areas will prevent accurate recognition of shapes of the latter while solving the inverse problem of cancer
detection~(IPCD). Our new computational algorithm also includes a penalization (Tikhonov-type) approach to
help ``synchronize'' optimization at both scales; refer, e.g., to
\cite{BukshtynovBook2023,EnglHankeNeubauer1996,AbdullaBukshtynovSeif2021,Bukshtynov2011,Bukshtynov2013} for
more details on using regularization methods for solving inverse problems. From this point, we abandon the word
``mesh'' and use ``scale'' instead, as various parameterization approaches discussed in Section~\ref{sec:math}
will lead to control spaces of reduced sizes without any reference to physical meshes.

In a nutshell, fine-scale optimization approximates the locations characterized by high and low electrical
conductivity. Projecting solutions onto the coarse scales provides dynamical (sharp-edge) filtering to the
fine-scale images optimized further for better matching the available data. The filtered images projected
back onto the fine scales preserve some information on recent changes obtained at the coarse scales. In our
current work, we extended the computational efficiency of the procedure for automated scale shifting to
accumulate optimal progress obtained at both scales. At the extent of how we identify the boundaries of the
cancerous spots by utilizing projections between scales with assigned controls at the coarse scale, this
approach has some relation to a group of level--set methods that apply multiscale techniques and adaptive
grids \cite{OsherSethian1988,GibouFedkiwOsher2018, TsaiOsher2003,ChenChengLinWang2009,Liu2018,Lien2005,Tai2004}.
It also employs some concepts of multiscale parameter estimation; refer, e.g., to
\cite{Grimstad2000,Grimstad2003,Cominelli2007,Lien2005} for some details. Here, we keep the main focus on
applying this computational approach to IPCD by the EIT technique. However, there are no known restraints
for employing the same methodology to a large diversity of problems in biomedical sciences, physics, geology,
chemistry, and other fields.

This paper proceeds as follows. In Section~\ref{sec:math}, we present the mathematical description of the
inverse EIT as an optimization problem to be solved at both fine and coarse scales by applying control space
reduction using PCA (fine scale), upscaling via dynamical partitioning (coarse scale), switching between scales,
and penalization for the improved performance. Description of models and detailed computational results are
presented in Section~\ref{sec:results}. Concluding remarks are provided in Section~\ref{sec:remarks}.

\section{Mathematical Description}
\label{sec:math}

\subsection{Optimization Model for Inverse EIT Problem}
\label{sec:math_model}

In the recent papers \cite{AbdullaBukshtynovSeif2021,KoolmanBukshtynov2021,ArbicBukshtynov2022}, the
inverse EIT problem is formulated as a PDE-constrained optimization problem with extensive numerical
analysis for 2D models by implementing various methods for solution space re-parameterization including
the PCA coupled with dynamical control space upscaling. Our current discussion of the inverse EIT problem
will use similar notations. Let $\Omega \subset \RR^n, \ n = 2, 3$, be an open and bounded set (domain)
representing the medium of our particular interest. We assume that function $\sigma(x): \, \Omega
\rightarrow \RR_+$ represents (isotropic) electrical conductivity at point $x \in \Omega$. We also assume
that $m$ electrodes $(E_{\ell})_{\ell=1}^m$ with contact impedances $(Z_{\ell})_{\ell=1}^m \in \RR^m_+$
are attached to the periphery $\partial \Omega$ of domain $\Omega$. Here, we employ the so-called
``voltage--to--current'' model where constant voltages (electrical potentials) $U = (U_{\ell})_{\ell=1}^m
\in \RR^m$ are applied to electrodes $(E_{\ell})_{\ell=1}^m$ to initiate electrical currents
$(I_{\ell})_{\ell=1}^m \in \RR^m$ through the same electrodes. We assume that electrical currents and
voltages satisfy the conservation of charge and ground (zero potential) conditions, respectively
\begin{equation}
  \sum_{\ell=1}^m I_{\ell}=0, \qquad \sum_{\ell=1}^m U_{\ell}=0.
  \label{eq:curr_volt_conds}
\end{equation}

We formulate the inverse EIT (conductivity) problem \cite{Calderon1980} as a PDE-constrained optimization
problem by considering minimization of the following objective (function)
\begin{equation}
  \mJ(\sigma) = \sum _{\ell=1}^m \left( I_{\ell} - I_{\ell}^* \right)^2,
  \label{eq:obj_fn}
\end{equation}
where $I^* = (I_{\ell}^*)_{\ell=1}^m \in \RR^m$ are measurements made for electrical currents $I_{\ell}$
computed as
\begin{equation}
  I_{\ell} = \int_{E_{\ell}}  \sigma(x) \Dpartial{u(x)}{n} \, ds,
  \quad \ell = 1, \ldots, m
  \label{eq:el_current}
\end{equation}
based on conductivity field $\sigma(x)$ set here as a control (variable). The distribution of electrical
potential $u(x): \, \Omega \rightarrow \RR$ is obtained as a solution of the following (elliptic) PDE
problem
\begin{subequations}
  \begin{alignat}{3}
    \bnabla \cdot \left[ \sigma(x) \bnabla u(x) \right] &= 0,
    && \qquad x \in \Omega \label{eq:forward_1}\\
    \Dpartial{u(x)}{n} &= 0, && \qquad x \in \partial \Omega -
    \bigcup\limits_{\ell=1}^{m} E_{\ell}, \ \ell= 1, \ldots, m
    \label{eq:forward_2}\\
    u(x) +  Z_{\ell} \sigma(x) \Dpartial{u(x)}{n} &= U_{\ell},
    && \qquad x \in E_{\ell}, \ \ell= 1, \ldots, m
    \label{eq:forward_3}
  \end{alignat}
  \label{eq:forward}
\end{subequations}
where $n$ is an external unit normal vector on $\partial \Omega$. A well-known fact is that the inverse
EIT problem to identify electrical conductivity $\sigma(x)$ over the discretized domain $\Omega$ with
available input data $I^*$ of size $m$ is highly ill-posed. Therefore, the formulation of our optimization
problem has to be adapted to the situation when the size of input data is increased through additional
measurements while keeping the size of the unknown parameters, i.e., elements in the discretized description
for $\sigma(x)$, fixed. As detailed in \cite{AbdullaBukshtynovSeif2021}, we use a ``rotation scheme'' by
setting $U^1=U, \ I^1=I$ and considering $m-1$ new permutations of boundary voltages
\begin{equation}
  U^j = (U_j, \ldots, U_{m}, U_1, \ldots, U_{j-1}), \quad
  j = 2, \ldots, m
  \label{eq:permutations}
\end{equation}
applied to electrodes $E_1, E_2, \ldots, E_{m}$, respectively. Using the ``voltage--to--current'' model
allows us to measure associated currents $I^{j*} = (I^{j*}_1,\ldots,I^{j*}_{m})$ and further increase
the total number of available measurements from $m^2$ up to $K m^2$ by applying \eqref{eq:permutations}
to $K$ different permutations of potentials within set $U$. Having a new set of $K m$ input data
$(I^{j*})^{Km}_{j=1}$ and the Robin condition \eqref{eq:forward_3} used together with \eqref{eq:el_current},
we finally consider an updated optimization problem of minimizing a new objective function
\begin{equation}
  \mJ (\sigma) = \sum_{j=1}^{Km} \sum_{\ell=1}^m \beta^j_{\ell}
  \left[ \int_{E_{\ell}} \dfrac{U^j_{\ell}-u^j(x;\sigma)}{Z_{\ell}}
  \, ds - I^{j*}_{\ell} \right]^2,
  \label{eq:obj_fn_final}
\end{equation}
where each function $u^j(\cdot; \sigma), \ j = 1, \ldots, Km$, solves problem
\eqref{eq:forward_1}--\eqref{eq:forward_3}. Added weights $\beta^j_{\ell} \geq 0$ in \eqref{eq:obj_fn_final},
in general, allow setting the level of importance for measurements $I^{j*}_{\ell}$ (when $\beta^j_{\ell} > 0$)
or excluding those measurements ($\beta^j_{\ell} = 0$) from all computations related to objectives and
associated gradients. We also note that the (forward) EIT problem \eqref{eq:forward_1}--\eqref{eq:forward_3}
together with \eqref{eq:el_current} may be used to generate various model examples and obtain synthetic data
for inverse EIT problems to mimic cancer-related diagnoses seen in reality.

Finally, as proposed in \cite{AbdullaBukshtynovSeif2021}, the solution of the optimization problem
\begin{equation}
  \hat \sigma(x) = \underset{\sigma}{\argmin} \ \mJ(\sigma)
  \label{eq:minJ_sigma}
\end{equation}
to minimize objective function \eqref{eq:obj_fn_final} subject to PDE constraint \eqref{eq:forward} could be
obtained by an iterative algorithm utilizing adjoint gradients with respect to control $\sigma$
\begin{equation}
  \bnabla_{\sigma} \mJ = - \sum_{j=1}^{Km} \bnabla \psi^j(x)
  \cdot \bnabla u^j(x)
  \label{eq:grad_sigma}
\end{equation}
computed based on solutions $\psi^j(\cdot; \sigma): \, \Omega \rightarrow \RR, \ j = 1, \ldots, Km$, of the
following adjoint PDE problem:
\begin{equation}
  \begin{aligned}
    \bnabla \cdot \left[ \sigma(x) \bnabla \psi(x) \right] &= 0,
    && \quad x \in \Omega\\
    \Dpartial{\psi(x)}{n} &= 0, && \quad x \in \partial \Omega -
    \bigcup\limits_{\ell=1}^{m} E_{\ell}\\
    \psi(x) + Z_{\ell} \Dpartial{\psi(x)}{n} &=
    2 \beta_{\ell} \left[ \int_{E_{\ell}}
    \dfrac{u(x)-U_{\ell}}{Z_{\ell}} \, ds
    + I^*_{\ell} \right], && \quad
    x \in E_{\ell}, \ \ell= 1, \ldots, m
  \end{aligned}
  \label{eq:adjoint}
\end{equation}

\subsection{Fine-Scale Optimization via PCA-based Parameterization}
\label{sec:fine_scale}

The optimization problem \eqref{eq:obj_fn_final}--\eqref{eq:minJ_sigma}, as discussed in
Section~\ref{sec:math_model} and when solved for spatially discretized state variable $u(x)$ and control
$\sigma(x)$ (both use the similar piecewise constant representation), is usually over-parameterized.
To resolve the problem of ill-posedness associated with this
issue we apply commonly used re-parameterization of the control space ($\sigma$-space) based on principal
component analysis to represent $\sigma(x)$ in terms of uncorrelated variables (components of vector $\xi$)
mapping $\sigma(x)$ and $\xi$ by
\begin{subequations}
  \begin{alignat}{1}
    \sigma &= \Phi \, \xi + \bar \sigma, \label{eq:map_pca_1}\\
    \xi &= \hat \Phi^{-1} (\sigma - \bar \sigma). \label{eq:map_pca_2}
  \end{alignat}
  \label{eq:map_pca}
\end{subequations}
In \eqref{eq:map_pca}, $\Phi$ is the linear transformation matrix constructed using $N_r$ sample solutions
$(\sigma^*_n)_{n=1}^{N_r}$ (realizations discretized in the same way as control $\sigma$) as its columns, and
$\hat \Phi^{-1}$ denotes the pseudo-inverse of $\Phi$. The prior mean $\bar \sigma$ is given by
$\bar \sigma = (1/N_r) \sum_{n=1}^{N_r} \sigma^*_n$; see \cite{Bukshtynov2015,Jolliffe2002} for details on
constructing a complete PCA representation. The optimization problem \eqref{eq:minJ_sigma} is now restated in terms
of new model parameters $\xi \in \RR^{N_{\xi}}$, $1 \leq N_{\xi} \leq N_r$, used in place of control $\sigma(x)$
as follows
\begin{equation}
  \hat \xi = \underset{\xi}{\argmin} \ \mJ(\xi)
  \label{eq:minJ_xi}
\end{equation}
subject to discretized PDE model \eqref{eq:forward} with control mapping \eqref{eq:map_pca} used for
computing $\mJ(\xi) = \mJ(\sigma(\xi))$. For solving problem \eqref{eq:minJ_xi}, gradients $\bnabla_{\xi} \mJ$
of objective $\mJ(\sigma)$ with respect to new control $\xi$ can be expressed as
\begin{equation}
  \bnabla_{\xi} \mJ = \Phi^T \, \bnabla_{\sigma} \mJ
  \label{eq:grad_xi}
\end{equation}
to project gradients $\bnabla_{\sigma} \mJ$ obtained by \eqref{eq:grad_sigma} from initial (physical)
$\sigma$-space onto the reduced-dimensional $\xi$-space.

Finally, we note that the purpose of our research is to create an approach with minimal dependence on the
structure and source of the prior information. In addition, the vector subspace for control $\sigma$ spanned by
$\Phi$ must not be too restrictive in connection to this prior. Therefore, we generated a sufficient number of
random realizations $\sigma^*$ in the assumption that we do not know the number and location of the cancerous regions;
see the detailed description of this method in Section~\ref{sec:comp_model}. We also refer to \cite{VolkovBukshtynov2018}
for the discussion of how this PCA-based parameterization approach is consistent with treating the uncertainty
of control $\sigma$ using a classical Bayesian formulation.

\subsection{Coarse-Scale Partitioning}
\label{sec:coarse_scale}

Here, we briefly review the methodology of the control $\sigma$-space re-parameterization via
partitioning introduced in \cite{KoolmanBukshtynov2021} and substantially upgraded as described
in Section~\ref{sec:switching}. Generally speaking, this algorithm defines a new control space
with a reasonably small number of parameters (controls) $(\zeta_j)_{j=1}^{N_{\zeta}} \in
\RR^{N_{\zeta}}_+$. \cite{ChunMS2022} describes the general concept and various practical algorithms
by which spatial controls may be grouped (partitioned); Section~\ref{sec:comp_model} will provide details
on its use with our current computational framework.

The coarse-scale phase of the optimization framework includes the following steps.
\begin{enumerate}
  \item Discretize control $\sigma(x)$ in problem \eqref{eq:minJ_sigma} over the fine mesh with $N$
    elements each of area (or volume, in 3D) $\Delta_i$.
  \item Represent $\sigma(x)$ by a finite set of controls $(\sigma_i)_{i=1}^N \in \RR^N_+$.
  \item Partition this set into $N_{\zeta}$ subsets $C_j, \ j = 1, \ldots, N_{\zeta}$, by selecting
    (with no repetition) $N_j$ controls for the $j$th subset and defining a map (fine--to--coarse partition)
    \begin{equation}
      \begin{aligned}
        \mM: \, (\sigma_i)_{i=1}^N \rightarrow \bigcup
        \limits_{j=1}^{N_{\zeta}} C_j, \quad
        C_j &= \{ \sigma_i : P_{i,j} = 1, \, i = 1, \ldots, N\},\\
        \sum_{j=1}^{N_{\zeta}} |C_j| &= \sum_{j=1}^{N_{\zeta}} N_j = N,
      \end{aligned}
      \label{eq:upscale_map}
    \end{equation}
    supplied with the partition (indicator) function
    \begin{equation}
      P_{i,j} = \left\{
      \begin{aligned}
        1, \quad \sigma_i \in C_j,\\
        0, \quad \sigma_i \notin C_j.
      \end{aligned}
      \right.
      \label{eq:part_ind}
    \end{equation}
  \item Compute new (upscaled) gradients $\bnabla_{\zeta} \mJ$ by summing up those components
    $\Dpartial{\mJ}{\sigma_i}$ of discretized gradients $\bnabla_{\sigma} \mJ$ that are related
    to controls $\sigma_i \in C_j$, i.e.,
    \begin{equation}
      \Dpartial{\mJ}{\zeta_j} =
      \sum_{i=1}^{N} P_{i,j} \Dpartial{\mJ}{\sigma_i} \Delta_i.
      \label{eq:grad_zeta}
    \end{equation}
\end{enumerate}

Here, we note that a practical application of this method to solving the inverse EIT problem presented
in \cite{KoolmanBukshtynov2021} assumes $N_{\zeta} = 2$. This assumption uses a simplistic approach for
partitioning (subsets $C_1$ and $C_2$) based on the current values of control $\sigma (x)$ to distinguish
regions with low (subset $C_1$) and high (subset $C_2$) conductivity throughout the entire domain~$\Omega$.
It may cause poor computational performance for models featuring several cancer-affected areas of different
sizes and conductivities. As such, a new partitioning methodology and updated scheme for switching between
fine and coarse scales are developed to allow higher variations in the geometry of reconstructed binary
images discussed in the next section.

\subsection{Switching Between Scales}
\label{sec:switching}

\begin{figure}[!htb]
  \begin{center}
  \includegraphics[width=1.0\textwidth]{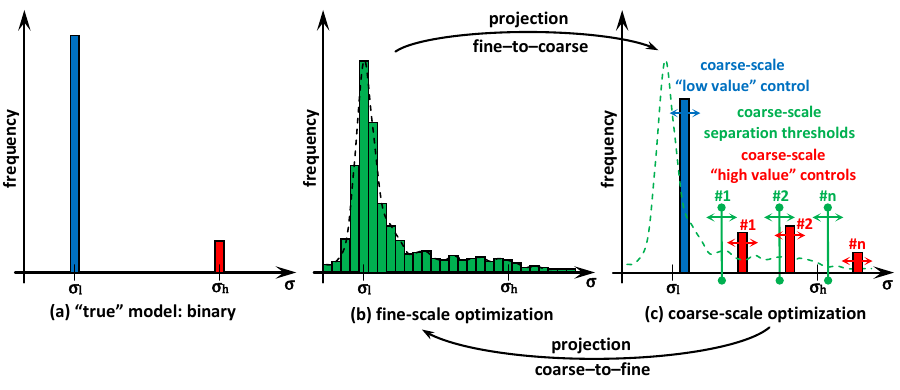}
  \end{center}
  \caption{A schematic illustration (adopted from~\cite{KoolmanBukshtynov2021}) of the general concept for the
    multiscale optimization framework modified to represent the proposed computational approach. In (a-c),
    $\sigma_l$ and $\sigma_h$ values represent two modes associated with healthy and cancer-affected regions
    within the domain $\Omega$, respectively. (a)~A typical histogram representing a binary distribution of true
    electrical conductivity $\sigma(x)$ used in EIT. (b)~An example of the Gaussian-type histogram typical for
    solution $\sigma^k(x)$ obtained after the $k$ iterations at a fine scale. (c)~A binary histogram for
    solution $\sigma^k(x)$ obtained after $k$th iteration at the coarse scale. Positions of blue and red bars
    are associated with current values of $\sigma^k_{low}$ and $\sigma^k_{high,n}$ ($1 \leq n \leq N_{\max}$)
    controls, and their heights are computed based on the fine-scale representation $\sigma(\xi^k)$ cut off
    by the current values of the coarse-scale separation threshold controls $\sigma^k_{th,n}$. See
    Section~\ref{sec:switching} for more details. Coarse--to--fine and fine--to--coarse projections are defined
    by \eqref{eq:proj_CtoF}--\eqref{eq:proj_CtoF_rlx} and \eqref{eq:sigma_coarse}--\eqref{eq:minJ_zeta_ini_2},
    respectively.}
  \label{fig:multiscale_cartoon}
\end{figure}

In the same fashion, as introduced in \cite{KoolmanBukshtynov2021}, the proposed computational approach
obtains the optimal solution $\hat \sigma(x)$ of a required binary type by employing multiscale optimization
at both fine and coarse scales, each with its own controls used interchangeably. We still consider a simplified
scheme for switching between scales: they are changed after completing $n_s$ optimization iterations. The
coarse-scale indicator function helps to account for these switches:
\begin{equation}
  \chi_c (k) = \left\{
  \begin{aligned}
    0, \quad (2k_s - 2) n_s &< k \leq (2k_s - 1) n_s,
    && \quad {\rm(fine \ scale)}\\
    1, \quad (2k_s - 1) n_s &< k \leq 2k_s n_s.
    && \quad {\rm(coarse \ scale)}
  \end{aligned}
  \right.
  \label{eq:scale_ind}
\end{equation}
In \eqref{eq:scale_ind}, $k_s = 1, 2, \ldots$ and $k = 0, 1, 2, \ldots$ denote the counts for switching cycles
and optimization iterations, respectively. We also use this function to define the termination condition
\begin{equation}
    \left| \dfrac{\mJ(\sigma^k) - \mJ(\sigma^{k-1})}{\mJ(\sigma^k)}
    \right| < (1-\chi_c) \epsilon_f + \chi_c \epsilon_c,
    \quad k \neq k_s n_s + 1
  \label{eq:termination}
\end{equation}
subject to chosen tolerances $\epsilon_f, \epsilon_c \in \RR_+$. The general concept of switching between
scales is shown schematically in Figure~\ref{fig:multiscale_cartoon}, including multiple controls (from $\#1$
to $\#n$) linked to solutions associated with different cancer-affected (with higher electrical conductivity)
regions; see Figure~\ref{fig:multiscale_cartoon}(c).

\subsubsection{Fine-Scale Optimization}
\label{sec:opt_fine}

While performing the fine-scale optimization ($\chi_c(k) = 0$), control $\sigma^k = \sigma^k(x)$ obtained after
the $k$th iteration as $\sigma(\xi^k)$ is updated by solving optimization problem \eqref{eq:minJ_xi} in the
reduced-dimensional $\xi$-space and by using map \eqref{eq:map_pca_1} as described in Section~\ref{sec:fine_scale},
i.e., $\sigma^k = \sigma(\xi^k)$. During the coarse-scale optimization phase ($\chi_c(k) = 1$), $\sigma(\xi^k)$
updated last time at the end of the fine-scale phase is used in partitioning discussed in Section~\ref{sec:coarse_scale}.
Apparently, updates made for fine-scale controls $\xi^k$ should ensure receiving as much information related to
recent changes in $\sigma^k$ during the coarse-scale phase as possible but not worsening the results $\sigma(\xi^k)$
previously obtained at the fine scale. To satisfy this requirement, solution $\sigma^k$, obtained at the end of the
coarse-scale phase, is projected onto $\xi$-space by using a convex combination of $\sigma(\xi^k)$ and $\sigma^k$, i.e.,
\begin{equation}
  \bar \sigma(\xi^k) = \alpha_{c \rightarrow f} \, \sigma(\xi^k)
  + (1 - \alpha_{c \rightarrow f}) \, \sigma^k, \quad
  \alpha_{c \rightarrow f} \in [0,1]
  \label{eq:proj_CtoF}
\end{equation}
followed by the re-initialization of control $\xi^k$ from $\bar \sigma(\xi^k)$ using map \eqref{eq:map_pca_2}.
As $\bar \sigma(\xi^k)$ and $\sigma^k$ have different (Gaussian and binary, respectively) distributions, the
coarse--to--fine projection \eqref{eq:proj_CtoF} also includes an extra step for projecting $\sigma^k$ to its
PCA equivalent
\begin{equation}
  \sigma^k_{PCA} = \Phi \hat \Phi^{-1} (\sigma^k - \bar \sigma)
  + \bar \sigma,
  \label{eq:proj_CtoF_pcaProj}
\end{equation}
before using it in \eqref{eq:proj_CtoF}; see \cite{VolkovBukshtynov2018,Bukshtynov2015} for details. An optimal
value of relaxation parameter $\alpha_{c \rightarrow f}$ is obtained by solving an additional 1D optimization
problem
\begin{equation}
  \alpha_{c \rightarrow f} =
  \hat \alpha =
  \underset{
  \begin{array}{cc}
     0 \leq \alpha \leq \alpha_{\max} \leq 1\\
     \mJ(\bar \sigma(\xi^k)) \leq \mJ(\sigma(\xi^k))
  \end{array}
  }{\argmin \ \alpha}
  \label{eq:proj_CtoF_rlx}
\end{equation}
that appeared to be highly nonlinear due to the inequality constraint to control the quality of fine scale
solutions $\sigma(\xi^k)$ in transition between subsequent switching cycles. We assume the tuning parameter
$\alpha_{\max}$ in \eqref{eq:proj_CtoF_rlx} depends on the problem. However, as shown in Section~\ref{sec:results},
it should not deviate too much from 1 to avoid enforced interventions into the ``natural flow'' of information
exchanged between the scales.

\subsubsection{Coarse-Scale Optimization}
\label{sec:opt_coarse}

The proposed procedure for running optimization at the coarse scale differs substantially from that used initially
in \cite{KoolmanBukshtynov2021} by assigning different controls to reconstruct $\sigma (x)$ at locations associated
with individual cancer-affected regions. The image quality (both shape and associated values of $\sigma (x)$ inside)
for each region depends on its size and the distance to the measuring electrodes. As such, the (cumulative) sensitivity
of objective $\mJ(\sigma)$ (with respect to changes in the part of control $\sigma (x)$ that is related to the location
of this region) varies between the regions, sometimes by several orders of magnitude.

First, we must specify the maximum number of expected cancer-affected (high conductivity) regions $N_{\max}$. By
considering the healthy part (low conductivity region) of domain $\Omega$ as a single region, partitioning
\eqref{eq:upscale_map}-\eqref{eq:part_ind} will attempt to create $N_{\zeta} = N_{\max} + 1$ subsets $C_j, \ j = 1,
\ldots, N_{\max}+1$. At the coarse scale, we define a new control vector $\zeta = (\zeta_j)_{j=1}^{2N_{\max}+1}$ in which
the first entry is the low value of (binary) electrical conductivity $\sigma (x)$ associated with a healthy region in
domain $\Omega$. The next $N_{\max}$ controls are the high values of $\sigma (x)$ related to areas affected by cancer,
i.e.,
\begin{equation}
  \zeta_1 = \sigma_{low}, \quad \zeta_2 = \sigma_{high,1},
  \quad \zeta_3 = \sigma_{high,2}, \quad \ldots, \quad
  \zeta_{N_{\max}+1} = \sigma_{high,N_{\max}}.
  \label{eq:ctrls_low_high}
\end{equation}
These controls are shown schematically in Figure~\ref{fig:multiscale_cartoon}(c) as one blue and multiple red bars,
respectively. The rest $N_{\max}$ components
\begin{equation}
  \zeta_{N_{\max}+2} = \sigma_{th,1}, \quad
  \zeta_{N_{\max}+3} = \sigma_{th,2}, \quad \ldots, \quad
  \zeta_{2N_{\max}+1} = \sigma_{th,N_{\max}}
  \label{eq:ctrls_thresh}
\end{equation}
take responsibility for the shape of those $N_{\max}$ cancerous regions. They are set as separation thresholds to
define boundaries between the low and high conductivity regions, as shown in green in Figure~\ref{fig:multiscale_cartoon}(c).
Such a structure of control $\zeta$ allows us to create a systematic representation of the coarse-scale solution
$\zeta^k$ for control $\sigma^k$ at the $k$th iteration based on the current fine-scale parameterization
$\sigma(\xi^k) = (\sigma_i(\xi^k))_{i=1}^N$, i.e.,
\begin{equation}
  \sigma^k_i = \left\{
  \begin{aligned}
    \sigma^k_{low},  & \quad \sigma_i(\xi^k) < \sigma^k_{th,n},\\
    \sigma^k_{high,n}, & \quad \sigma_i(\xi^k) \geq \sigma^k_{th,n},
  \end{aligned}
  \right. \quad i = 1, \ldots, N, \quad 1 \leq n \leq N_{\max}.
  \label{eq:sigma_coarse}
\end{equation}
Here, $n = n(i)$ denotes the number of a particular cancer-affected region defined subject to the partitioning map
$\mM^k$ currently established and used for the $k$th iteration. We also note that
\begin{equation}
  \begin{aligned}
    0 < \sigma^k_{low} &< \underset{1 \leq n \leq N_{\max}}{\min} \ \sigma^k_{high,n}, \\
    \underset{1 \leq i \leq N}{\min} \ \sigma_i(\xi^k) < \sigma^k_{th,n}
    &< \underset{1 \leq i \leq N}{\max} \ \sigma_i(\xi^k), \ n = 1,
    \ldots , N_{\max}.
  \end{aligned}
  \label{eq:sigma_coarse_bounds}
\end{equation}
Simply, \eqref{eq:sigma_coarse} provides a rule for creating fine--to--coarse partition $\mM^k$ in \eqref{eq:upscale_map}
where $N_{\zeta} = N_{\max} + 1$ based on the current state of control $\zeta^k$ (at the $k$th iteration). During the
coarse-scale optimization phase ($\chi_c(k) = 1$), control $\sigma^k$ is updated by solving the following
($2N_{\max}+1$)-dimensional optimization problem in the $\zeta$-space
\begin{equation}
  \hat \zeta =
  \underset{\zeta}{\argmin} \ \mJ(\zeta)
  \label{eq:minJ_zeta}
\end{equation}
subject to constraints (bounds) provided in \eqref{eq:sigma_coarse_bounds} and then setting $\sigma^k = \sigma(\zeta^k)$.
When solving problem \eqref{eq:minJ_zeta} during the first switching cycle (i.e., $k = n_s$), $\zeta^k$ could be initially
approximated by some constants, e.g.,
\begin{equation}
  \begin{aligned}
    \sigma^k_{th,n} &= \sigma_{ini} = \frac{1}{2}
    \left[ \underset{1 \leq i \leq N}{\max} \ \sigma_i(\xi^k) +
    \underset{1 \leq i \leq N}{\min} \ \sigma_i(\xi^k) \right],\\
    \sigma^k_{low} &= \underset{1 \leq i \leq N}{\mean}
    \left\{ \sigma_i(\xi^k) : \ \sigma_i(\xi^k) < \sigma_{ini} \right\},\\
    \sigma^k_{high,n} &= \underset{1 \leq i \leq N}{\mean}
    \left\{ \sigma_i(\xi^k) : P_{i,n+1} = 1, \ \sigma_i(\xi^k) \geq
    \sigma_{ini} \right\}, \quad n = 1, \ldots, N_{\max}.
  \end{aligned}
  \label{eq:minJ_zeta_ini_1}
\end{equation}
This initialization procedure for coarse-scale controls $\zeta^k$ is due to the practical approach used for
creating and updating maps $\mM^k$ \eqref{eq:upscale_map}; refer to \cite{ChunMS2022} for more details.
Fine--to--coarse switching when $k = (2k_s - 1) n_s, \ k_s > 1$, could be performed by utilizing the
corresponding values of control $\zeta$ obtained at the end of the previous coarse-scale phase, i.e.,
\begin{equation}
  \zeta^k = \zeta^{k-n_s-1}.
  \label{eq:minJ_zeta_ini_2}
\end{equation}
In fact, formulas \eqref{eq:sigma_coarse}--\eqref{eq:minJ_zeta_ini_2} provide a complete description of the
fine--to--coarse projection for control $\sigma(x)$ used in our approach.

Finally, while solving \eqref{eq:minJ_zeta} presumably by approaches that require computing gradients, their
first $N_{\max}+1$ components
\begin{equation}
  \Dpartial{\mJ(\zeta)}{\zeta_1} = \Dpartial{\mJ}{\sigma_{low}},
  \quad \Dpartial{\mJ(\zeta)}{\zeta_{n+1}} = \Dpartial{\mJ}{\sigma_{high,n}},
  \quad n = 1, \ldots N_{\max}
  \label{eq:grad_zeta_a}
\end{equation}
could be easily obtained by using the gradient summation formula \eqref{eq:grad_zeta} after completing the
partitioning map $\mM^k$ \eqref{eq:upscale_map}--\eqref{eq:part_ind} and employing \eqref{eq:sigma_coarse}.
On the other hand, the rest components may be approximated by a finite difference scheme, e.g., of the first
order:
\begin{equation}
  \begin{aligned}
    \Dpartial{\mJ(\zeta)}{\zeta_j} &=
    \Dpartial{\mJ}{\sigma_{th,n}} = \dfrac{1}{\delta_{\zeta}} \left[ \mJ
    \left(\sigma^k(\ldots,\zeta_j +
    \delta_{\zeta}, \ldots) \right) - \mJ
    \left( \sigma^k(\ldots,\zeta_j,\ldots)
    \right)\right] + \mO(\delta_{\zeta}),\\
    n &= 1, \ldots, N_{\max}, \quad j = n + N_{\max} +1.
  \end{aligned}
  \label{eq:grad_zeta_b}
\end{equation}
Parameter $\delta_{\zeta}$ in \eqref{eq:grad_zeta_b} is to be set experimentally, pursuing a trade-off
between being reasonably small to ensure accuracy and large enough to protect the gradient components
from being zero.

\subsubsection{Enhanced Scale Switching via Tuned PCA}
\label{sec:pca_switch}

We further improve the computational performance of our multiscale optimization algorithm by changing
the description of the fine-scale optimization space. We perform this dynamically by varying the parameter
$N^k_{\xi}$ used in place of fixed $N_{\xi}$ to govern the PCA-based parameterization \eqref{eq:map_pca}
applied on fine meshes. As discussed in Section~\ref{sec:kappa}, computational results evaluated at
different stages of the optimization process show the variance in the quality of obtained gradients
$\bnabla_{\xi} \mJ$ and the values of the objective function $\mJ(\xi)$ when the number of principal
components (size of vector $\xi$) changes. In practical computations, the complete PCA representation
uses the linear transformation matrix $\Phi$ in \eqref{eq:map_pca_1} constructed using the truncated
singular value decomposition~(TSVD) with the number of principal components reduced to $1 \leq N_{\xi}
\leq N_r$. Usually, $N_{\xi}$ is fixed to a constant number which is high enough to provide the
reduced-dimensional $\xi$-space with sufficient degrees of freedom to allow detailing of the obtained
solutions at a small scale. In our algorithm, however, it looks reasonable to change this parameter
following the solution updates at both fine and coarse scales. Below we describe the suggested procedures
to search for the optimal number of PCA components $1 \leq N_{\xi}^k \leq N_{\xi}$ performed during the
scale switching.

We start with fine--to--coarse switching ($k = (2k_s - 1) n_s$) with a straightforward idea of ``adjusting''
the fine-scale parameterization used to assist optimization performed at the coarse scale by ensuring a
better fit of the created binary image to the used data. Here, we suggest solving the following 1D integer
optimization problem
\begin{equation}
  N_{\xi}^k = \hat N_t =
  \underset{
  \begin{array}{cc}
     N_t \in \ZZ_{+}\\
     1 \leq N_t \leq N_{\xi}
  \end{array}
  }{\argmin \ \mJ(\sigma(\xi^k_{N_t}))}
  \label{eq:pca_FtoC}
\end{equation}
to minimize objective function $\mJ(N_t) = \mJ(\sigma(\xi^k_{N_t}))$ evaluated by using the part of the
fine--to--coarse projection described by \eqref{eq:sigma_coarse}, \eqref{eq:sigma_coarse_bounds}, and
\eqref{eq:minJ_zeta_ini_2}, and the current fine-scale control $\xi^k$ ``truncated'' in the following
way:
\begin{equation}
  \xi^k_{N_t} = \left[ \xi_1 \ \xi_2 \ \ldots \ \xi_{N_t} \ 0 \ \ldots \ 0 \right] \in \RR^{N_{\xi}}.
  \label{eq:xi_trunc}
\end{equation}
We note that \eqref{eq:xi_trunc} provides a computationally efficient approach to ``re-truncate'' PCA
without changing the structure of matrix $\Phi$ in \eqref{eq:map_pca_1} by removing its $N_{\xi}-N_t$ last
columns.

Next, during coarse--to--fine switching, we suggest choosing the optimal number of PCA components by another
integer optimization problem
\begin{equation}
  N_{\xi}^k = \hat N_t =
  \underset{
  \begin{array}{cc}
     N_t \in \ZZ_+\\
     1 \leq N_t \leq N_{\xi}\\
     \mJ(\sigma(\xi^k_{N_t})) \leq \mJ(\sigma(\xi^{k-n_s}))
  \end{array}
  }{\argmax \ N_t}
  \label{eq:pca_CtoF}
\end{equation}
solved before finding an optimal value of relaxation parameter $\alpha_{c \rightarrow f}$ by \eqref{eq:proj_CtoF_rlx}.
We need to explain the difference in structures of optimization problems \eqref{eq:pca_CtoF} and \eqref{eq:pca_FtoC}.
We note that the solution of \eqref{eq:pca_CtoF} is used to re-parameterize the fine-scale representation of control
$\sigma(\xi^k)$, which is heavily involved in the process of identifying the boundaries of regions with different
(averaged) values of $\sigma(x)$. Therefore, we prefer to keep the number of principal components in use as high as
possible (namely, to maximize it) to guarantee further progress with optimization at the fine scale. Also, condition
$\mJ(\sigma(\xi^k_{N_t})) \leq \mJ(\sigma(\xi^{k-n_s}))$ in \eqref{eq:pca_CtoF} ensures that the truncated version
$\xi^k_{N_t}$ of the control $\xi^{k-n_s}$ obtained at the end of the last fine-scale phase will not worsen the
fine-scale solution at iteration $k - n_s$.

Finally, we reiterate that optimization problems \eqref{eq:pca_CtoF} and \eqref{eq:proj_CtoF_rlx} are solved during
the same coarse--to--fine switching. The solution to the problem \eqref{eq:proj_CtoF_rlx}, relaxation parameter
$\alpha_{c \rightarrow f}$, defines the optimal proportion of the information obtained at the coarse scale that should
be added to the fine-scale solution. We consider this part of the ``communication'' established between scales as the
most important. As such, we see it logical to re-parameterize the fine scale by solving \eqref{eq:pca_CtoF} before
proceeding to \eqref{eq:proj_CtoF_rlx}. However, we acknowledge that other schemes to search for the optimal number
of PCA components may be designed to improve the computational performance of the proposed method. A complete
computational workflow to perform the described optimization over multiple scales is provided in
Algorithm~\ref{alg:main_opt} of Appendix~A.

\subsection{Coarse-Scale Regularization}
\label{sec:penal}

The performance of the proposed computational framework to perform optimization at multiple scales could be evaluated
by the accuracy in reconstructing electrical conductivity $\sigma(x)$ within the regions considered as affected by
cancer. We reiterate that the binary images of those regions are created based on two groups of controls, namely
\eqref{eq:ctrls_low_high} and \eqref{eq:ctrls_thresh}, reconstructed at the coarse scale. Those are the low or high
values of $\sigma(x)$ related to healthy and cancerous areas, respectively, and the controls involved in creating the
boundaries of the regions suspected of cancer. In practical applications, obtaining images with correct shapes is vital
for accurately locating and controlling the dynamics of cancer- and treatment-related processes. However, this type of
``shape'' optimization may lack sensitivity as the shapes are reconstructed by implicit interpretation of the fine-scale
images through their projection onto the coarse scales rather than using explicit parameterization.

To assist in developing the correct shapes while performing optimization \eqref{eq:minJ_zeta} during the coarse-scale
phase, we may assume that some prior knowledge exists for making predictions on the true values of controls $\sigma_{low}$
and $\sigma_{high}$ given by two constant values $\bar \sigma_l$ and $\bar \sigma_h$ ($\bar \sigma_l < \bar \sigma_h$),
respectively. Therefore, we define a penalization term
\begin{equation}
  \mJ_c = \chi_c (k) \ \beta_c \left[ (\zeta_1 - \bar \sigma_l)^2
  + \sum_{n=1}^{N_{\max}} (\zeta_{n+1} - \bar \sigma_h )^2 \right],
  \label{eq:reg_c}
\end{equation}
where $\beta_c \in \RR^+$ is an adjustable parameter. $\mJ_c$ augments the (core) objective function $\mJ$ given by
\eqref{eq:obj_fn_final} with a new term $\mJ_c$, i.e.,
\begin{equation}
  \bar \mJ = \mJ + \mJ_c
  \label{eq:obj_reg_full}
\end{equation}
being active while performing optimization during the coarse-scale phase ($\chi_c (k) = 1$) only. The structure of the
new objective function $\bar \mJ$ allows the evaluation of corresponding gradients
\begin{equation}
  \begin{aligned}
    \chi_c (k) &= 0: \quad &\bnabla_{\sigma} \bar \mJ &= \bnabla_{\sigma} \mJ,\\
    \chi_c (k) &= 1: \quad &\Dpartial{\bar \mJ(\zeta)}{\zeta_1} &=
    \Dpartial{\mJ(\zeta)}{\zeta_1} + 2 \beta_c (\zeta_1 - \bar \sigma_l),\\
    &&\Dpartial{\bar \mJ(\zeta)}{\zeta_{n+1}} &=
    \Dpartial{\mJ(\zeta)}{\zeta_{n+1}} + 2 \beta_c (\zeta_{n+1} - \bar \sigma_h),
    \quad n = 1, \ldots, N_{\max}
    \end{aligned}
  \label{eq:obj_grad_reg}
\end{equation}
to support the same multiscale computational framework discussed in Sections~\ref{sec:math_model}--\ref{sec:switching}.
Finally, we notice that this (Tikhonov-type) penalization has proven to have an additional effect of regularizing the
reconstruction procedure against noise possibly contained in the measured data \cite{Bukshtynov2011,Bukshtynov2013}.

\section{Computational Results}
\label{sec:results}

\subsection{Computations for Models in 2D}
\label{sec:comp_model}

\begin{figure}[!htb]
  \begin{center}
  \mbox{
  \subfigure[model~\#1]
    {\includegraphics[width=0.33\textwidth]{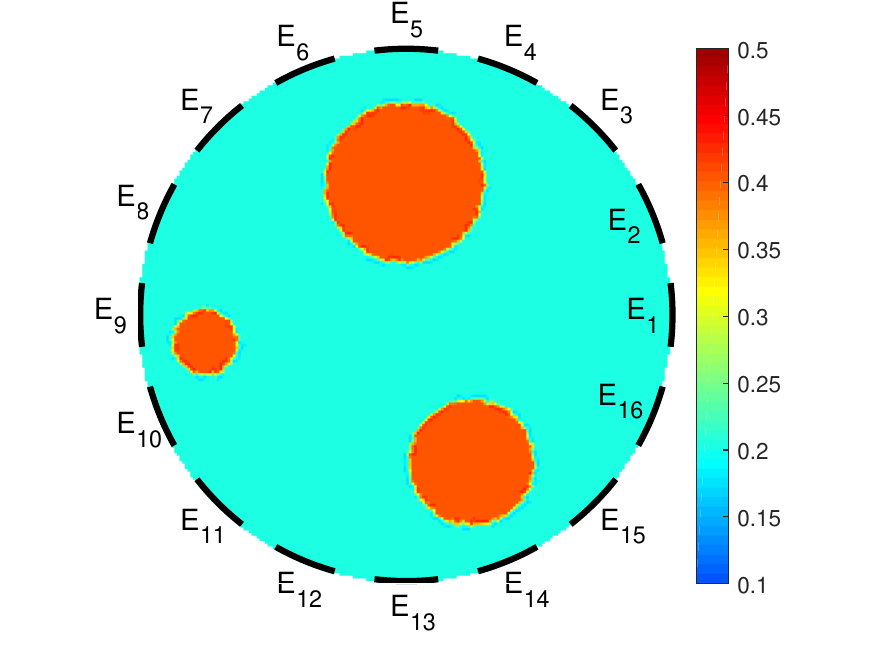}}
  \subfigure[$U^1$ applied]
    {\includegraphics[width=0.33\textwidth]{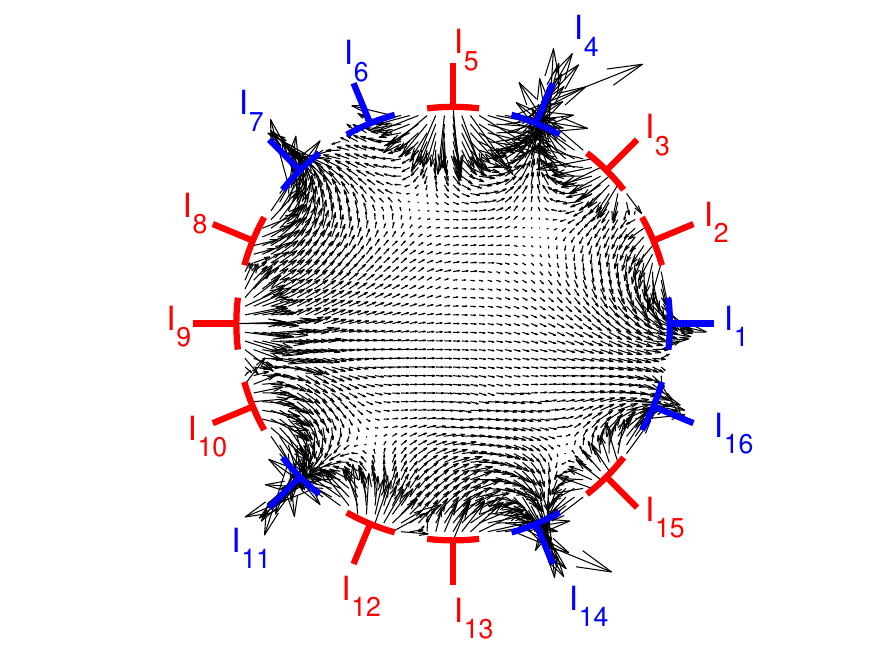}}
  \subfigure[$U^2$ applied]
    {\includegraphics[width=0.33\textwidth]{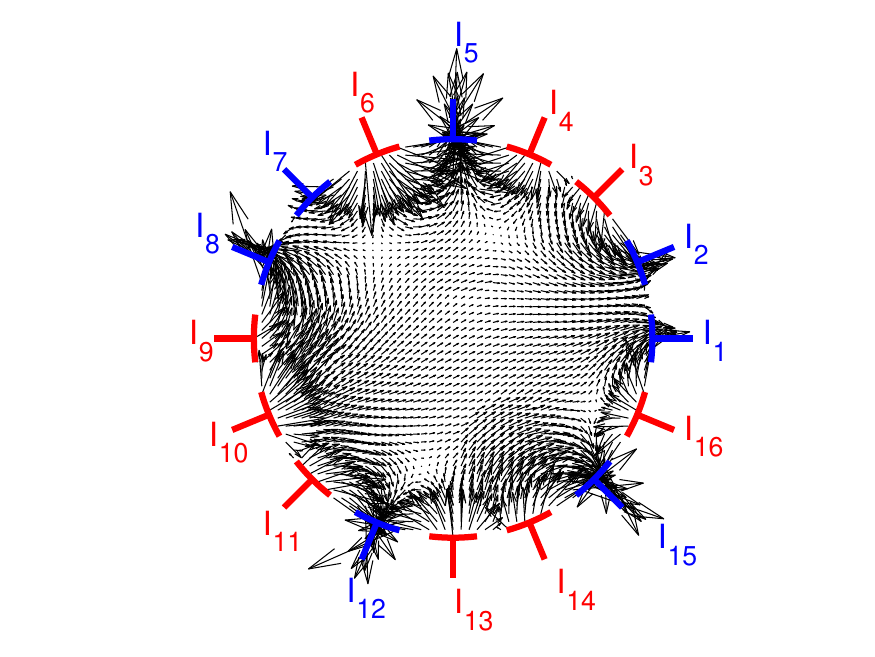}}}
  \end{center}
  \caption{(a)~EIT model~\#1: true electrical conductivity $\sigma_{true}(x)$ and equispaced geometry of electrodes $E_{\ell}$
    placed over boundary $\partial \Omega$. (b,c)~Electrical currents $I_l$ (positive in red, negative in blue) induced at
    electrodes $E_{\ell}$. Black arrows show the distribution of flux $\sigma(x) \bnabla u(x)$ of electrical potential $u$ in
    the interior of domain $\Omega$.}
  \label{fig:geometry}
\end{figure}

The computational part of the complete optimization framework integrates facilities for solving the EIT problem
\eqref{eq:forward}, adjoint problem \eqref{eq:adjoint}, and evaluation of the gradients according to \eqref{eq:grad_sigma},
\eqref{eq:grad_xi}, \eqref{eq:grad_zeta_a}--\eqref{eq:grad_zeta_b}, and \eqref{eq:obj_grad_reg}. These facilities are
incorporated using {\tt FreeFEM} \cite{FreeFem2012}, an open--source, high--level integrated development environment for
obtaining numerical solutions for PDEs based on the finite element method~(FEM). Numerical solutions for forward and
adjoint PDE problems assume spatial discretization implemented using FEM triangular finite elements. We apply P2
(quadratic) and P0 (constant) piecewise representations for electrical potential $u(x)$ and conductivity field $\sigma(x)$,
respectively, and solve systems of algebraic equations obtained after such discretization with {\tt UMFPACK}, a solver for
nonsymmetric sparse linear systems \cite{UMFPACK}. All computations for all models used in the current paper are performed
using 2D domain
\begin{equation}
  \Omega = \left\{ x \in \RR^2 : \ x_1^2 + x_2^2 < r^2_{\Omega} \right\},
  \label{eq:domain2D}
\end{equation}
which is a disc of radius $r_{\Omega} = 0.1$ with $m = 16$ equidistant electrodes $E_{\ell}$ with half-width $w = 0.12$ rad
covering approximately 61\% of boundary $\partial \Omega$ as shown in Figure~\ref{fig:geometry}(a).

Electrical potentials
\begin{equation}
  U = (U_{\ell})_{\ell=1}^{16} = \left\{ -3, +1, +2, -5, +4, -1, -3, +2, +4, +3, -3, +3, +2, -4, +1, -3 \right\}
  \label{eq:meas_a}
\end{equation}
are applied to electrodes $(E_{\ell})_{\ell=1}^{16}$ as seen in \eqref{eq:permutations} following the ``rotation scheme''
discussed in Section~\ref{sec:math_model}. These potentials are chosen to be consistent with the ground potential condition
\eqref{eq:curr_volt_conds}. Using PCA, tuned during the scale switching as discussed in Section~\ref{sec:pca_switch},
allows a relatively small number of principal components $N_{\xi}$ to operate on fine scales during the main course of the
optimization. As such, we will use only one permutation of the potentials $(U_{\ell})_{\ell=1}^{16}$ as shown in \eqref{eq:meas_a}
with the total number of measurements $m^2 = 256$, with $K=1$ in \eqref{eq:obj_fn_final}. Figures~\ref{fig:geometry}(b) and
\ref{fig:geometry}(c) show the examples of the distribution of flux $\sigma(x) \bnabla u(x)$ of electrical potential $u$ in
the interior of domain $\Omega$ and measured currents $(I_{\ell}^*)_{\ell=1}^{16}$ during EIT for two subsequent sets of
potentials
\begin{equation}
  \begin{aligned}
    U^1 = U &= \left\{ -3, +1, +2, -5, +4, -1, -3, +2, +4, +3, -3, +3, +2, -4, +1, -3 \right\},\\
    U^2 &= \left\{ -3, -3, +1, +2, -5, +4, -1, -3, +2, +4, +3, -3, +3, +2, -4, +1 \right\},
  \end{aligned}
  \label{eq:meas_b}
\end{equation}
respectively. To determine the Robin part of the boundary conditions in \eqref{eq:forward_3}, we equally set the electrode
contact impedance $Z_{\ell} = 0.1$.

Physical domain $\Omega$ is discretized using mesh totaling $N = 7726$ triangular FEM elements inside $\Omega$. This mesh is
then used to construct gradients $\bnabla_{\sigma} \mJ$, $\bnabla_{\xi} \mJ$, and $\bnabla_{\zeta} \mJ$, perform optimization
as described in Algorithm~\ref{alg:main_opt}, and compute maps $\mM^k$ using the partitioning methodology discussed in
Section~\ref{sec:coarse_scale}. We also apply the concept of spatial grouping provided in \cite{ChunMS2022} to determine (FEM)
elements in each subset $C_j$ (high-conductivity regions), $j = 2, \ldots, N_{\zeta}$, based on their location inside $\Omega$.
In particular, we define the ``neighboring'' principle when the neighbor elements go the $j$th subset at $k$th iteration if
they share at least one vertex and have conductivity $\sigma_i$ above the current threshold $\sigma^k_{th,j-1}$. For solving
optimization problems \eqref{eq:minJ_xi} and \eqref{eq:minJ_zeta}, our framework employs Sparse Nonlinear OPTimizer~{\tt SNOPT},
a software package for solving large-scale nonlinear optimization problems \cite{SNOPTManual}.

For all models used in this paper, the actual (true) electrical conductivity $\sigma_{true}(x)$ we seek to reconstruct is
given by
\begin{equation}
  \sigma_{true}(x) = \left\{
  \begin{aligned}
    \sigma_c, & \quad x \in \Omega_c,\\
    \sigma_h, & \quad x \in \Omega_h.
  \end{aligned}
  \right.
  \label{eq:sigma_true}
\end{equation}
In \eqref{eq:sigma_true}, $\sigma_c = 0.4$ and $\sigma_h = 0.2$ define cancer-affected regions of sub-domain $\Omega_c$
(spots of different sizes and complexity of their geometry depending on the model) and healthy tissue part $\Omega_h$,
respectively. The initial guess for control $\sigma(x)$ uses a constant approximation to $\sigma_{true}$ given by
$\sigma_0 = \frac{1}{2} \left(\sigma_h + \sigma_c \right) = 0.3$. Termination tolerances in \eqref{eq:termination} are set
to $\epsilon_c = 0$ (to avoid early termination at the coarse scale) and $\epsilon_f = 10^{-9}$. Optimization also terminates
after reaching the limit of 2,000 iterations or 100,000 objective function evaluations to assume the practical feasibility
of computational time. To enforce bounds established for coarse-scale control $\sigma^k_{th,n}$ in
\eqref{eq:sigma_coarse_bounds}, in all our computations, we used fine--to--coarse partition \eqref{eq:sigma_coarse}
redefined as
\begin{equation}
  \begin{aligned}
  \sigma^k_i = &\left\{
  \begin{aligned}
    \sigma^k_{low},  & \quad \sigma_i(\xi^k) < (1-\sigma^k_{th,n}) \,
    \underset{i}{\min} \ \sigma_i(\xi^k) + \sigma^k_{th,n} \,
    \underset{i}{\max} \ \sigma_i(\xi^k),\\
    \sigma^k_{high,n}, & \quad {\rm otherwise},
  \end{aligned}
  \right.\\
  &i = 1, \ldots N, \quad 1 \leq n \leq N_{\max}
  \end{aligned}
  \label{eq:sigma_coarse_new}
\end{equation}
while ensuring $0 < \sigma^k_{th,n} < 1$.

Finally, all computations in this paper use a PCA-based map \eqref{eq:map_pca} for discretized control $\sigma(x)$
established between $N$-dimensional $\sigma$-space and $N_{\xi}^k$-dimensional $\xi$-space, as described in
Sections~\ref{sec:fine_scale} and \ref{sec:pca_switch}. A set of $N_r = 1000$ realizations $(\sigma^*_i)_{i=1}^{1000}$
is generated using uniformly distributed random numbers. E.g., each realization $\sigma^*_i$ ``contains'' from one to
seven ``cancer-affected'' areas with $\sigma_c = 0.4$. Each area is located randomly within domain $\Omega$ and
represented by a circle of randomly chosen radius $0 < r \leq 0.3 r_{\Omega}$. We also apply TSVD by choosing the
initial number of principal components $N_{\xi}^0 = 662$ by retaining 662 basis vectors in the PCA description. This
value corresponds to the preservation of respectively 99\% of the ``energy'' in the full set of basis vectors; see
\cite{Bukshtynov2015,KoolmanBukshtynov2021} for more details.

\subsection{Model~\#1: Analysis of Gradient Validation}
\label{sec:kappa}

We start with model~\#1, used to demonstrate the applicability of the proposed computational framework discussed in
Section~\ref{sec:math} and check its overall performance while solving the inverse EIT problem. This model represents
a typical situation for a cancer-affected biological tissue containing several spots suspicious of tumor and, as such,
having elevated electrical conductivity. Model~\#1, featuring three circular-shaped cancerous regions of various sizes,
is shown in Figure~\ref{fig:geometry}(a).

\begin{figure}[!htb]
  \begin{center}
  \mbox{
  \subfigure[]
    {\includegraphics[width=0.5\textwidth]{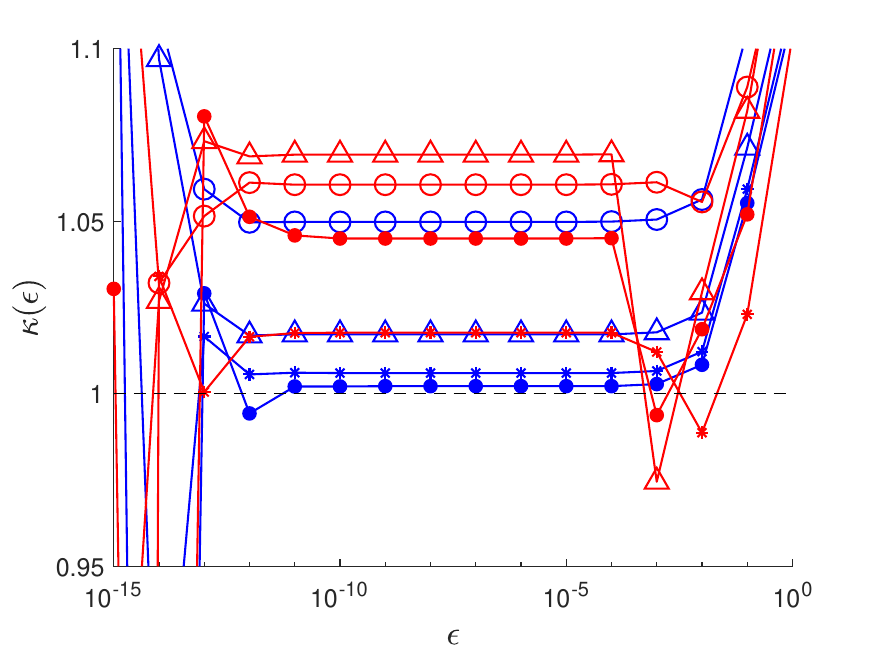}}
  \subfigure[]
    {\includegraphics[width=0.5\textwidth]{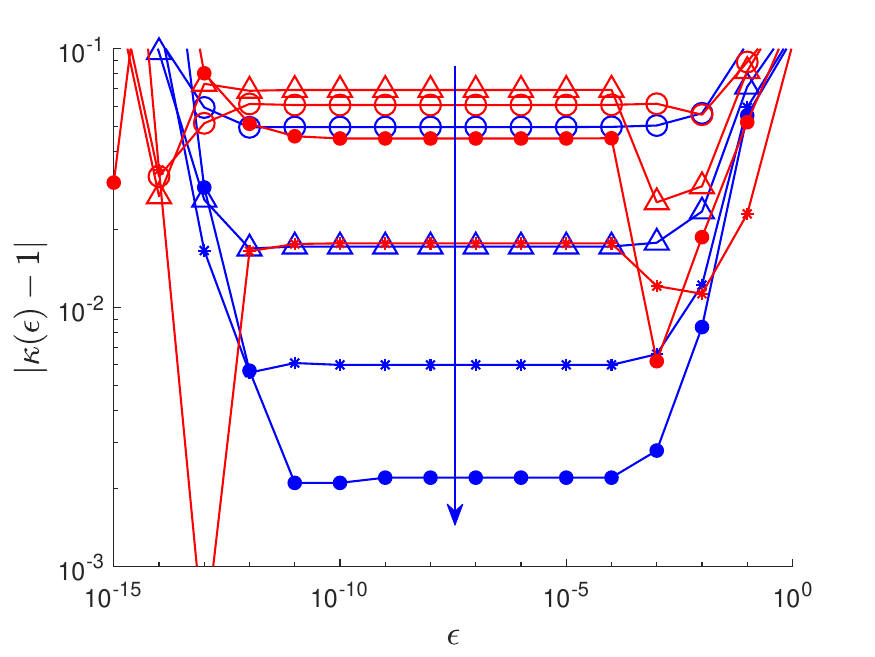}}}
  \end{center}
  \caption{The behavior of (a)~$\kappa (\epsilon)$ and (b)~$\log_{10} |\kappa(\epsilon) -1 |$ as a function of $\epsilon$
    while checking the consistency of gradients (blue)~$\bnabla_{\sigma} \mJ$ and (red)~$\bnabla_{\xi} \mJ$ computed
    for model~\#1 with different spatial discretization: (open circles)~$N =$ 712, (triangles)~$N =$ 2,032,
    (asterisks)~$N =$ 7,726, and (filled circles)~$N =$ 29,348. The blue arrow in (b) shows the direction in which the
    number of FEM elements increases.}
  \label{fig:kappa_cheap}
\end{figure}

First, we present results demonstrating the consistency of the fine-scale gradients before ($\bnabla_{\sigma} \mJ$) and
after ($\bnabla_{\xi} \mJ$) PCA-based parameterization. Figure~\ref{fig:kappa_cheap} shows the results of a diagnostic
test ($\kappa$-test) commonly employed to verify the correctness of the discretized gradients; see, e.g.,
\cite{Bukshtynov2011,Bukshtynov2013,BukshtynovBook2023}. It consists in computing the directional differential, e.g.,
$\mJ'(\xi; \delta \xi) = \left\langle \bnabla_{\xi} \mJ, \delta \xi \right\rangle_{L_2}$, for some selected variations
(perturbations) $\delta \xi$ in two different ways: namely, using a finite--difference approximation versus using
\eqref{eq:grad_sigma} with \eqref{eq:grad_xi} and then examining the ratio of the two quantities, i.e.,
\begin{equation}
  \kappa (\epsilon) = \dfrac{\frac{1}{\epsilon} \left[ \mJ(\xi + \epsilon \, \delta \xi) - \mJ(\xi) \right]}
  {\int_{\Omega} \bnabla_{\xi} \mJ \delta \xi \, dx}
  \label{eq:kappa}
\end{equation}
for a range of values of $\epsilon$. If these gradients are computed correctly, then for intermediate values of $\epsilon$,
$\kappa(\epsilon)$ will be close to the unity. Figure~\ref{fig:kappa_cheap}(a) demonstrates such behavior over a range of
$\epsilon$ spanning about 10 orders of magnitude for both $\bnabla_{\sigma} \mJ$ (in blue) and $\bnabla_{\xi} \mJ$ (in red).
As can be expected, the quantity $\kappa(\epsilon)$ deviates from the unity for very small values of $\epsilon$ due to the
subtractive cancelation (round--off) errors and also for large values of $\epsilon$ due to the truncation errors (both of
which are well--known effects). In addition, the quantity $\log_{10} | \kappa(\epsilon) - 1 |$ plotted in
Figure~\ref{fig:kappa_cheap}(b) shows how many significant digits of accuracy are captured in a given gradient evaluation.
Remarkably, this figure shows clear evidence that non-parameterized gradients $\bnabla_{\sigma} \mJ$ provide better
$\kappa$-test results (by moving values of $\kappa (\epsilon)$ closer to the unity) in case the spatial discretization of
domain $\Omega$ is refined. The blue arrow shows the direction in which the number of FEM elements increases. The reason is
that in the ``optimize--then--discretize'' paradigm adopted by the current computational framework, such refinement of
discretization leads to a better approximation of the continuous gradient \cite{ProtasBewleyHagen2004}. The same conclusion
does not apply to the gradients $\bnabla_{\xi} \mJ$ obtained with PCA-based transformation applied. For constructing this PCA,
all 1000 realizations $(\sigma^*_i)_{i=1}^{1000}$ were precomputed and stored using spatial discretization with $N =$ 7,726.
Therefore, refining spatial mesh has little effect on the gradients' quality: we will use the same number of FEM elements
(7,726) to construct mesh for all computations in the rest of this paper.

\begin{figure}[!htb]
  \begin{center}
  \mbox{
  \subfigure[$k = 0$]
    {\includegraphics[width=0.5\textwidth]{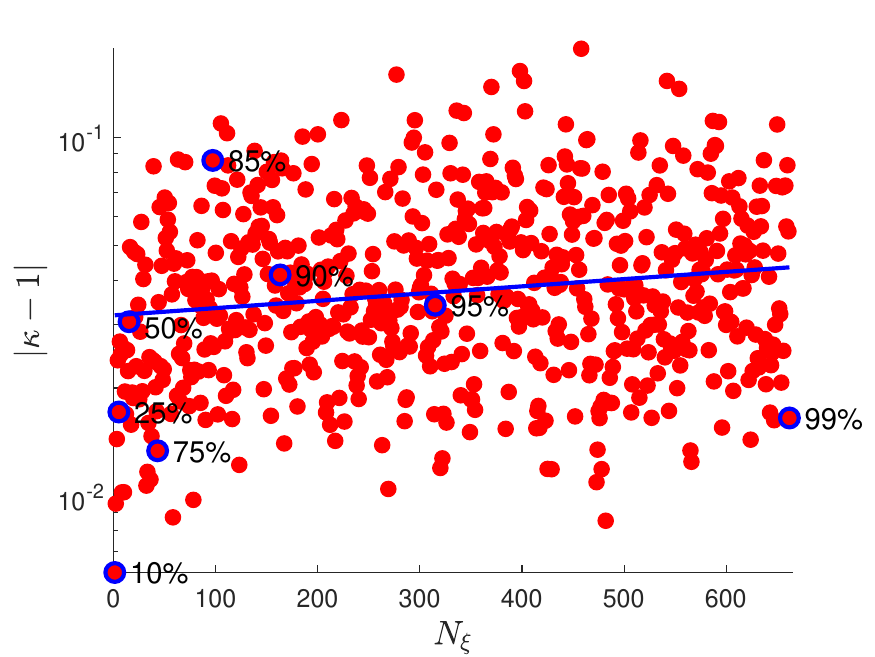}}
  \subfigure[$k = 200$]
    {\includegraphics[width=0.5\textwidth]{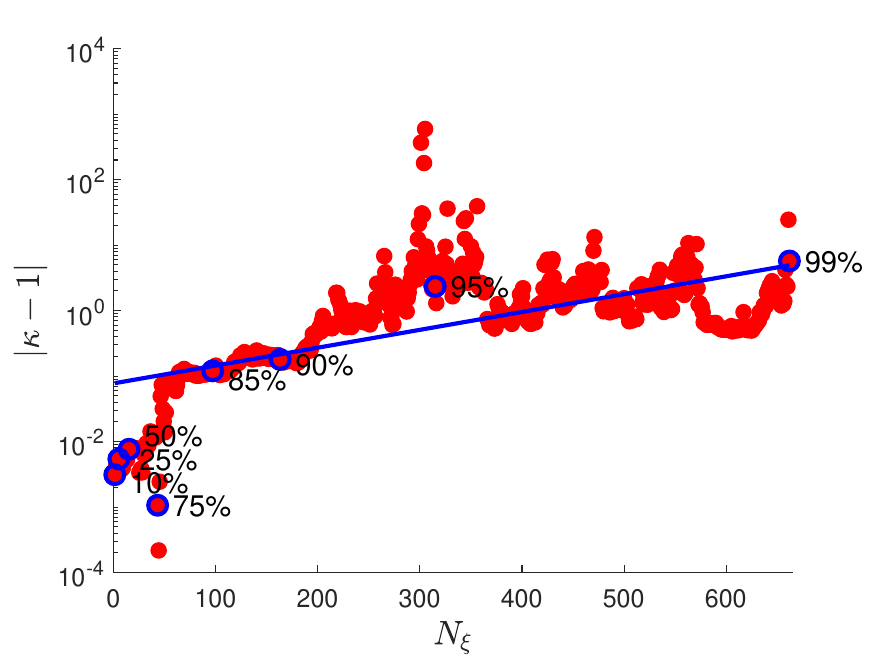}}}
  \end{center}
  \caption{The behavior of $\log_{10} |\kappa(\epsilon) -1 |$ as a function of $N_{\xi}$ for checking the consistency of
    gradients $\bnabla_{\xi} \mJ$ computed for model~\#1 with different sizes of the $\xi$-space ($N_{\xi} = 1, \ldots, 662$)
    when (a)~$k = 0$ and (b)~$k = 200$. For both plots, blue circles identify the results related to preserving a
    particular portion (in percent) of the ``energy'' in the full set of basis vectors used to construct the PCA transform.
    The blue lines describe the linear regression between the points.}
  \label{fig:kappa_cheap_xi}
\end{figure}

As we are particularly interested in applying PCA for parameterizing the control $\sigma(x)$ at the fine scale, the natural
question arises about the relation between the gradients' consistency and the number of principal components $N_{\xi}^k$ in
use. We repeat the $\kappa$-test multiple times for all values of $N_{\xi}$ (from 1 to 662), examining gradients $\bnabla_{\xi} \mJ$
at the beginning ($k = 0$) and after $k = 200$ optimization iterations; see Figures~\ref{fig:kappa_cheap_xi}(a) and
\ref{fig:kappa_cheap_xi}(b) for the respective outcomes. Although these plots depict different structures, both suggest the
same conclusion. The consistency (or accuracy, as we may suggest) of the reduced-dimensional gradients $\bnabla_{\xi} \mJ$
depends on $N_{\xi}$, and in general, it decreases when the size of the $\xi$-space increases.

\begin{figure}[!htb]
  \begin{center}
  \mbox{
  \subfigure[]
    {\includegraphics[width=0.5\textwidth]{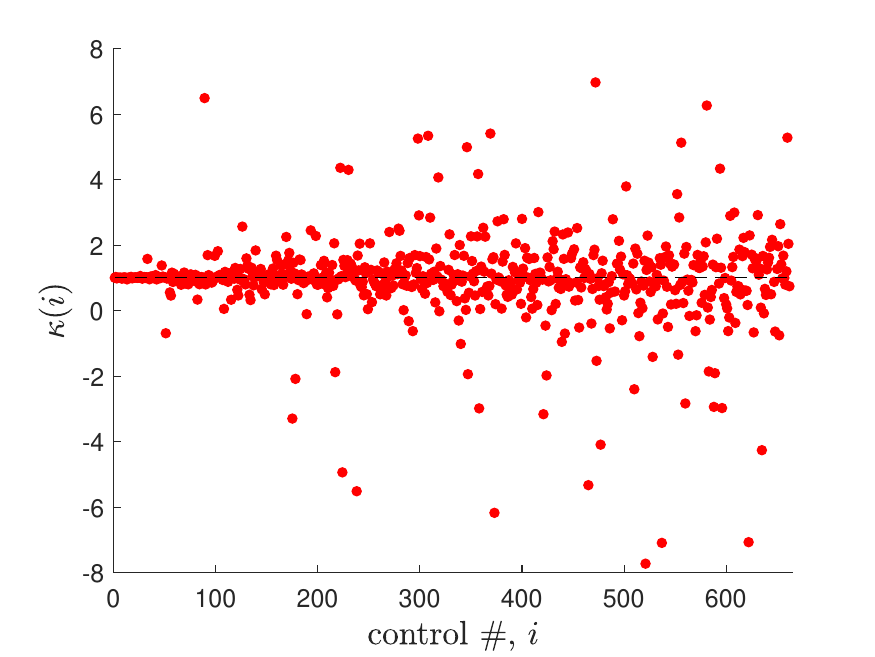}}
  \subfigure[]
    {\includegraphics[width=0.5\textwidth]{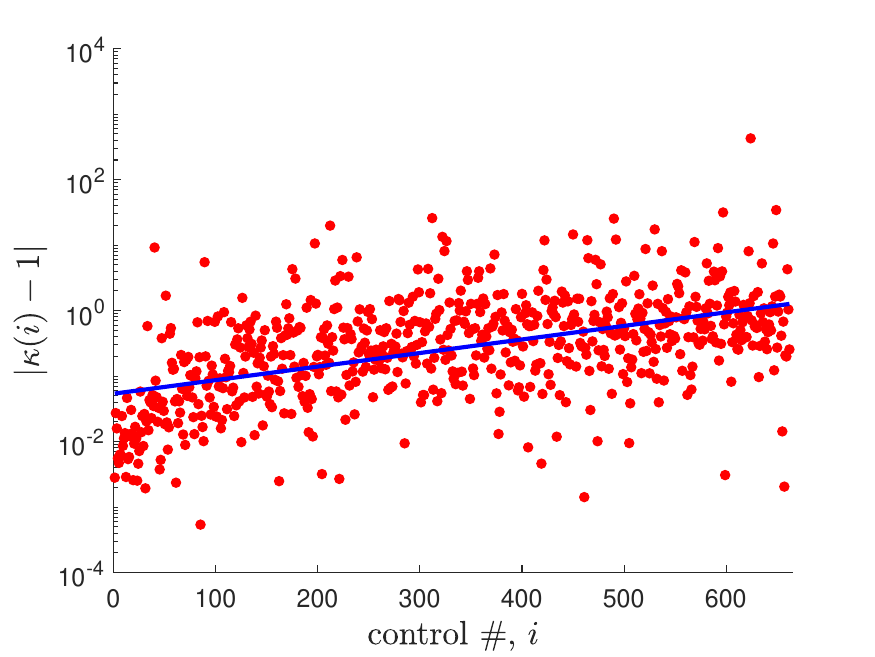}}}
  \end{center}
  \caption{The behavior of (a)~$\kappa (i)$ and (b)~$\log_{10} |\kappa(i) -1 |$ as a function of $i$, the number of
    the perturbed component $\xi_i$ in control vector $\xi$, while applying the ``expensive'' $\kappa$-test
    for $k = 0$. In (b), the blue line describes the linear regression between the points.}
  \label{fig:kappa_exp}
\end{figure}

We support this conclusion by applying the so-called ``expensive'' $\kappa$-test to check the sensitivity for all components
of the control vector $\xi$ by perturbing them individually using \eqref{eq:kappa} with fixed $\epsilon = 10^{-8}$; see
\cite{BukshtynovBook2023} for more details. In cases when the control vector components are associated with spatial locations,
this test helps locate elements with inaccurately computed gradient components due to the lack of sensitivity. The corresponding
points lay far outside the ``cloud'', which is positioned more or less symmetrically around 1. As shown in Figure~\ref{fig:kappa_exp},
this symmetry is not observed. Moreover, the points in Figure~\ref{fig:kappa_exp}(a) are more dispersed when the number $i$ of
the perturbed component $\xi_i$ in vector $\xi$ increases. It confirms that adding more ``energy'' to the PCA transform for
enabling the identification of small features in the physical $\sigma$-space does not necessarily lead to a better quality
of used gradients.

\begin{figure}[!htb]
  \begin{center}
  \mbox{
  \subfigure[$k = 0$]
    {\includegraphics[width=0.33\textwidth]{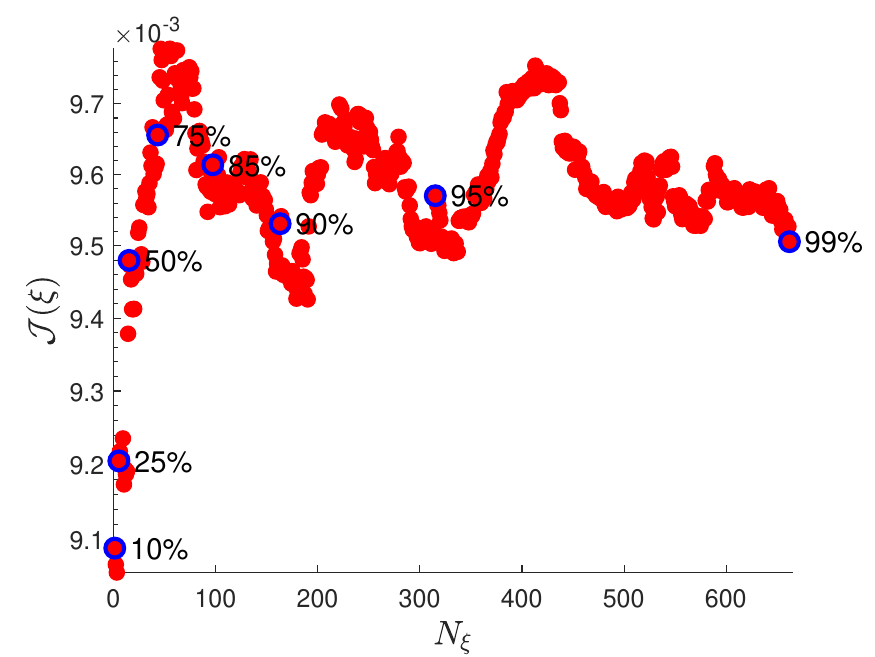}}
  \subfigure[$k = 200$]
    {\includegraphics[width=0.33\textwidth]{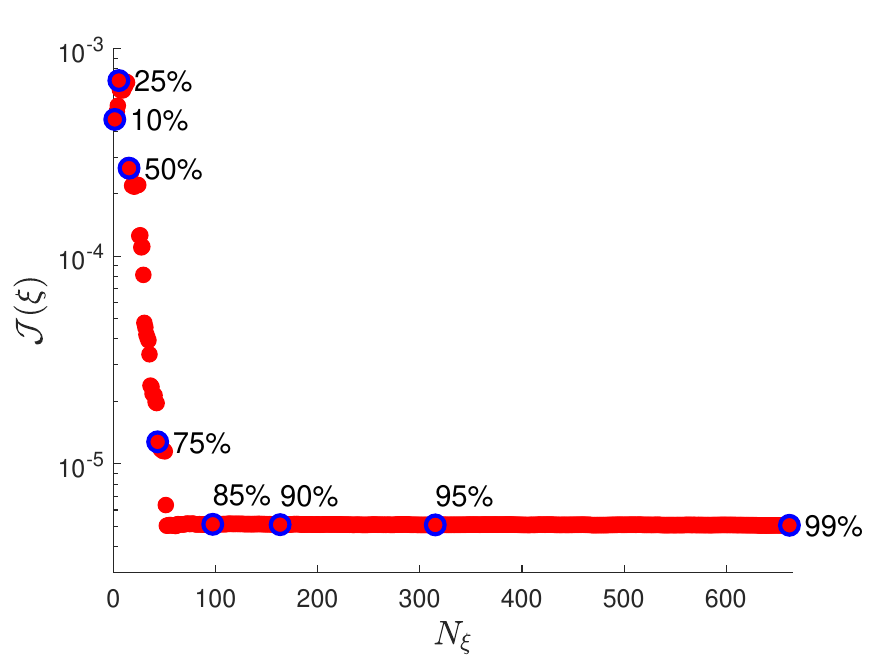}}
  \subfigure[$k = 200$ (close look)]
    {\includegraphics[width=0.33\textwidth]{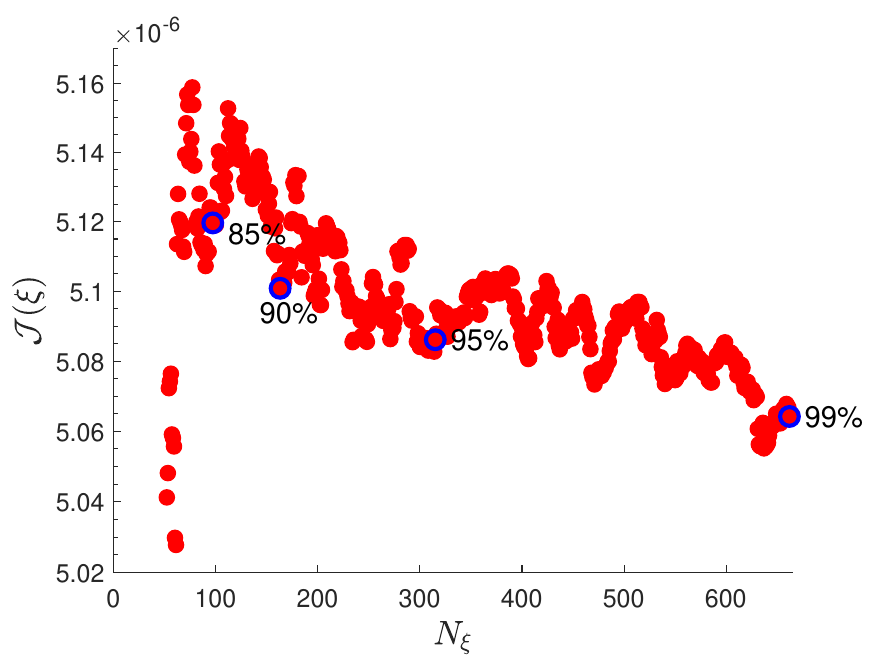}}}
  \end{center}
  \caption{Objective function $\mJ(\xi)$ evaluated for (a)~$k = 0$ and (b,c)~$k = 200$ as a function of principal
    components $N_{\xi}$. For all plots, blue circles identify the results related to preserving a particular portion
    (in percent) of the ``energy'' in the full set of basis vectors used to construct the PCA transform. (c)~A close
    look at the results in (b) that form a ``plateau'' for $k > 50$.}
  \label{fig:obj_xi}
\end{figure}

These results motivated us to explore the effect on the optimization results if we change the description of the
fine-scale space by tuning dynamically the number of principal components $N_{\xi}$ in the PCA-based parameterization
discussed in Section~\ref{sec:pca_switch}. We get even more convincing results after evaluating objective functions
$\mJ(\xi)$ for $k = 0$ and $k = 200$, see Figures~\ref{fig:obj_xi}(a) and \ref{fig:obj_xi}(b,c), respectively, while
using different values of $N_{\xi}$ to define the number of ``active'' components in the control vectors $\xi^0$ and
$\xi^{200}$. The comparison reveals that the ``optimal'' value of $N_{\xi}^k$ changes throughout the optimization
process depending on $k$: while it is close to 1 at the beginning ($k = 0$), it is between 50 and 60 after 200 iterations.
Unless stated otherwise, we apply the designed algorithms for enhanced scale switching by tuning PCA, as discussed
in Section~\ref{sec:pca_switch}, to all numerical experiments in this paper.

\subsection{Model~\#1: Main Computational Results}
\label{sec:m8_main}

In this section, we evaluate the performance of the proposed computational framework described in Section~\ref{sec:math}.
Here, we perform optimization using multiple scales supplied with multilevel parameterization as described in
Algorithm~\ref{alg:main_opt} in application to model~\#1. We refer to Section~\ref{sec:comp_model} for the main parameters
describing the model, spatial discretization, initiation and termination of the optimization process, and constructing the
PCA transform. In addition, for all computations, we switch scales after $n_s = 5$ iterations and use only one set of
measurements at the coarse scale (namely, $U^1$ in (39)), while fine-scale optimization employs all 16 sets, namely,
from $U^1$ to $U^{16}$.

\begin{figure}[!htb]
  \begin{center}
  \mbox{
  \subfigure[$\hat \sigma$: $N_{\max} = 1$]
    {\includegraphics[width=0.33\textwidth]{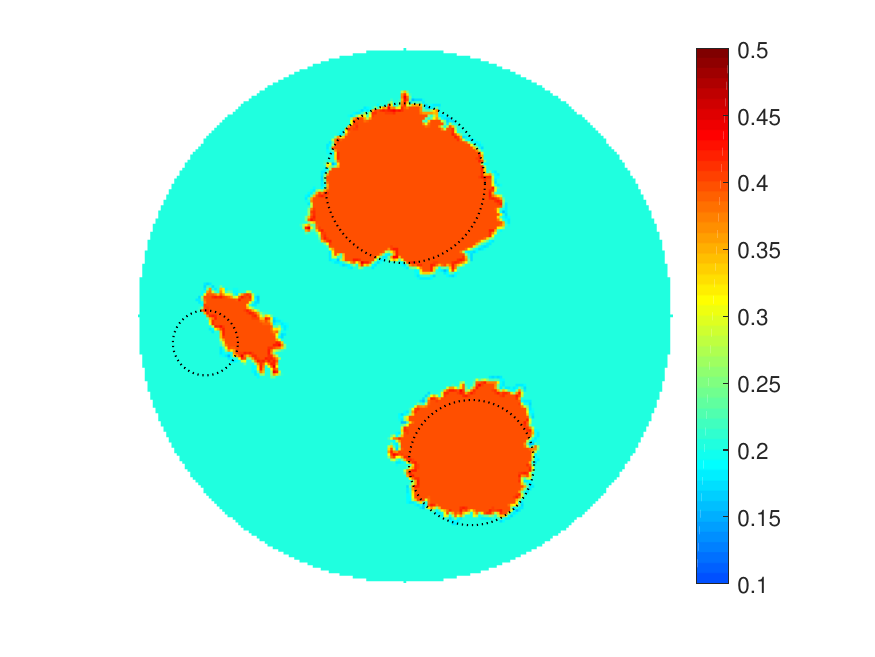}}
  \subfigure[objectives]
    {\includegraphics[width=0.33\textwidth]{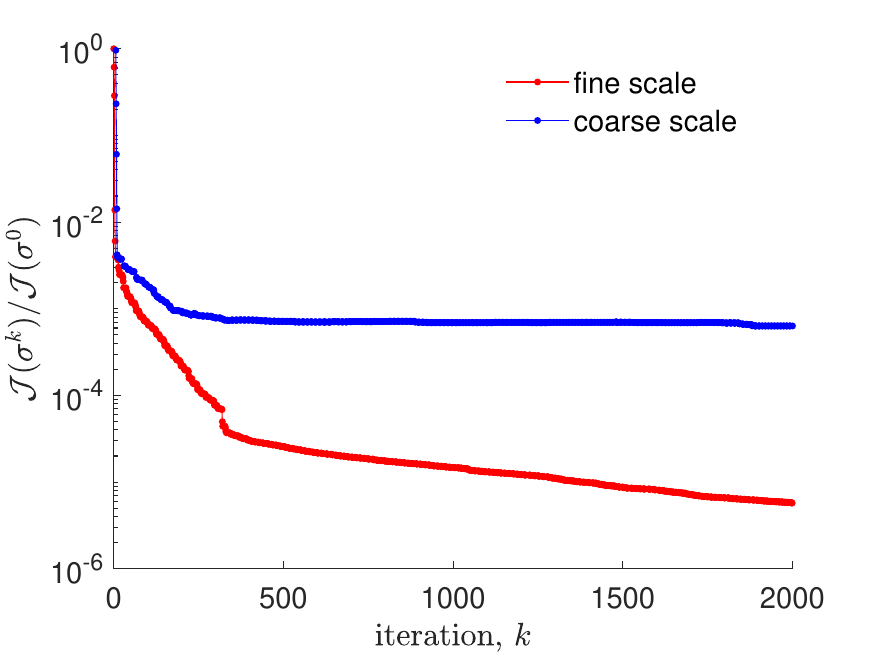}}
  \subfigure[coarse-scale controls]
    {\includegraphics[width=0.33\textwidth]{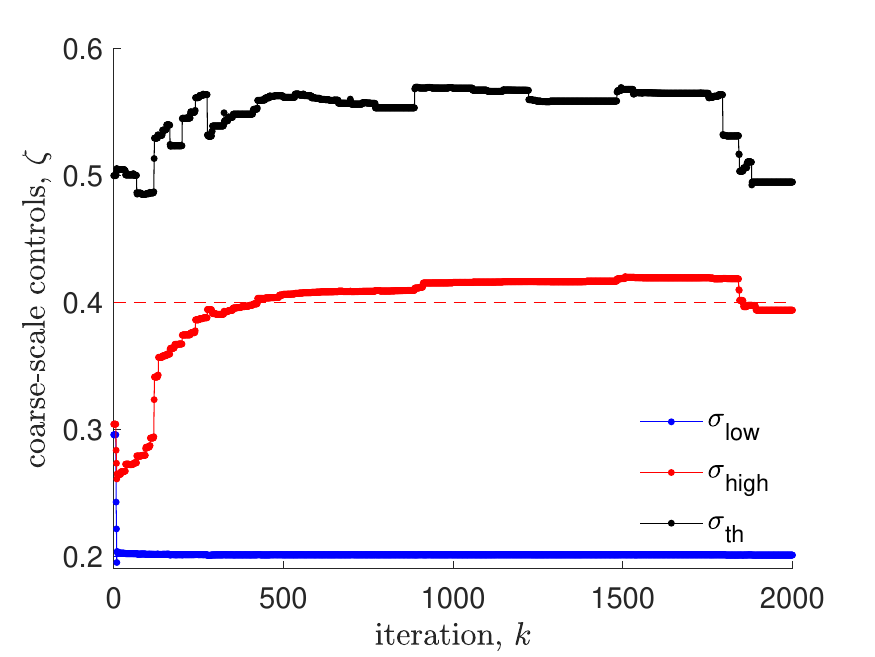}}}
  \mbox{
  \subfigure[$\hat \sigma$: $N_{\max} = 1$ \& PCA]
    {\includegraphics[width=0.33\textwidth]{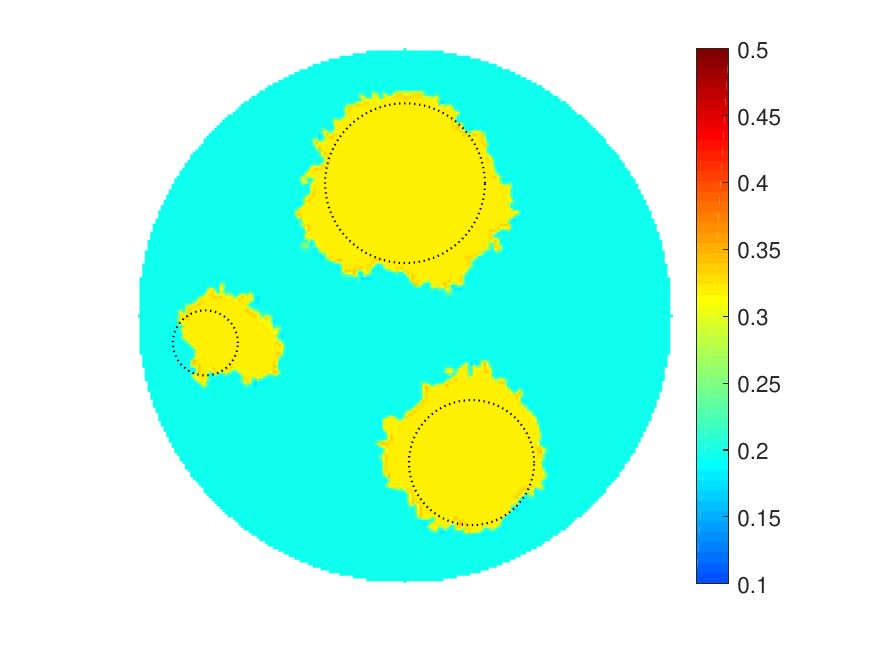}}
  \subfigure[objectives]
    {\includegraphics[width=0.33\textwidth]{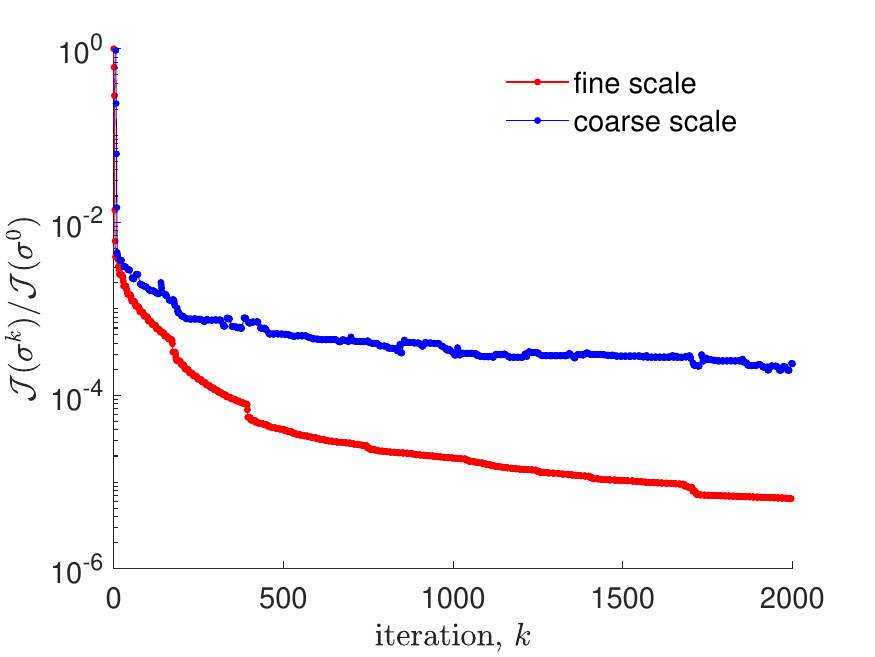}}
  \subfigure[coarse-scale controls]
    {\includegraphics[width=0.33\textwidth]{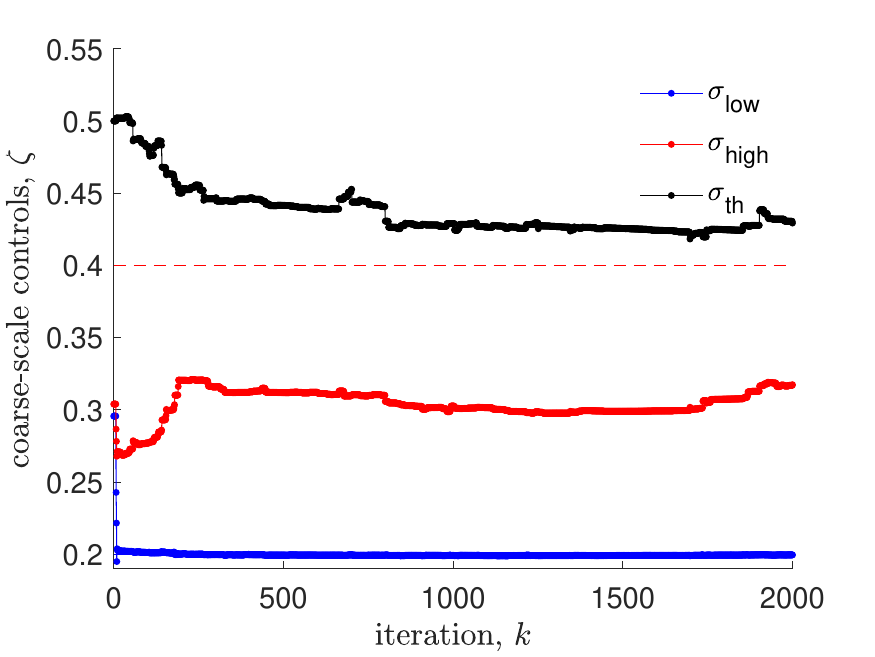}}}
  \end{center}
  \caption{Optimization results for model~\#1 (a-c)~without using and (d-f)~with the enhanced
    scale switching by the tuned PCA when $N_{\max} = 1$. Plots in (a,d)~show the images obtained
    after applying a multiscale framework by Algorithm~\ref{alg:main_opt} with added dashed circles
    to represent the location of cancer-affected regions taken from known $\sigma_{true}(x)$
    in Figure~\ref{fig:geometry}(a). Graphs in (b,e)~present normalized objective functions
    $\mJ(\sigma^k)/\mJ(\sigma^0)$ as functions of iteration count $k$ evaluated at fine (in red) and
    coarse (in blue) scales. Changes in the coarse-scale controls $\zeta^k = [\sigma^k_{low} \
    \sigma^k_{high} \ \sigma^k_{th}]$ are shown in (c,f) with $\sigma_c = 0.4$ (red dashed line),
    $\sigma^k_{low}$ in blue, $\sigma^k_{th}$ in black, and $\sigma^k_{high}$ in red.}
  \label{fig:m8_comp_a}
\end{figure}

Figure~\ref{fig:m8_comp_a} shows the first results obtained under the condition that all cancerous spots are treated as
one region with elevated electrical conductivity (i.e., the maximum number of expected cancerous spots $N_{\max} = 1$)
without using and with the enhanced scale switching by the tuned PCA. Comparison of the respective optimal solutions
$\hat \sigma$ in Figures~\ref{fig:m8_comp_a}(a) versus \ref{fig:m8_comp_a}(d) reveals the results that differ in
performance using different metrics. E.g., the shapes of two big spots in both cases appear reasonably accurate. However,
the case using switching assisted by the tuned PCA underestimated the high-conductivity coarse-scale control;
$\hat \sigma_{high} = 0.394$ vs.~$\hat \sigma_{high} = 0.317$ as the last values on plots (red curves) of
Figures~\ref{fig:m8_comp_a}(c) and \ref{fig:m8_comp_a}(f), respectively. At the same time, this case resulted in more
accurate shape recognition for the smallest spot. As an additional characteristic to describe the performance, we could
also assess the gap between the objectives evaluated at fine and coarse scales, see Figures~\ref{fig:m8_comp_a}(b) and
\ref{fig:m8_comp_a}(e), as a measure of effectiveness in ``communication'' between the scales. Here, the use of PCA
contributes positively as this gap tends to decrease throughout the entire optimization.

\begin{figure}[!htb]
  \begin{center}
  \mbox{
  \subfigure[$\hat \sigma$: $N_{\max} = 5$]
    {\includegraphics[width=0.33\textwidth]{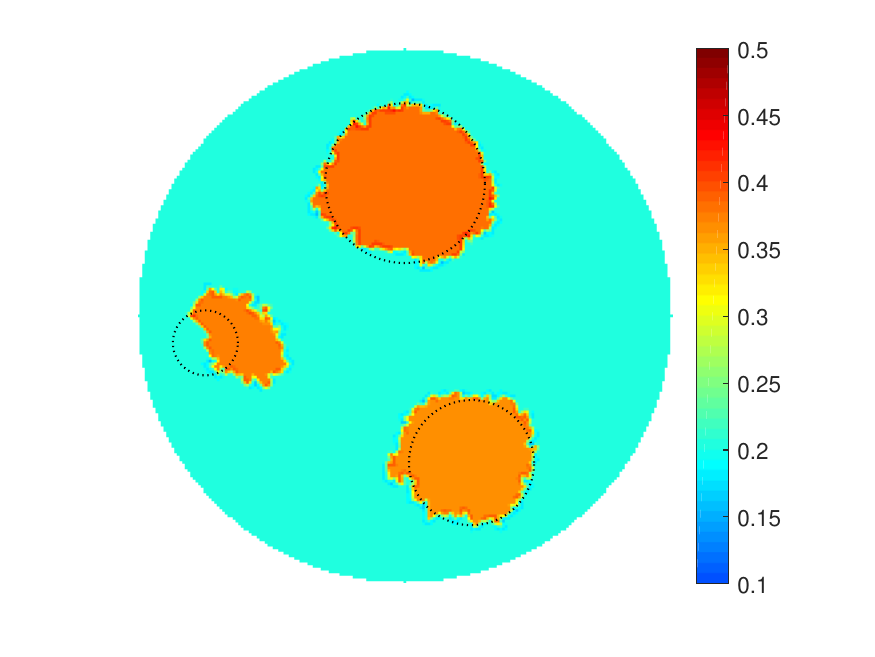}}
  \subfigure[objectives]
    {\includegraphics[width=0.33\textwidth]{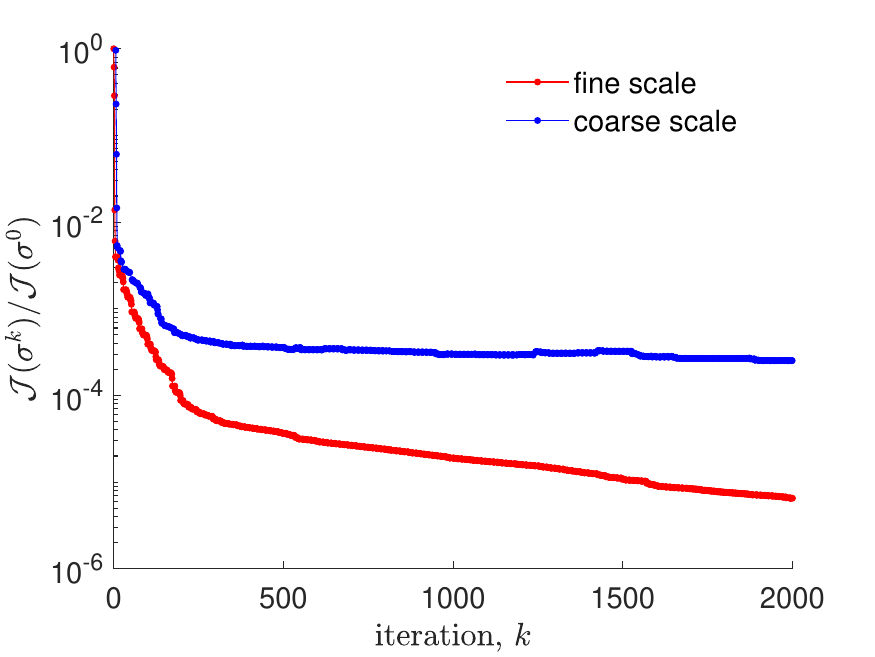}}
  \subfigure[coarse-scale controls]
    {\includegraphics[width=0.33\textwidth]{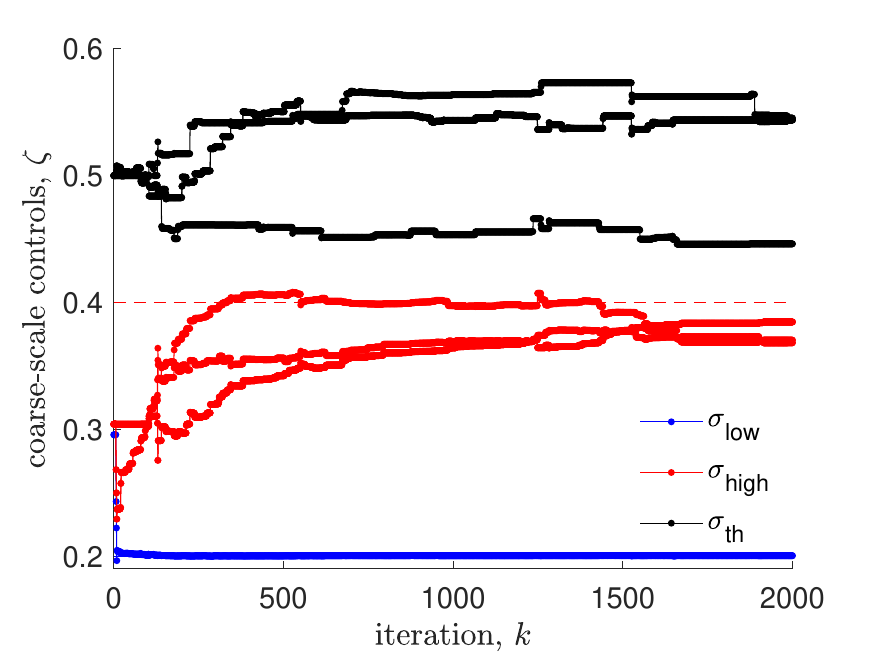}}}
  \mbox{
  \subfigure[$\hat \sigma$: $N_{\max} = 5$ \& PCA]
    {\includegraphics[width=0.33\textwidth]{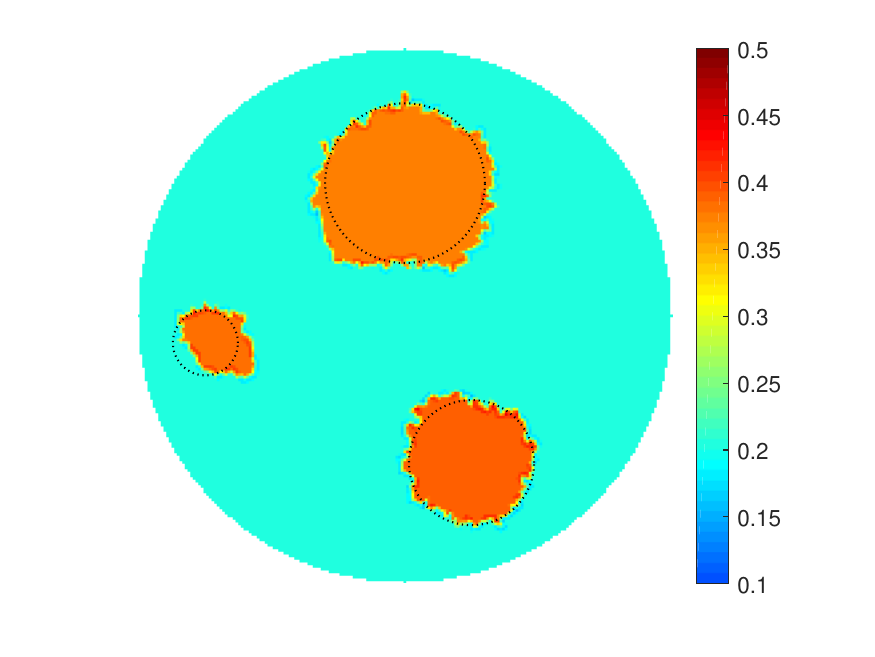}}
  \subfigure[objectives]
    {\includegraphics[width=0.33\textwidth]{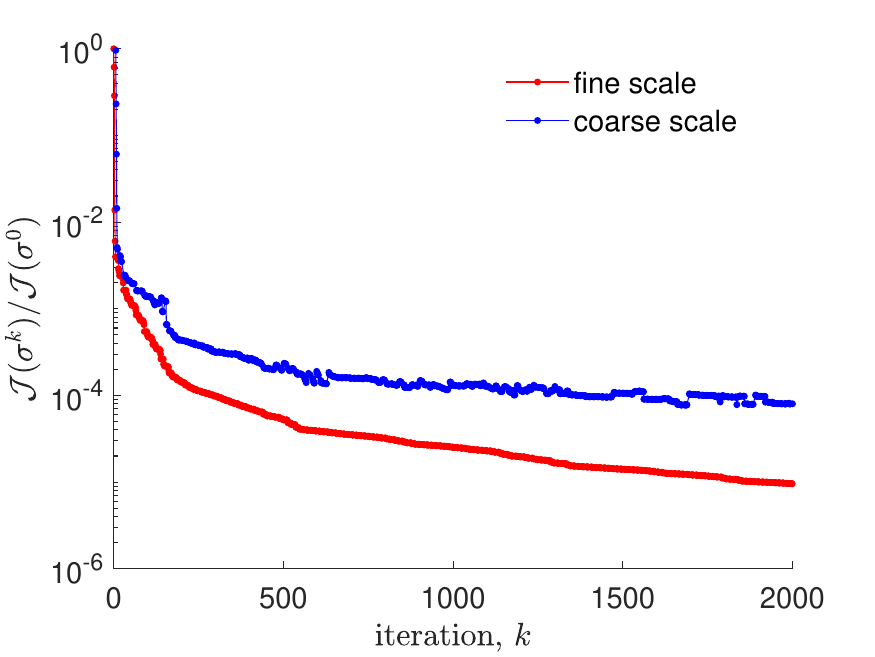}}
  \subfigure[coarse-scale controls]
    {\includegraphics[width=0.33\textwidth]{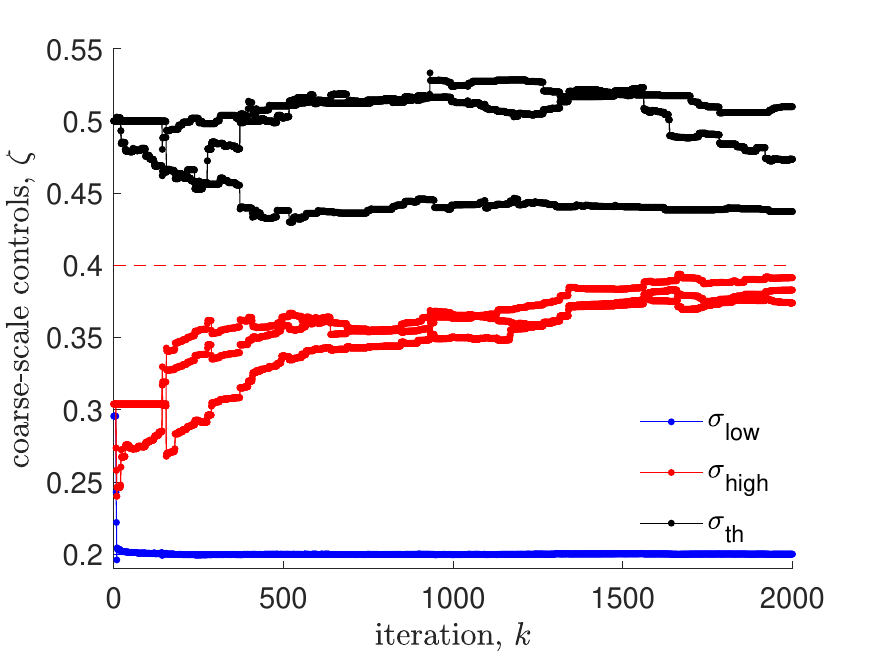}}}
  \end{center}
  \caption{Optimization results for model~\#1 (a-c)~without using and (d-f)~with the enhanced
    scale switching by the tuned PCA when $N_{\max} = 5$. Plots in (a,d)~show the images obtained
    after applying a multiscale framework by Algorithm~\ref{alg:main_opt} with added dashed circles
    to represent the location of cancer-affected regions taken from known $\sigma_{true}(x)$
    in Figure~\ref{fig:geometry}(a). Graphs in (b,e)~present normalized objective functions
    $\mJ(\sigma^k)/\mJ(\sigma^0)$ as functions of iteration count $k$ evaluated at fine (in red) and
    coarse (in blue) scales. Changes in the coarse-scale controls $\zeta^k = [\sigma^k_{low} \
    \sigma^k_{high,n} \ \sigma^k_{th,n}]$ ($n = 1, \ldots, 3$) are shown in (c,f) with $\sigma_c = 0.4$
    (red dashed line), $\sigma^k_{low}$ in blue, $\sigma^k_{th,n}$ in black, and $\sigma^k_{high,n}$ in
    red.}
  \label{fig:m8_comp_b}
\end{figure}

We explain the deficiency in the results obtained with $N_{\max} = 1$ as the method's inability to adequately
interpret reconstructed shapes based on the information enclosed in the fine-scale images. Including all spots
into a single area with only one control $\sigma_{high}$ assigned to represent high conductivity is compensated
either by defective shapes (i.e., separation thresholds $\sigma^k_{th}$) or values of $\sigma^k_{high}$. After
increasing the number of expected cancerous spots, say to $N_{\max} = 5$ (here and later, we assume that we do
not know a priori an exact number of individual regions), these defects naturally disappear. Although switching
from $N_{\max} = 1$ to $N_{\max} = 5$ alone shows progress in reconstructing accurate shapes, added PCA-assisted
scale switching makes the reconstruction almost perfect in both shapes and $\sigma_{high}$ values; refer to images
and the history of the coarse-scale controls in Figures~\ref{fig:m8_comp_b}(a,c) and \ref{fig:m8_comp_b}(d,f),
respectively. Figures~\ref{fig:m8_comp_b}(b) and \ref{fig:m8_comp_b}(e) also show continuing progress in exchanging
information between the scales: the case with $N_{\max} = 5$ and tunable PCA has evidently the best performance,
confirmed by the minimal gap between objectives at both scales.

\begin{figure}[!htb]
  \begin{center}
  \mbox{
  \subfigure[$\hat \sigma$: fine scale only]
    {\includegraphics[width=0.33\textwidth]{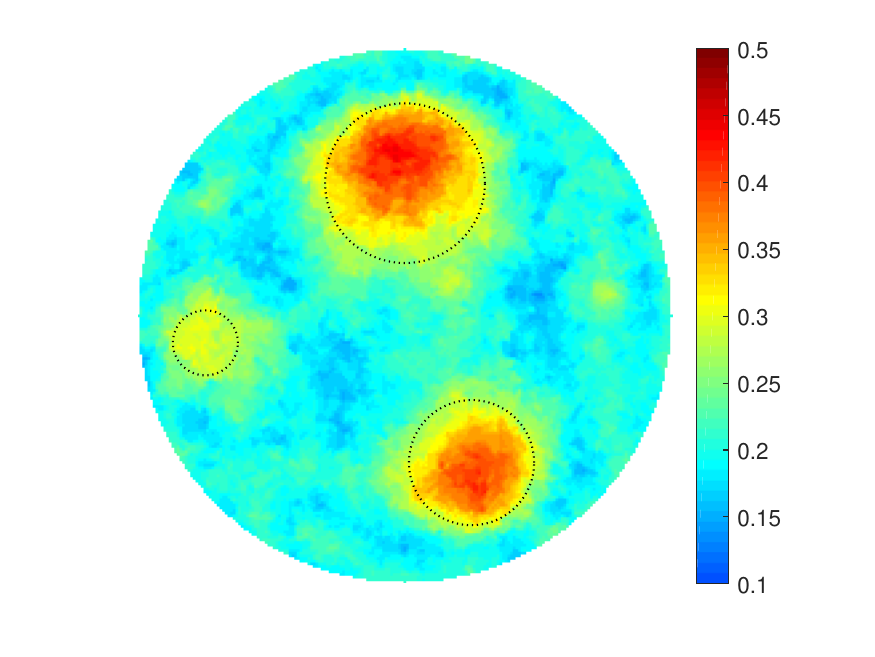}}
  \subfigure[solution error]
    {\includegraphics[width=0.33\textwidth]{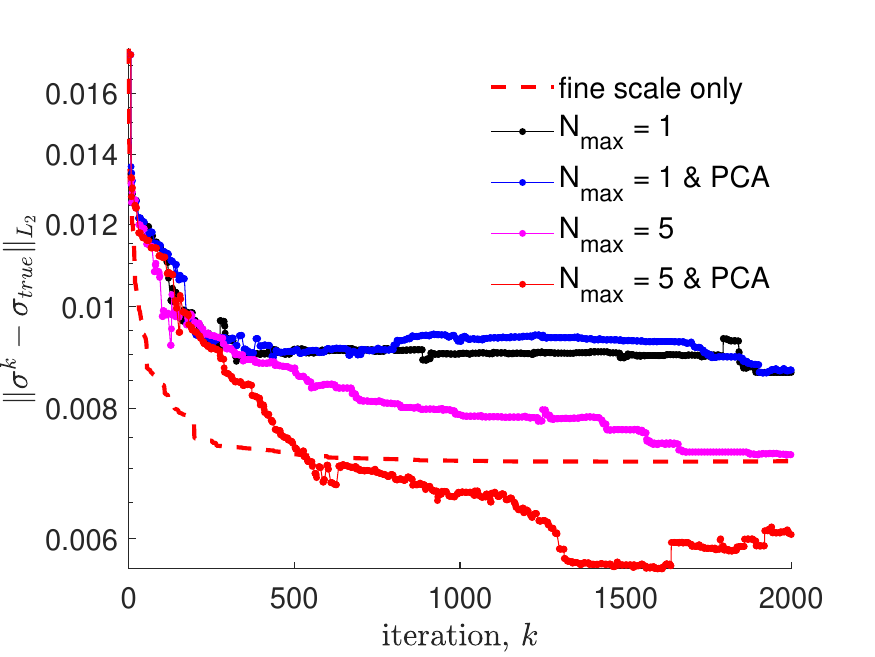}}
  \subfigure[PCA tuning]
    {\includegraphics[width=0.33\textwidth]{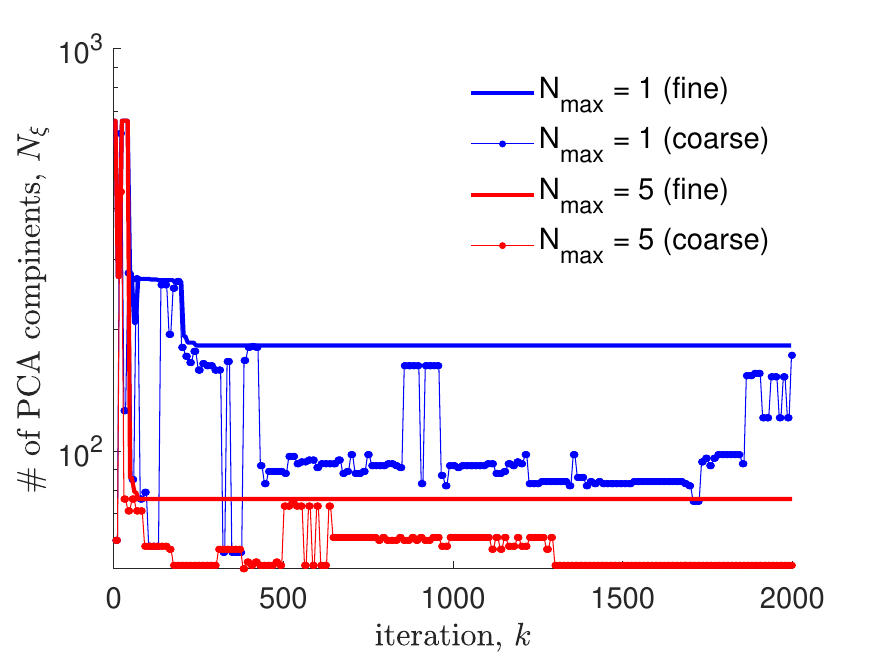}}}
  \end{center}
  \caption{(a)~The image of model~\#1 obtained after optimizing only at the fine scale. The dashed circles
    are added to represent the location of cancer-affected regions taken from known $\sigma_{true}(x)$ in
    Figure~\ref{fig:geometry}(a). (b)~Solution errors $\| \sigma^k - \sigma_{true} \|_{L_2}$ as functions of
    iteration count $k$ evaluated at the coarse scale for cases $N_{\max} = 1$ and $N_{\max} = 5$
    with/without the use of the tunable PCAs. The four cases are compared with the result obtained
    at the fine scale (dashed line). (c)~Numbers of principal PCA components as functions of iteration count $k$
    for cases $N_{\max} = 1$ and $N_{\max} = 5$ updated dynamically for fine and coarse scales.}
  \label{fig:m8_comp_c}
\end{figure}

To complete the performance analysis, we refer to Figure~\ref{fig:m8_comp_c}(a), which shows the results of
reconstructing all three cancerous spots by optimization performed only at the fine scale. The appearance of
these spots leaves no doubts about the complications mentioned above related to identifying the correct location
and boundaries between healthy and cancer-affected areas. However, due to the high resolution of such images,
they still provide a better match for the available measurements (data); refer to Figures~\ref{fig:m8_comp_a}(b,e)
and \ref{fig:m8_comp_b}(b,e). As a posteriori assessment of the quality of the reconstructed images, we use the
$L_2$-norm error $\| \sigma^k - \sigma_{true} \|_{L_2}$ evaluated only at the coarse scale.
Figure~\ref{fig:m8_comp_c}(b) demonstrates that $N_{\max} = 1$ cases with and without tunable PCA have about the
same error at the end of optimization. Case with $N_{\max} = 5$ has a much better result compared to the error
of the fine-scale reconstruction, and $N_{\max} = 5$ with added PCA-assisted scale switching has minimal error
much below the threshold established by the fine-scale solution. Finally, we refer to Figure~\ref{fig:m8_comp_c}(c),
showing the history of tuning PCA for both $N_{\max} = 1$ and $N_{\max} = 5$. We notice that the case with an
underestimated number of expected cancerous spots ($N_{\max} = 1$) tends to use higher numbers of principal
components to compensate for the lack of information at the coarse scale and the smaller size of its control space.
It is also evident that the fine scales tune their PCAs actively at the beginning when the high-resolution images
have large-scale changes in their structures. The coarse scales, however, show sensitivity to the updated PCAs
throughout the major part of the optimization runs.

\begin{figure}[!htb]
  \begin{center}
  \mbox{
  \subfigure[$\hat \sigma$:  $\alpha_{\max} = 0.95$]
    {\includegraphics[width=0.33\textwidth]{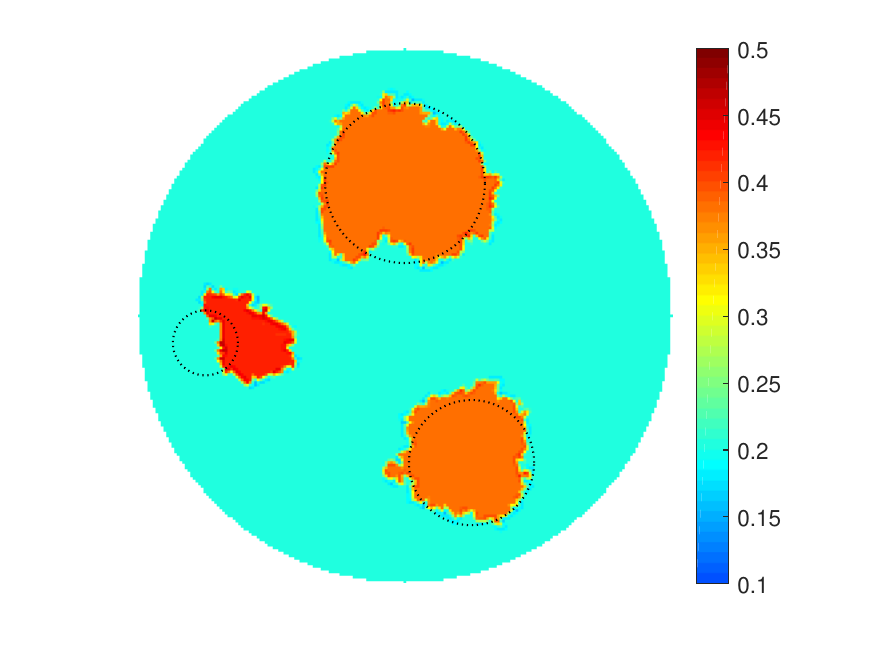}}
  \subfigure[$\hat \sigma$:  $\alpha_{\max} = 0.9$]
    {\includegraphics[width=0.33\textwidth]{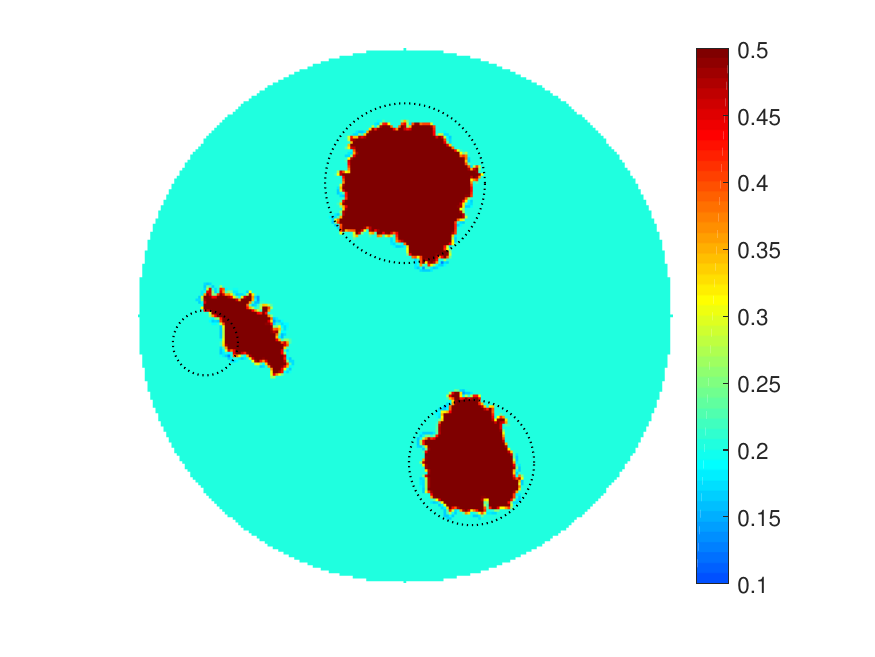}}
  \subfigure[solution error]
    {\includegraphics[width=0.33\textwidth]{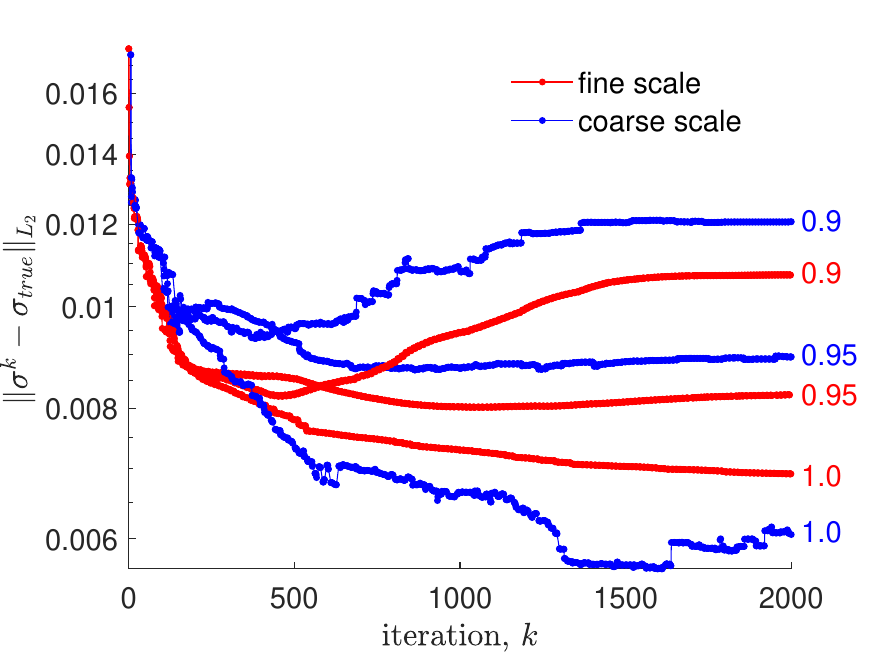}}}
  \end{center}
  \caption{Optimization results for model~\#1 obtained with (a)~$\alpha_{\max} = 0.95$ and (b)~$\alpha_{\max} = 0.9$
    with the enhanced scale switching by the tuned PCA when $N_{\max} = 5$. The dashed circles are added to represent
    the location of cancer-affected regions taken from known $\sigma_{true}(x)$ in Figure~\ref{fig:geometry}(a).
    (c)~Solution errors $\| \sigma^k - \sigma_{true} \|_{L_2}$ as functions of iteration count $k$ evaluated at
    coarse (in blue) and fine (in red) scales for cases $\alpha_{\max} = 1.0$, 0.95, and 0.9.}
  \label{fig:m8_thresh}
\end{figure}

In the next turn, we would like to address some issues related to setting parameter $\alpha_{\max}$ to a proper
value while solving the optimization problem \eqref{eq:proj_CtoF_rlx} for finding an optimal value of relaxation
parameter $\alpha_{c \rightarrow f}$ in the coarse-to-fine projection. Setting $\alpha_{\max}$ to 1 accepts that
in some iterations, solutions obtained at the coarse scale may not contribute to the next-phase fine-scale solutions.
However, the user may set this parameter to a lower value (i.e.,  $0 \leq \alpha_{\max} < 1$) to forcedly ensure
communication between scales for every coarse-to-fine switching by mixing at least $1 - \alpha_{\max}$ part of the
coarse-scale solution with the fine-scale one. Indeed, it depends on the problem, and any improper interference in
the dynamics of such communication may lead to worsened performance. We illustrate this fact in Figure~\ref{fig:m8_thresh},
where (a) and (b) plots depict the images obtained with $\alpha_{\max} = 0.95$ and $\alpha_{\max} = 0.9$, respectively.
In this particular problem, ``enforced communication'' resulted in overshooting, affecting solutions at the fine scale,
and as an implicit consequence, at the coarse scale. Figure~\ref{fig:m8_thresh}(c) compares the solution errors
evaluated at both scales for $\alpha_{\max} = 1.0$, 0.95, and 0.9 to exemplify this phenomenon in applications to
our current problem. Although setting $\alpha_{\max}$ to values between 1 and 0.95 seems to provide minimal harm,
we will keep $\alpha_{\max} = 1$ for the rest numerical experiments.

\begin{figure}[!htb]
  \begin{center}
  \mbox{
  \subfigure[model~\#1 (modified)]
    {\includegraphics[width=0.33\textwidth]{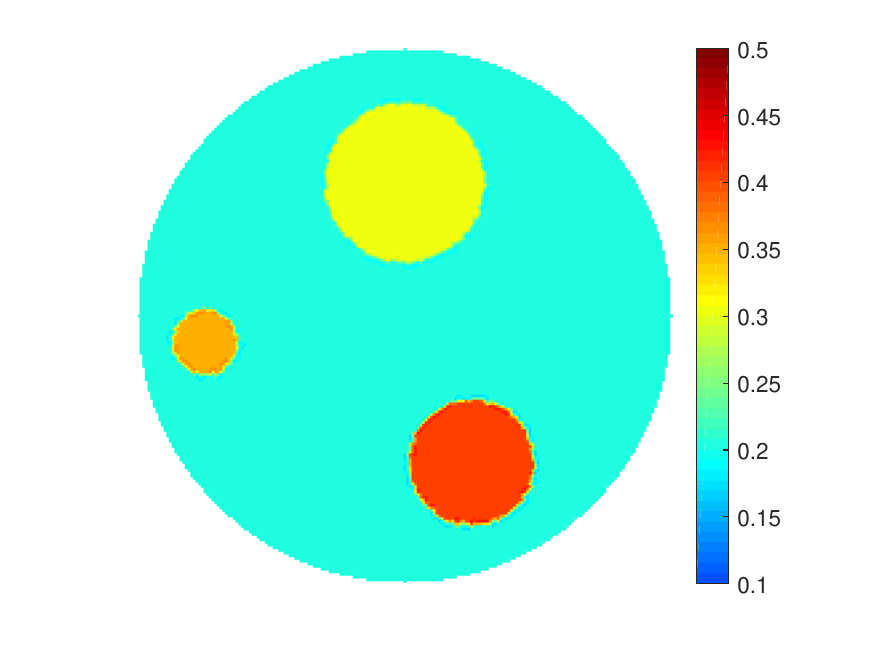}}
  \subfigure[$\hat \sigma$: fine scale only]
    {\includegraphics[width=0.33\textwidth]{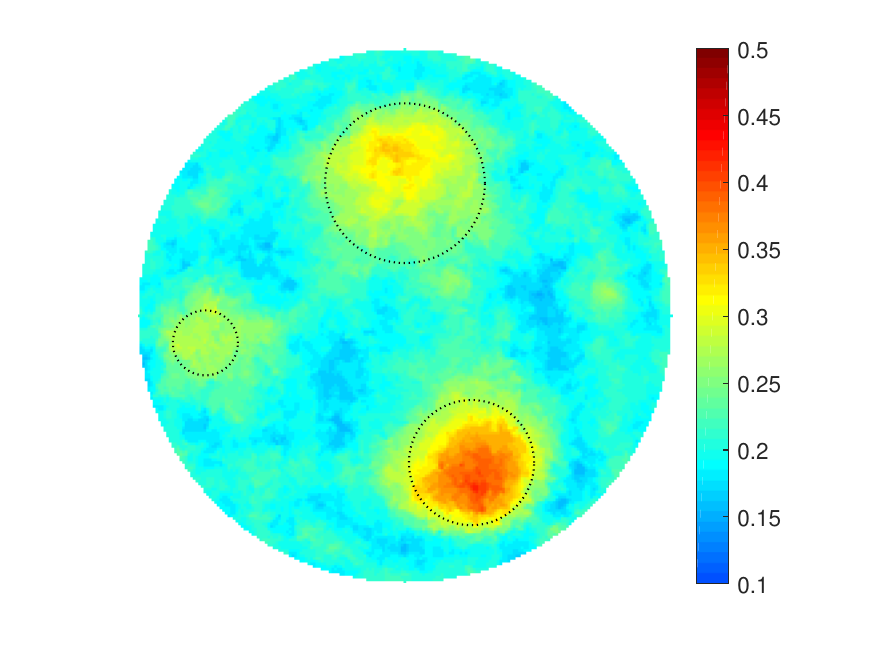}}
  \subfigure[solution error]
    {\includegraphics[width=0.33\textwidth]{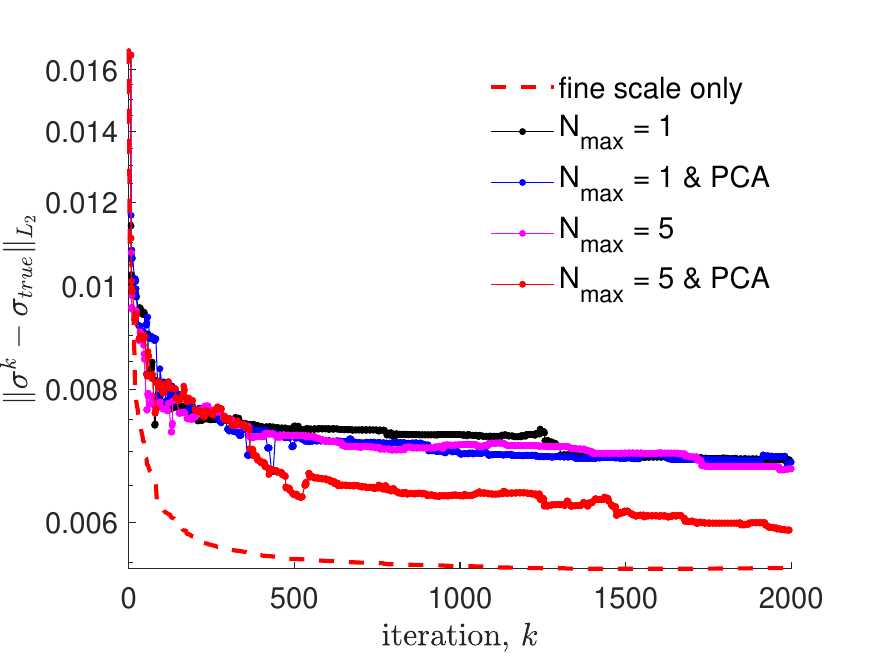}}}
  \mbox{
  \subfigure[$\hat \sigma$: $N_{\max} = 1$]
    {\includegraphics[width=0.25\textwidth]{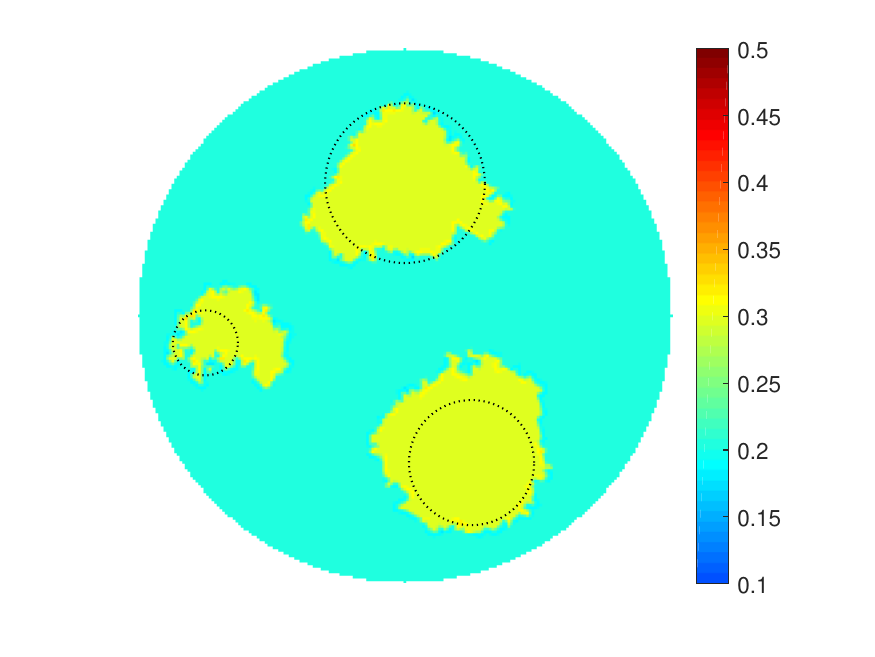}}
  \subfigure[$\hat \sigma$: $N_{\max} = 1$ \& PCA]
    {\includegraphics[width=0.25\textwidth]{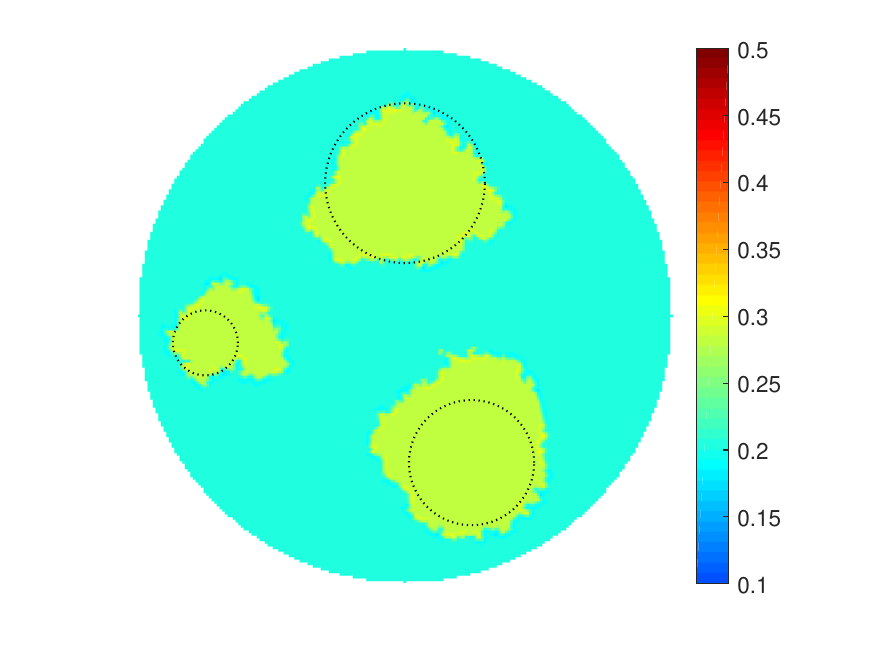}}
  \subfigure[$\hat \sigma$: $N_{\max} = 5$]
    {\includegraphics[width=0.25\textwidth]{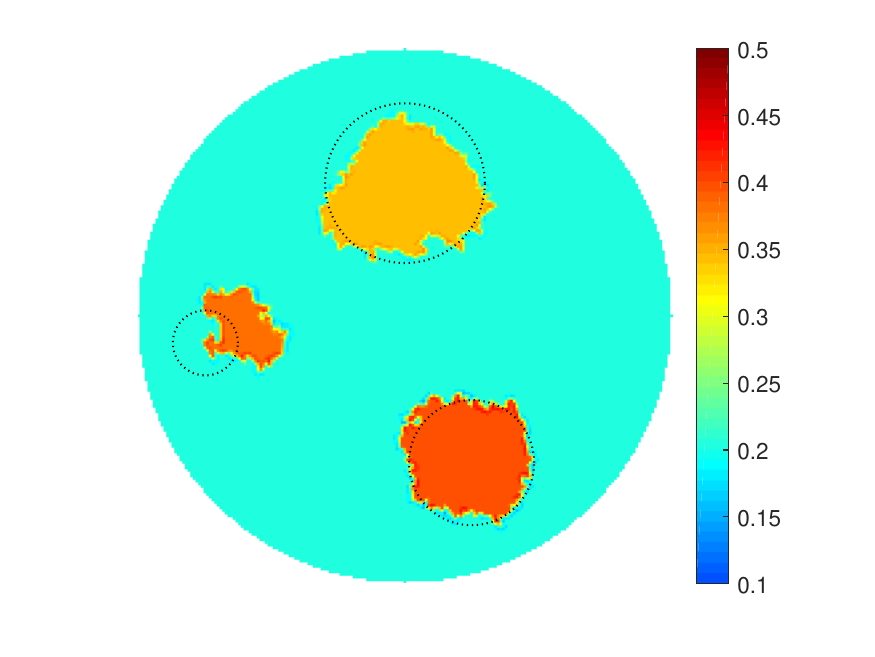}}
  \subfigure[$\hat \sigma$: $N_{\max} = 5$ \& PCA]
    {\includegraphics[width=0.25\textwidth]{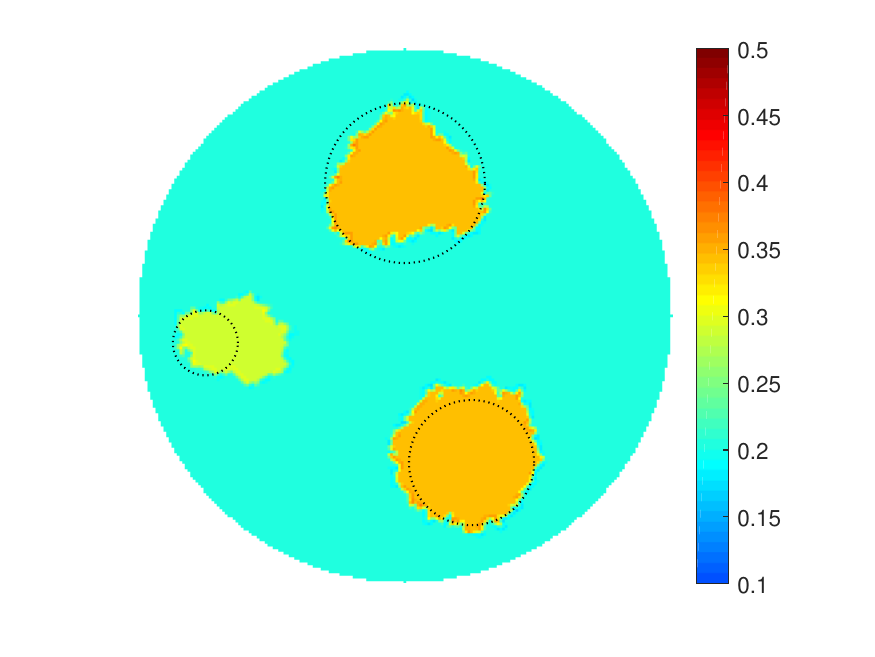}}}
  \end{center}
  \caption{(a)~EIT model~\#1 (modified): true electrical conductivity $\sigma_{true}(x)$. (b)~The image of modified
    model~\#1 obtained after optimizing only at the fine scale. (c)~Solution errors $\| \sigma^k - \sigma_{true} \|_{L_2}$
    as functions of iteration count $k$ evaluated at the coarse scale for cases $N_{\max} = 1$ and $N_{\max} = 5$
    with/without the use of the tunable PCAs. The four cases are compared with the result obtained at the fine scale
    (dashed line). (d-g)~Optimization results for the modified model~\#1 (d,f)~without using and (e,g)~with the
    enhanced scale switching by the tuned PCA when (d,e)~$N_{\max} = 1$ and (f,g)~$N_{\max} = 5$. Plots in (b,d-g) show
    the images obtained after applying a multiscale framework by Algorithm~\ref{alg:main_opt} with added dashed circles
    to represent the location of cancer-affected regions taken from known $\sigma_{true}(x)$ in (a).}
  \label{fig:m8_diff_cols}
\end{figure}

To further evaluate the performance of the proposed multiscale optimization framework, we modified our model~\#1 by
changing the high values of the electrical conductivity inside the cancerous spots while keeping the same their mutual
positioning and sizes. Figure~\ref{fig:m8_diff_cols}(a) displays the modified model~\#1, where we set $\sigma_c$ to 0.3,
0.4, and 0.35 for the big, medium-size, and small spots, respectively. As seen in Figure~\ref{fig:m8_diff_cols}(b),
the fine-scale-only image provides a limited ability to identify the boundaries for regions of small sizes or if
$\sigma_c$ does not deviate too much from $\sigma_h = 0.2$ (a big spot with $\sigma_c = 0.3$). We ran optimization
four times to compare the images obtained when $N_{\max} = 1$ and 5 and with/without tunable PCA for each case; refer
to Figures~\ref{fig:m8_diff_cols}(d-g). As mentioned before, failure to reconstruct correctly high conductivities is
natural when all spots are treated as a single area, $N_{\max} = 1$ in (d) and (e) plots, with only one control $\sigma_{high}$
assigned. In addition, the situation with the modified model~\#1 is aggravated by the fact that all three spots have
different conductivity. Increasing $N_{\max}$ to 5 in Figure~\ref{fig:m8_diff_cols}(f) changes the outcome significantly
-- all three $\sigma_{high}$ controls are fairly accurate: $\hat \sigma_{high,1} = 0.339$ (compared to the true value of 0.3),
$\hat \sigma_{high,2} = 0.399$ (true value of 0.4), and $\hat \sigma_{high,3} = 0.378$ (true value of 0.35). This result
is good, especially considering the accuracy in reconstructing shapes of big spots. We also notice that our 1000 sample
solutions (realizations) $(\sigma^*_n)_{n=1}^{1000}$ used to construct PCA transformation do not contain spots with
variable conductivities. Adding tuned PCA in Figure~\ref{fig:m8_diff_cols}(g) worsens the image for final ``colors''
related to different conductivities; however, it improves the accuracy in boundary positioning. Figure~\ref{fig:m8_diff_cols}(c)
summarizes the performance conclusions by comparing the solution errors for all four cases and identifying the last image
($N_{\max} = 5$ \& PCA) as the best solution for the complicated problem set by our modified problem~\#1.

\subsection{Effect of Noisy Data}
\label{sec:noise}

In this section, we address a well-known issue of noise that might be present in the measurements due to improper
electrode-medium contacts, possible electrode misplacement, wire interference, etc. The effect of noise influencing
the solutions of the inverse problems has already been investigated by many researchers both theoretically and numerically
to mitigate its negative impact on the quality of the obtained images. Here, we compare the effect of noise in
reconstructions obtained by applying the proposed multiscale optimization framework described in Chapter~\ref{sec:math} and
Algorithm~\ref{alg:main_opt} with the maximum number of expected cancerous spots $N_{\max} = 5$ and with tunable PCA-based
scale switching.

\begin{figure}[!htb]
  \begin{center}
  \mbox{
  \subfigure[$\hat \sigma$: fine only, $0.25\%$]
    {\includegraphics[width=0.25\textwidth]{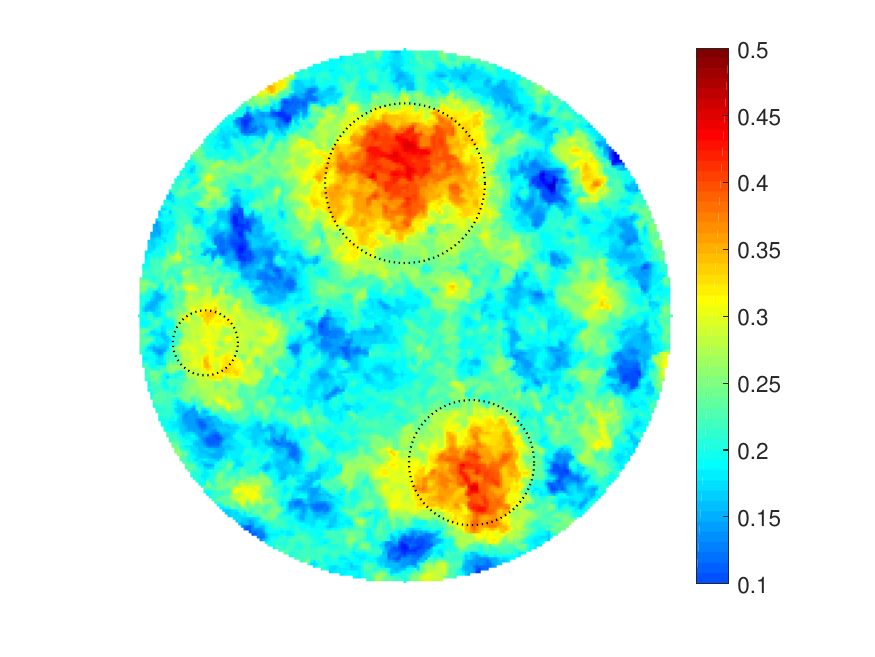}}
  \subfigure[$\hat \sigma$: fine only, $0.5\%$]
    {\includegraphics[width=0.25\textwidth]{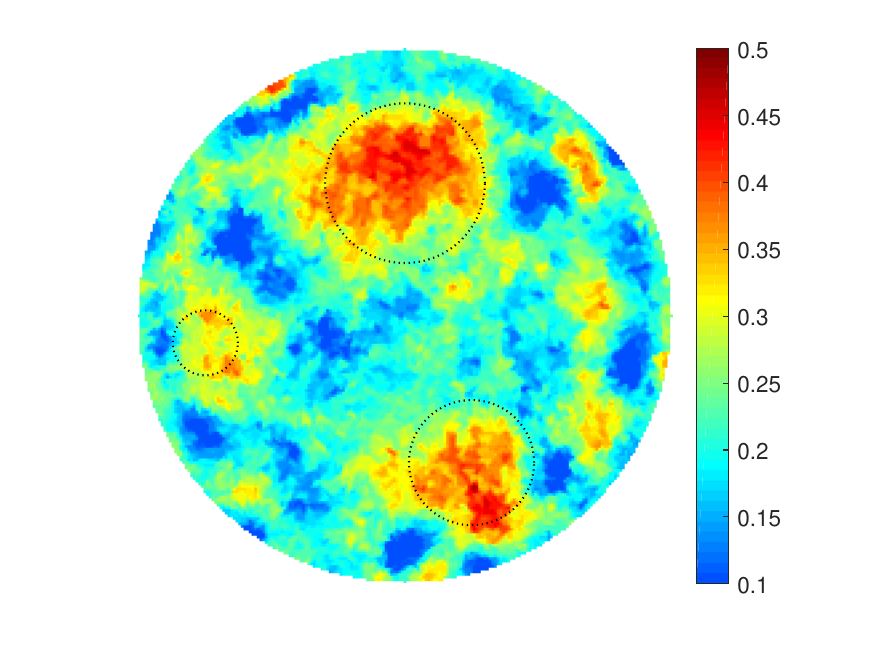}}
  \subfigure[$\hat \sigma$: fine only, $1.0\%$]
    {\includegraphics[width=0.25\textwidth]{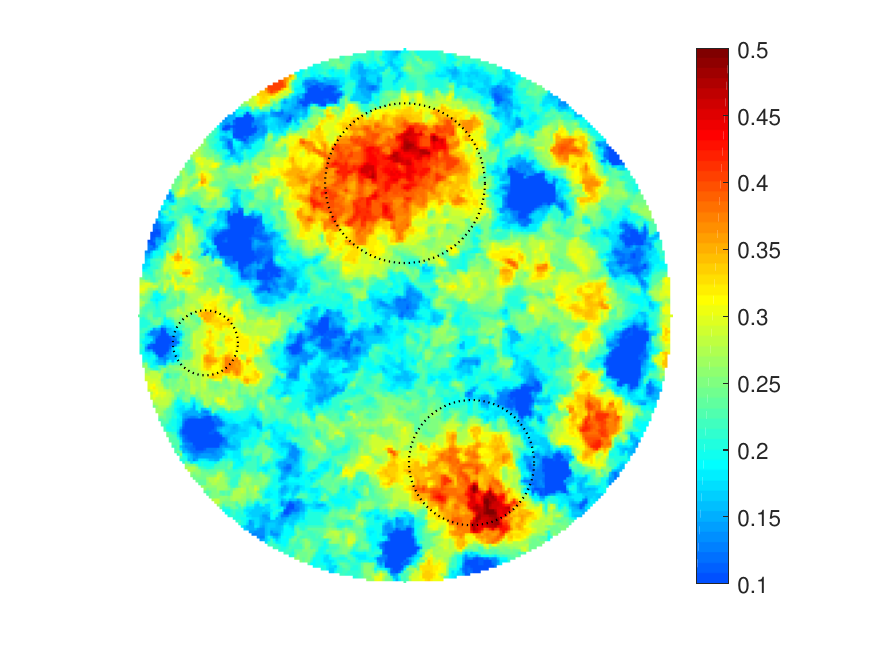}}
  \subfigure[$\hat \sigma$: fine only, $2.0\%$]
    {\includegraphics[width=0.25\textwidth]{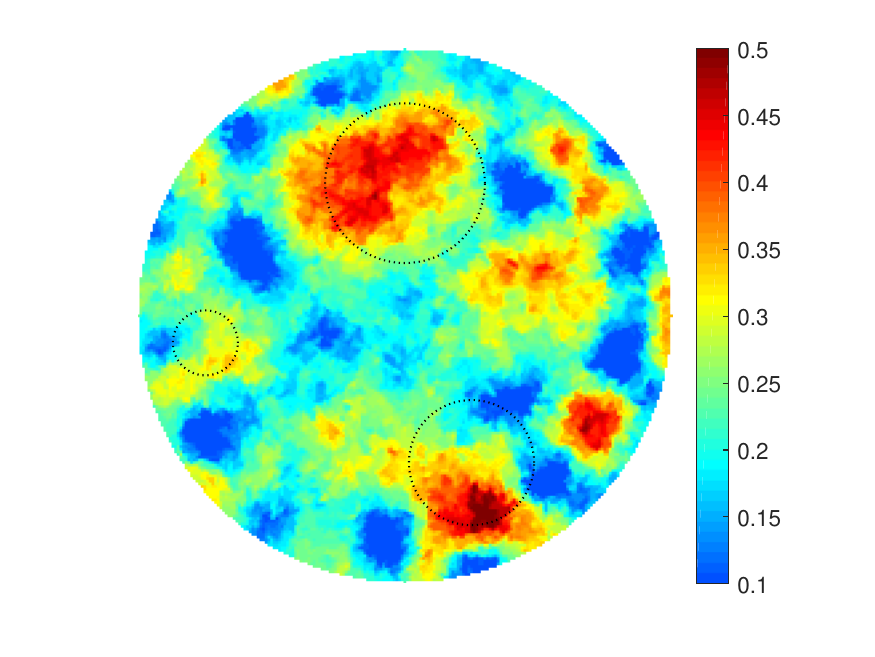}}}  \mbox{
  \subfigure[$\hat \sigma$: coarse scale, $0.25\%$]
    {\includegraphics[width=0.25\textwidth]{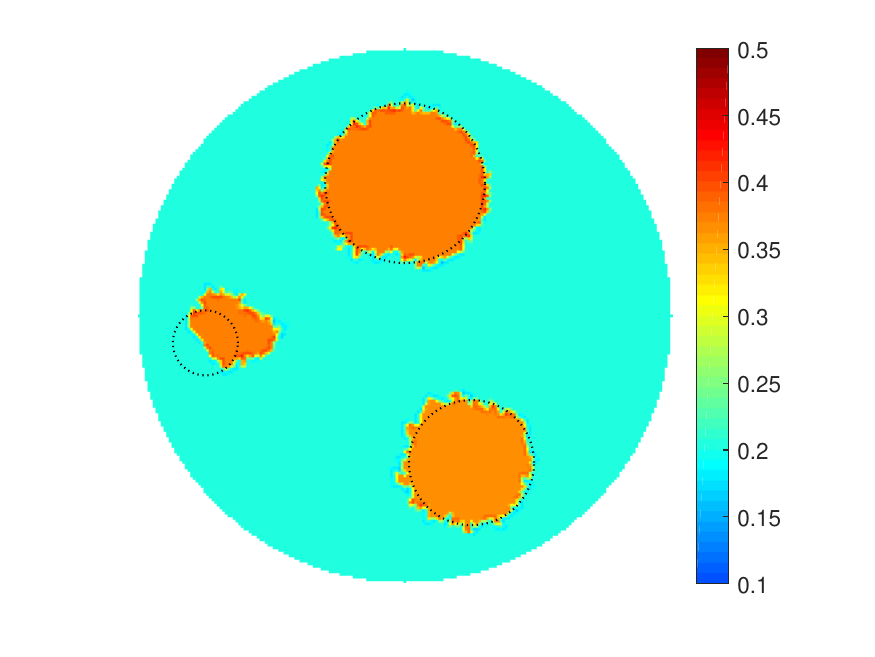}}
  \subfigure[$\hat \sigma$: coarse scale, $0.5\%$]
    {\includegraphics[width=0.25\textwidth]{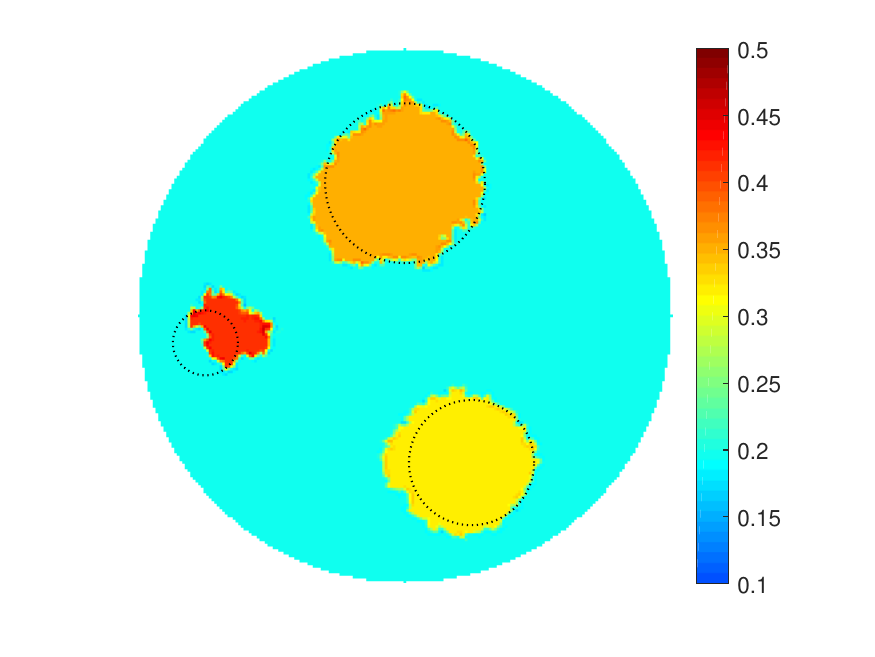}}
  \subfigure[$\hat \sigma$: coarse scale, $1.0\%$]
    {\includegraphics[width=0.25\textwidth]{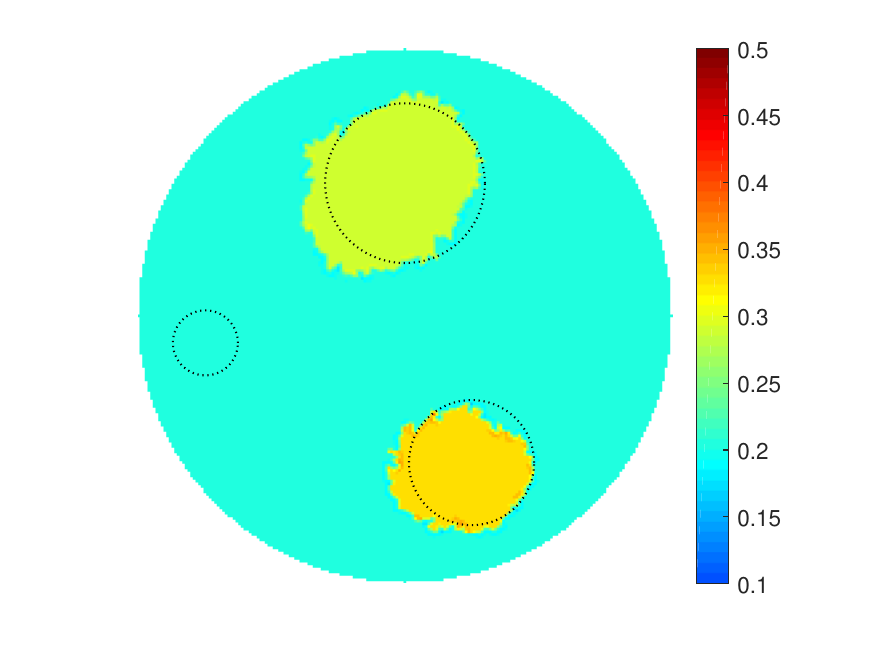}}
  \subfigure[$\hat \sigma$: coarse scale, $2.0\%$]
    {\includegraphics[width=0.25\textwidth]{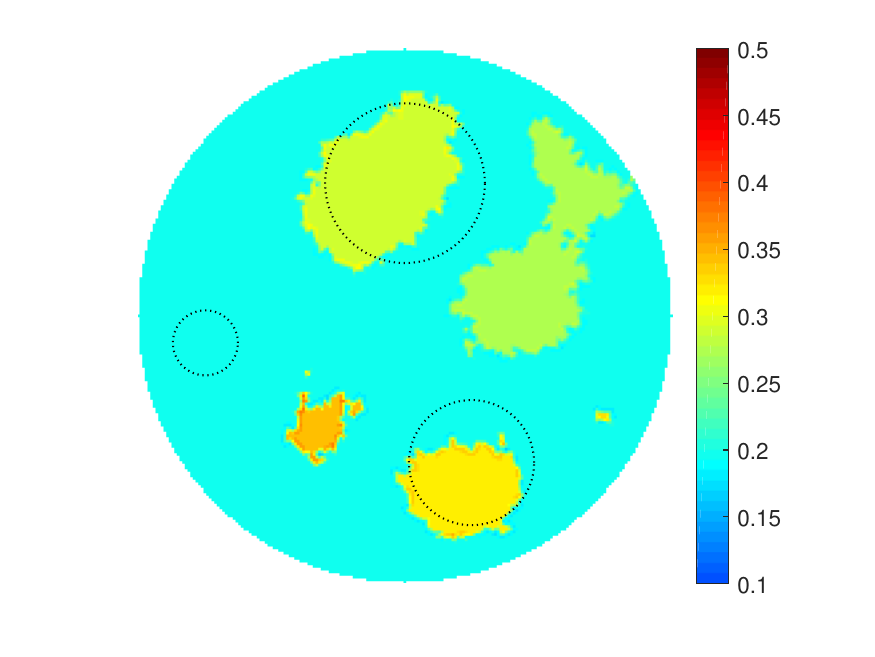}}}
  \mbox{
  \subfigure[$\hat \sigma$: fine scale, $0.25\%$]
    {\includegraphics[width=0.25\textwidth]{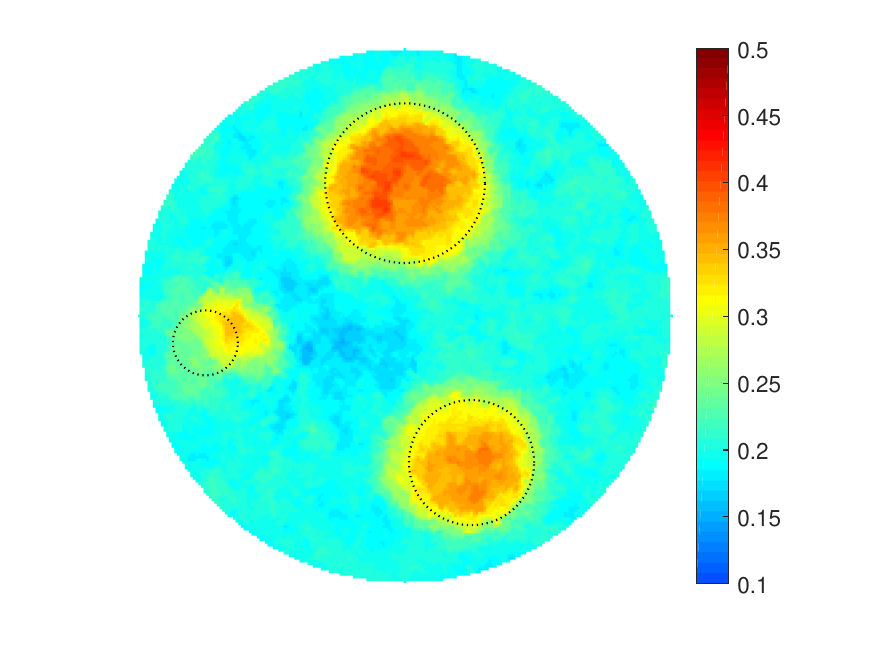}}
  \subfigure[$\hat \sigma$: fine scale, $0.5\%$]
    {\includegraphics[width=0.25\textwidth]{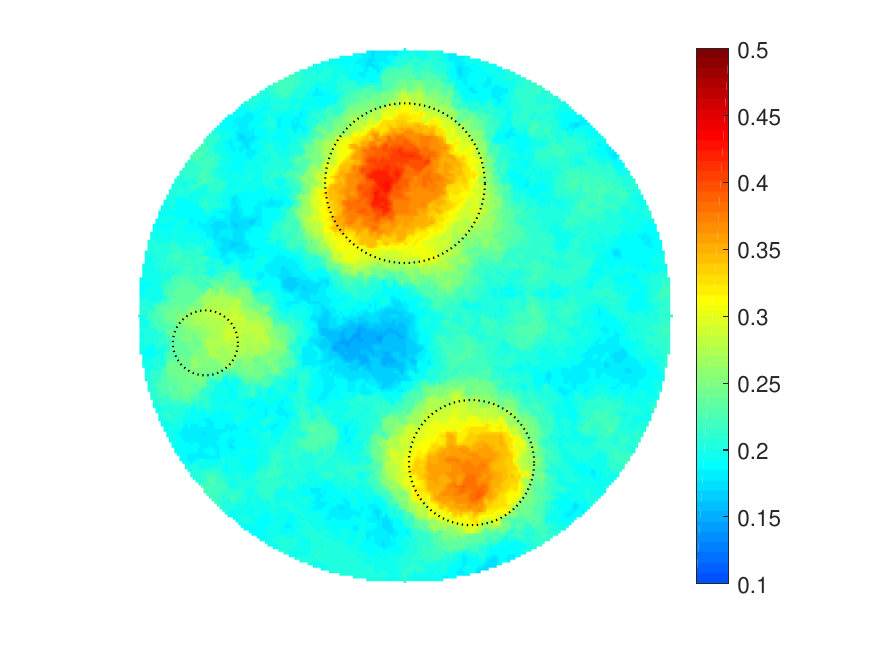}}
  \subfigure[$\hat \sigma$: fine scale, $1.0\%$]
    {\includegraphics[width=0.25\textwidth]{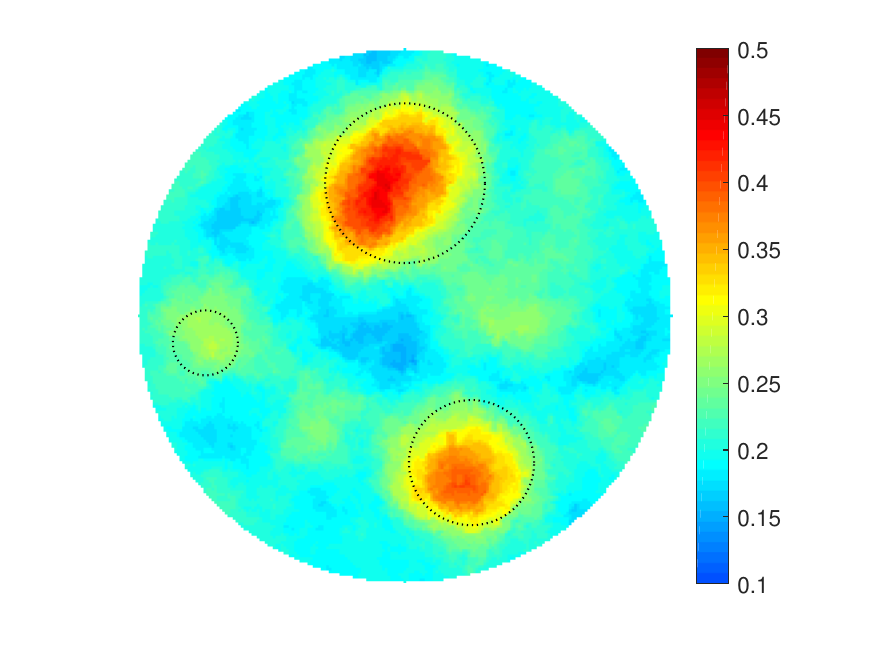}}
  \subfigure[$\hat \sigma$: fine scale, $2.0\%$]
    {\includegraphics[width=0.25\textwidth]{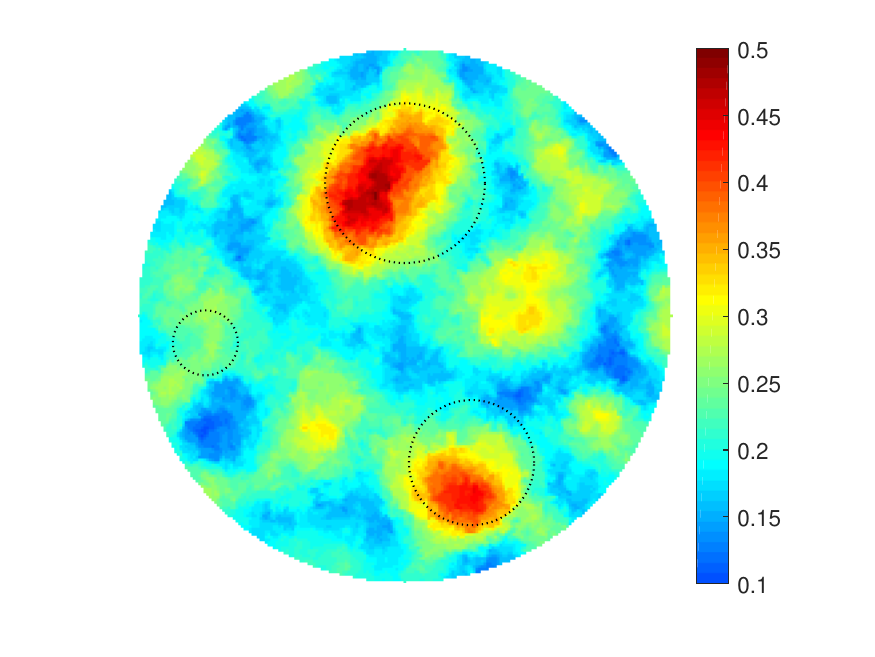}}}
  \end{center}
  \caption{Optimization results for model~\#1 obtained with measurements contaminated with
    (a,e,i)~0.25\%, (b,f,j)~0.5\%, (c,g,k)~1.0\%, and (d,h,l)~2.0\% noise. The images in (a-d)
    are obtained using only fine scales. (e-h) and (i-l) show images obtained at the coarse
    and fine scales, respectively, by using multiscale optimization with tuned PCA scale
    switching. The dashed circles are added to all images to represent the location of cancer-affected
    regions taken from known $\sigma_{true}(x)$ in Figure~\ref{fig:geometry}(a).}
  \label{fig:noise_sol}
\end{figure}

In Figure~\ref{fig:noise_sol}, we revisit our original model~\#1, with the true electrical conductivity $\sigma_{true}(x)$
provided in Figure~\ref{fig:geometry}(a), now with measurements contaminated with 0.25\%, 0.5\%, 1\%, and 2\% normally
distributed noise. As expected, we see that various levels of noise lead to oscillatory instabilities in the images
reconstructed by using only fine scales and fixed parameterization via PCA; refer to Figures~\ref{fig:noise_sol}(a-d).
If used practically in screening procedures, such imaging will obviously result in multiple cases of false positive outcomes.
On the other hand, as seen in Figures~\ref{fig:noise_sol}(e-h), the proposed computational algorithm with multilevel
parameterization demonstrates its stable ability to provide clear and accurate images with the appearance of false
positive or false negative results for some regions only with noise higher than 1-2\%. Also, as shown in
Figures~\ref{fig:noise_sol}(i-l), this new approach prevents the fine-scale solutions from the negative impacts caused by
propagated noise -- a noticeable distortion in the fine-scale images starts when the noise level passes beyond the 2\% level.
This result is impressive as it concludes the ability of the coarse-scale solutions ``properly disclosed'' to the fine scales
to improve ``noise resistance'' in general for solutions at both scales.

\begin{figure}[!htb]
  \begin{center}
  \mbox{
  \subfigure[objectives: coarse scale]
    {\includegraphics[width=0.33\textwidth]{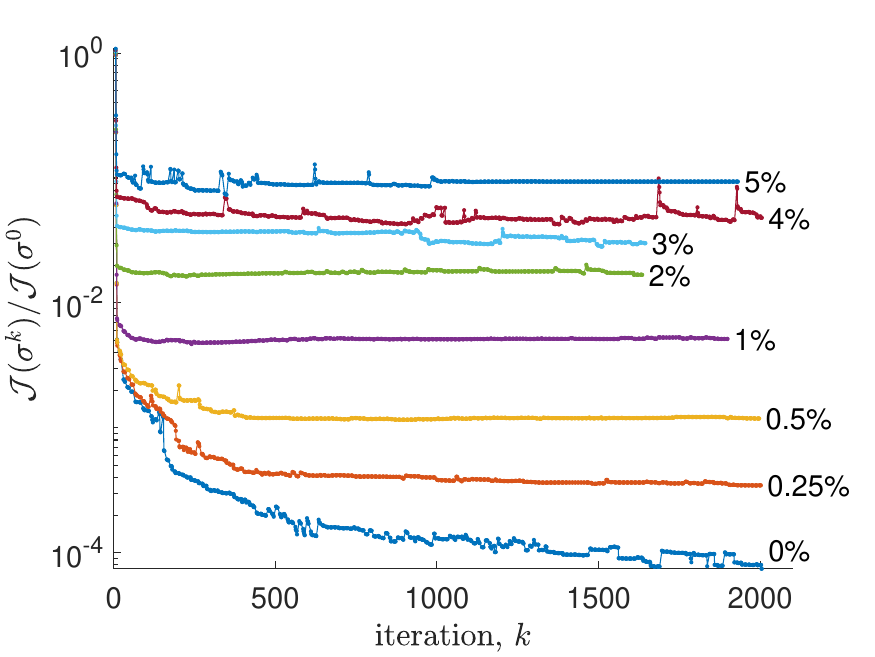}}
  \subfigure[solution error: coarse scale]
    {\includegraphics[width=0.33\textwidth]{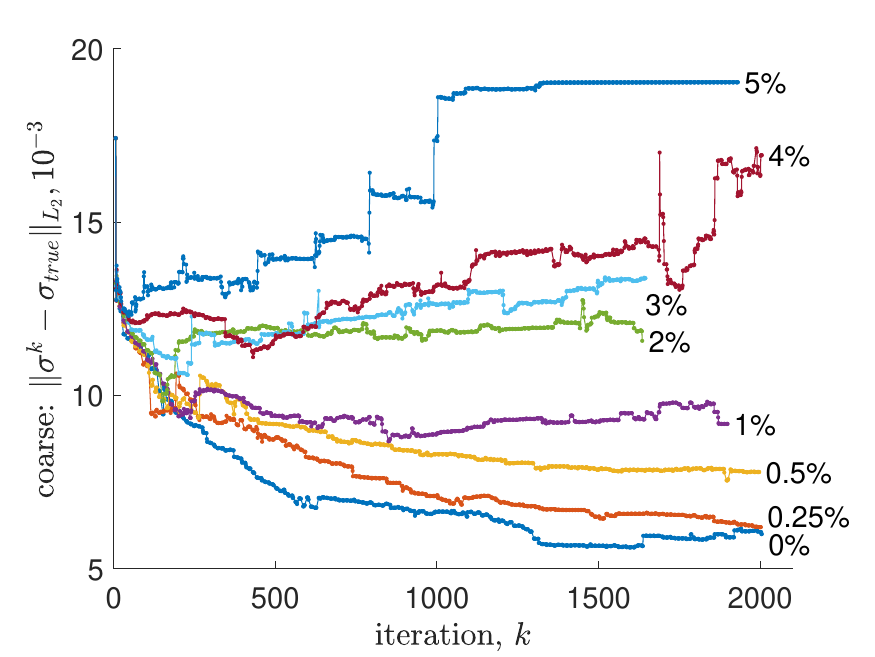}}
  \subfigure[solution error: fine scale]
    {\includegraphics[width=0.33\textwidth]{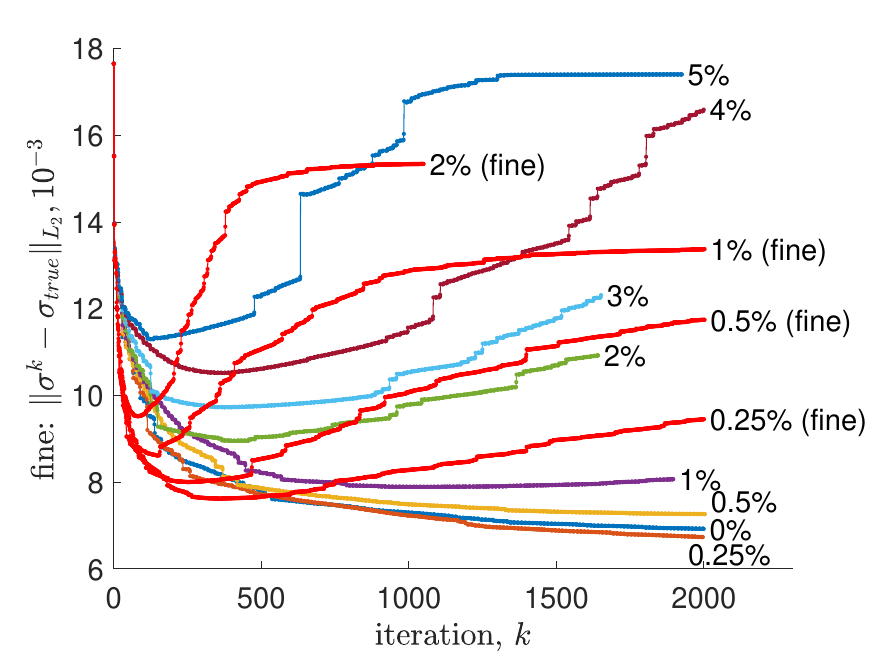}}}
  \end{center}
  \caption{(a)~Normalized objective functions $\mJ(\sigma^k)/\mJ(\sigma^0)$ and (b,c)~solution errors
    $\| \sigma^k - \sigma_{true} \|_{L_2}$ as functions of iteration count $k$ evaluated at (b)~coarse
    and (c)~fine scales for various noise levels from 0\% to 5\%. In (c), red curves represent solution
    errors for the fine-scale-only images obtained for noise ranging from 0.25\% to 2\%.}
  \label{fig:noise_comp}
\end{figure}

Finally, Figure~\ref{fig:noise_comp} provides a more thorough comparative study of the noise influence on the results
of obtaining EIT images at both coarse and fine scales. E.g., Figures~\ref{fig:noise_comp}(a,b) compare the objective
functions and solution errors evaluated at the coarse scale only for various noise levels up to 5\%. From both plots,
it is evident that the noise of 1\% and below has little effect on the quality of the reconstructed binary distributions
in both aspects, namely, an ability to match data properly (even noisy) and to generate a solution to the optimization
problem with relatively small error. Figure~\ref{fig:noise_comp}(b) also demonstrates that noise of 2\% and above forces
optimization to overfitting, a known effect in inverse problems supplied with highly noisy data. These conclusions are
consistent with the results seen in Figure~\ref{fig:noise_comp}(c). It shows the solution error evaluated at the fine
scale that deviates within a small interval for the noise levels from 0\% to 1\%. Overfitting starts with a noise of
2\% and up. This figure also shows the solution errors evaluated for fine-scale-only images (red curves) obtained with
different noise levels to compare them with those when a proposed multiscale optimization with the tuned PCA scale
switching is in use. It adds more to support our previous statement on the gained ability of the new approach to improve
its resistance to noise at both fine and coarse scales.

\subsection{Model~\#2: Real Breast Cancer Case}
\label{sec:model_2_new}

\begin{figure}[!htb]
  \begin{center}
  \mbox{
  \subfigure[mammogram image]
    {\includegraphics[width=0.25\textwidth]{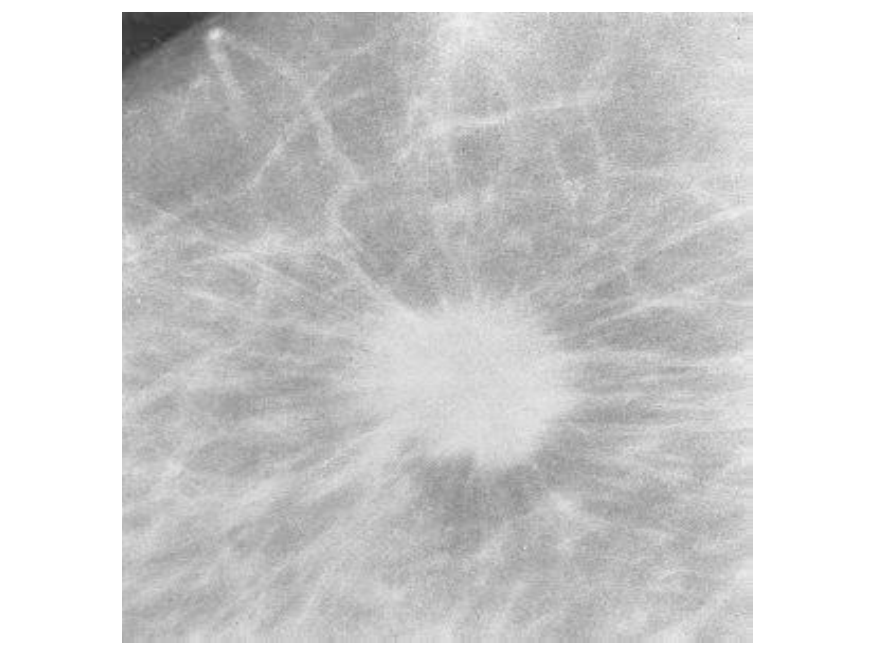}}
  \subfigure[re-mapping]
    {\includegraphics[width=0.25\textwidth]{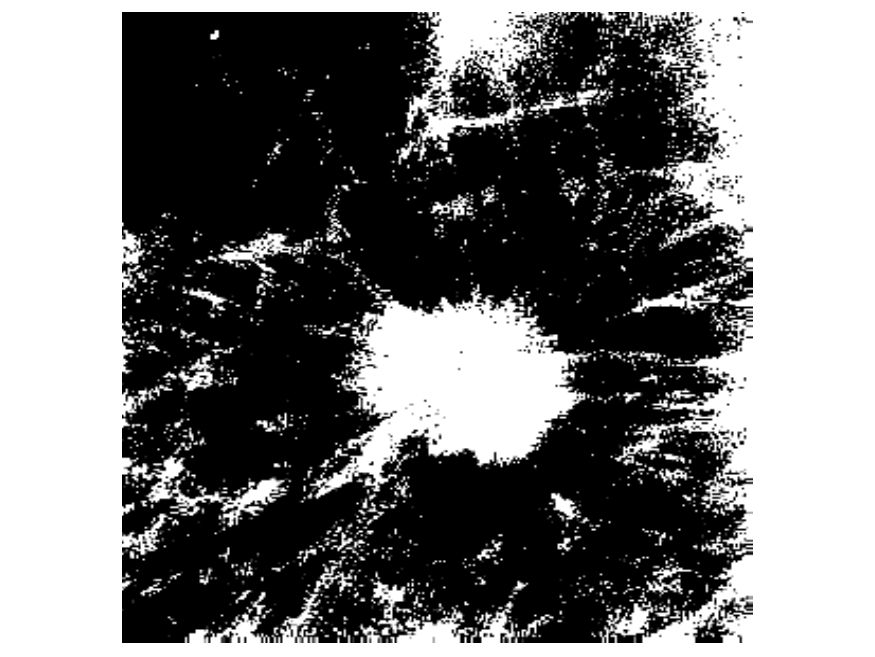}}
  \subfigure[isolating]
    {\includegraphics[width=0.25\textwidth]{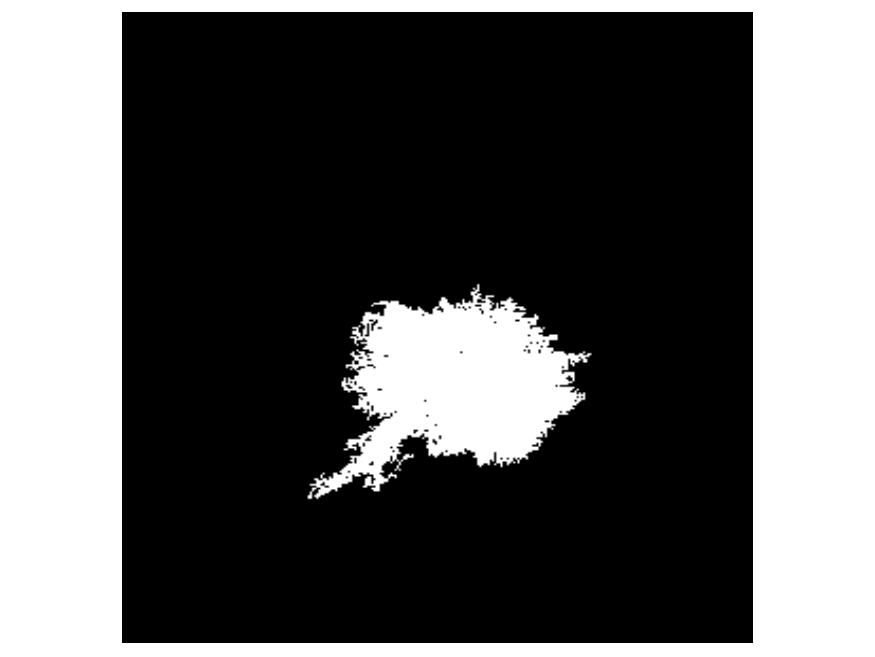}}
  \subfigure[final version]
    {\includegraphics[width=0.25\textwidth]{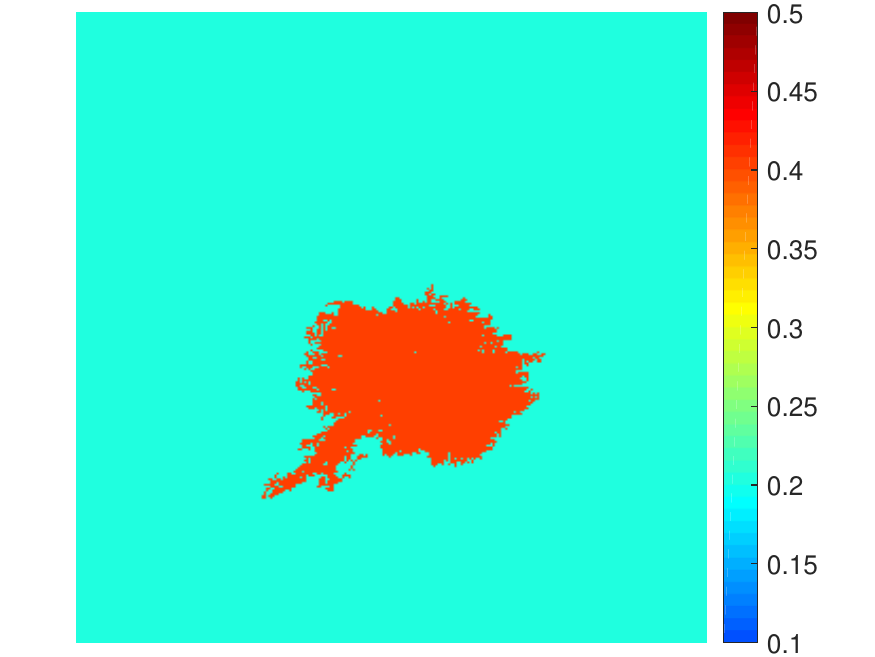}}}
  \end{center}
  \caption{(a)~The real breast cancer image for invasive ductal carcinoma (source: Bassett et al., 2003, \cite{Bassett2003}).
    (b-d)~A conversion process to create model~\#2 using filtering techniques in {\tt MATLAB} that include
    (b)~re-mapping, (c)~isolating, (d)~smoothing, and exporting to the data format.}
  \label{fig:m7_create}
\end{figure}

In the second part of our numerical experiments with the proposed optimization framework, we are particularly interested
in applying it to the cases seen in the medical practice during cancer-related screening procedures. We created our next
model (\#2) based on a mammogram image of a real breast cancer case available in \cite{Bassett2003}; refer to
Figure~\ref{fig:m7_create}(a). This model shows an invasive ductal carcinoma with an irregular shape and spiculated margins.
Both of these properties are characteristic of malignant masses and are two criteria radiologists would look for when
assessing a scan for breast cancer. Due to our incapability to produce actual measuring of the electrical currents, we
have to convert the mammogram image to its binary version and obtain synthetic data in place of the real measurements,
as discussed in Section~\ref{sec:math_model}. We performed this conversion by using filtering techniques in
{\tt MATLAB}\textsuperscript{\textregistered} that include three main stages, namely, re-mapping, isolating, and smoothing:
\begin{enumerate}
  \item[(1)] First, the black and white shades in the original image are re-mapped using {\tt MATLAB}'s function {\tt imadjust},
    supplied with a cutoff value. All pixels with the shades below that value turn black, and all above turn white to enhance
    the region of interest by increasing the contrast within a specific range; refer to Figure~\ref{fig:m7_create}(b).
  \item[(2)] Because tissues have variable density naturally, it is expected that the previous step leaves behind some small
    spots that are not part of the cancerous region, including blood vessels and small ducts. We remove them to leave only
    the cancerous one(s) by ``erasing'' shapes with a total pixel count below a specified value using {\tt MATLAB}'s function
    {\tt bwareaopen}; see Figure~\ref{fig:m7_create}(c).
  \item[(3)] It is natural that throughout steps (1) and (2), the edges of the cancerous region become jagged. In pursuit of
    better shapes (to be reconstructed), we find the region edges using {\tt MATLAB}'s function {\tt edge} and the {\tt Canny}
    edge detection method. Then these edges are dilated, and the region is re-filled to reveal the same shape as was detected
    prior, but with slightly smoother boundaries. The final version is then exported to a data file in a format consistent
    with the input requirements of the computational framework. The healthy tissue is represented by a value of $\sigma_h = 0.2$,
    whereas the cancerous region is assigned a value of $\sigma_c = 0.4$; refer to Figure~\ref{fig:m7_create}(d).
\end{enumerate}

\begin{figure}[!htb]
  \begin{center}
  \mbox{
  \subfigure[model~\#2]
    {\includegraphics[width=0.33\textwidth]{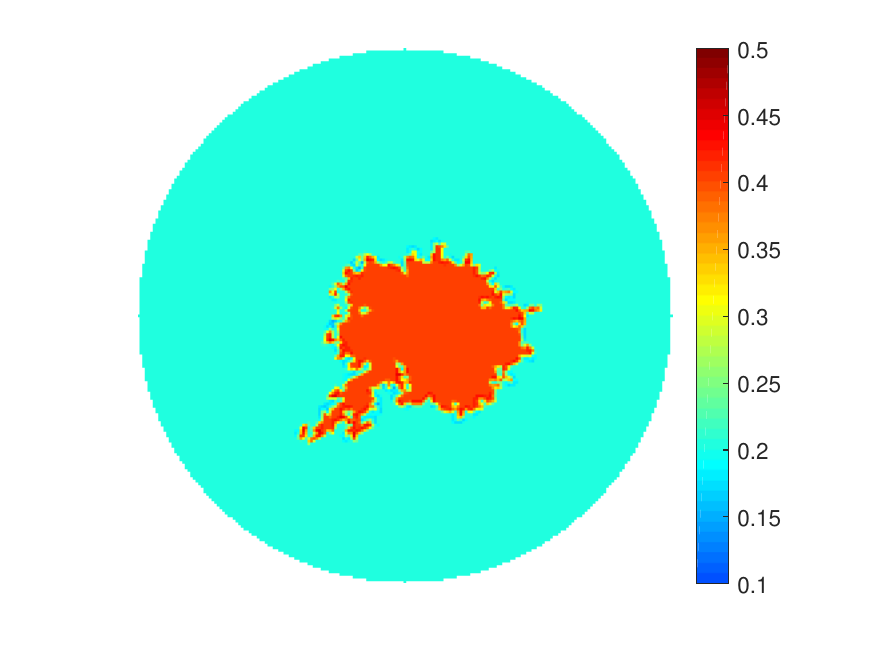}}
  \subfigure[$\hat \sigma$: fine scale only]
    {\includegraphics[width=0.33\textwidth]{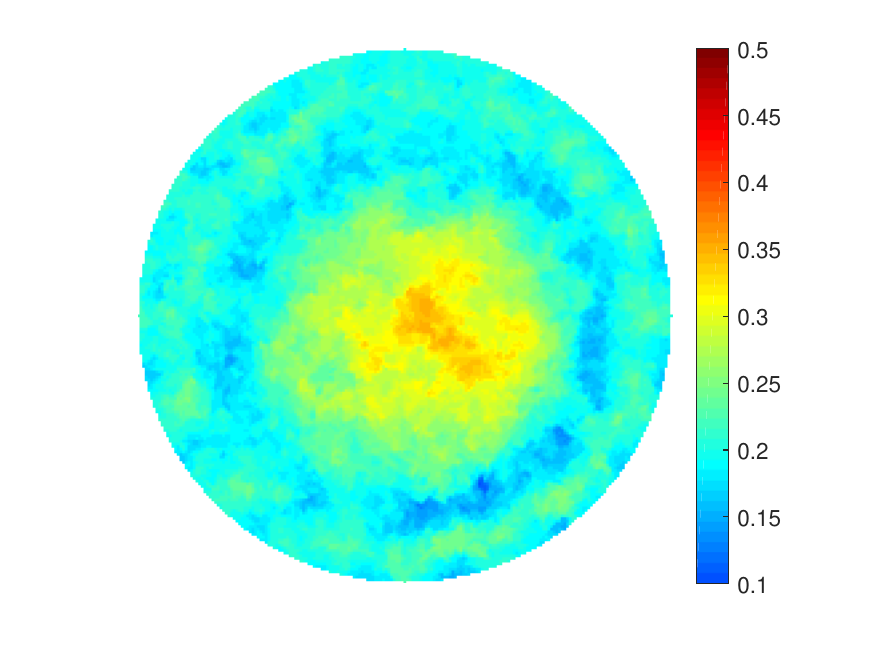}}
  \subfigure[$\hat \sigma$: $N_{\max} = 5$ \& PCA]
    {\includegraphics[width=0.33\textwidth]{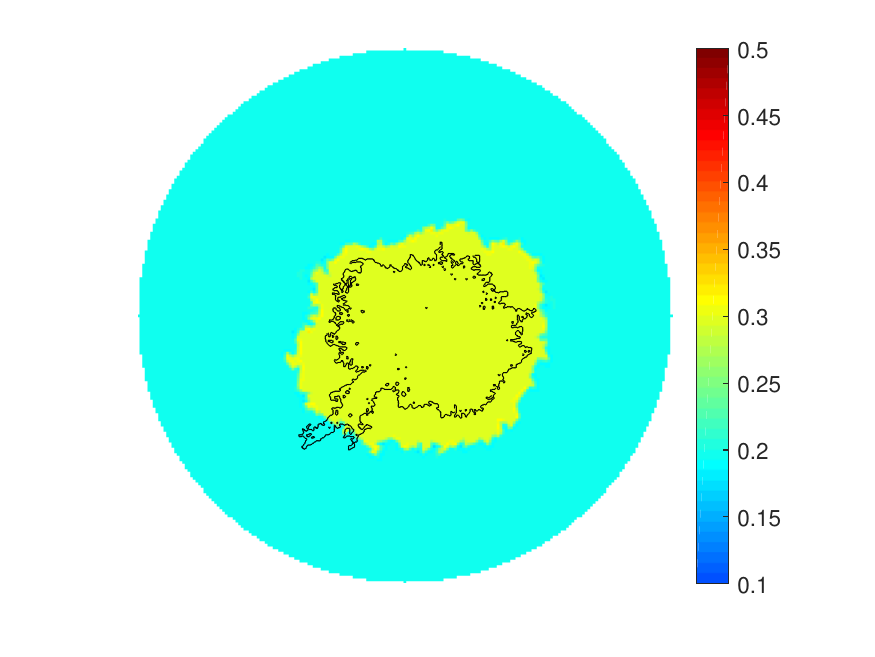}}}
  \mbox{
  \subfigure[objectives]
    {\includegraphics[width=0.33\textwidth]{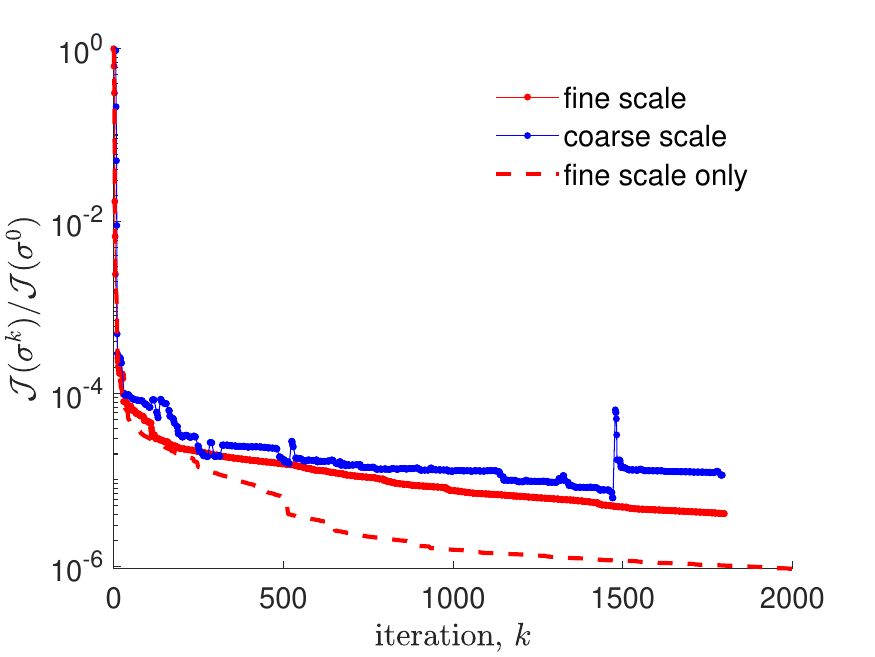}}
  \subfigure[solution error]
    {\includegraphics[width=0.33\textwidth]{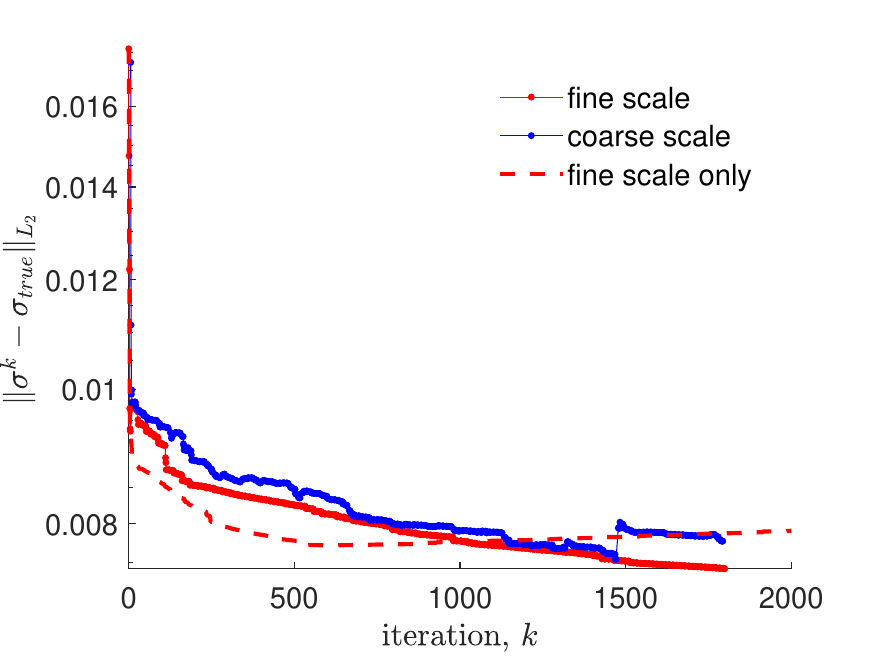}}}
  \end{center}
  \caption{(a)~EIT model~\#2: true electrical conductivity $\sigma_{true}(x)$. (b,c)~The images of model~\#2 obtained
    after optimizing (b)~only at the fine scale and (c)~at multiple scales with the enhanced scale switching by the
    tuned PCA when $N_{\max} = 5$. In (c), the added black curve represents the boundary of the cancer-affected
    region taken from known $\sigma_{true}(x)$ in (a). (d)~Normalized objective functions $\mJ(\sigma^k)/\mJ(\sigma^0)$
    and (e)~solution errors $\| \sigma^k - \sigma_{true} \|_{L_2}$ as functions of iteration count $k$ evaluated at
    fine (in red) and coarse (in blue) scales compared to the results in the fine-scale-only images (dashed lines).}
  \label{fig:m7_result}
\end{figure}

After creating model~\#2, see Figure~\ref{fig:m7_result}(a), we ran our computational framework with the complete suite of
the proposed multiscale optimization functionalities. Figures~\ref{fig:m7_result}(b,c) demonstrate the images obtained while
running optimization at the fine scale only and in the multiscale mode with PCA tuning activated during the scale switching
procedures and the maximum number of principal components $N_{\max} = 5$, respectively. Without any doubt, the solution to
the inverse EIT problem of model~\#2 is very challenging due to the nontrivial shape of the cancerous spot at the center.
Although the presence of cancer is evident, the fine-scale-only image in Figure~\ref{fig:m7_result}(b) does not provide much
help in identifying the exact location and boundaries. However, this solution fits data better than the one obtained at
multiple scales; see Figure~\ref{fig:m7_result}(d). Contrary to that, the binary image in Figure~\ref{fig:m7_result}(c)
has better resolution and is more informative. Although it seems a bit weaker in fitting data, the solution error in
Figure~\ref{fig:m7_result}(e) demonstrates better performance at both scales.

\begin{figure}[!htb]
  \begin{center}
  \mbox{
  \subfigure[objectives]
    {\includegraphics[width=0.5\textwidth]{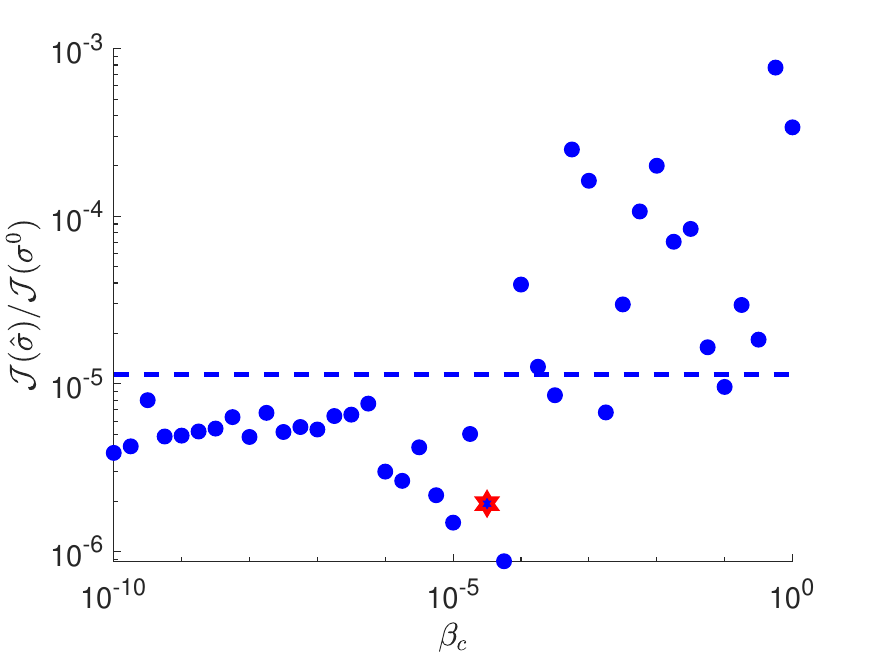}}
  \subfigure[solution error]
    {\includegraphics[width=0.5\textwidth]{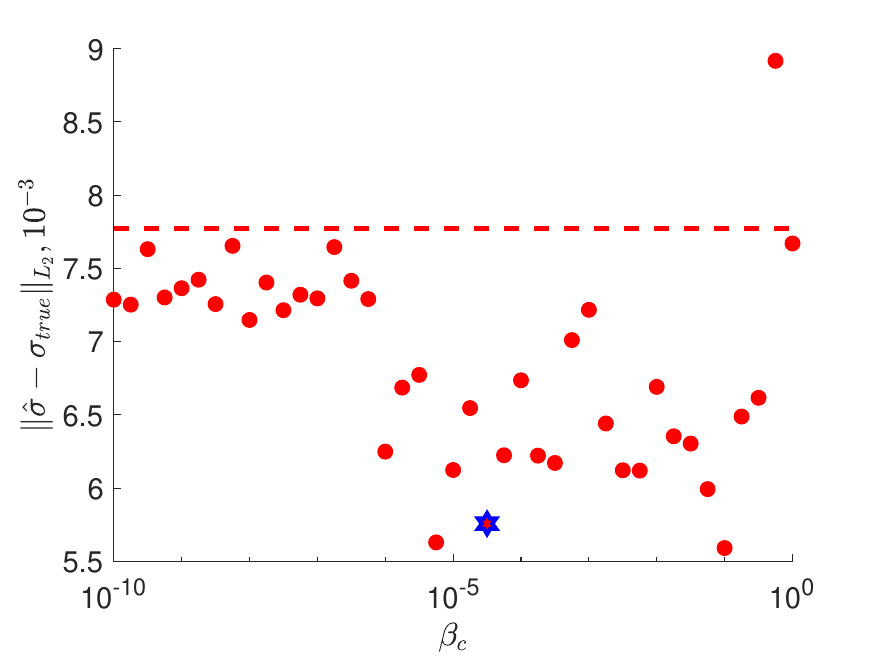}}}
  \end{center}
  \caption{(a)~Normalized objective functions $\mJ(\hat \sigma)/\mJ(\sigma^0)$ and (b)~solution errors
    $\| \hat \sigma - \sigma_{true} \|_{L_2}$ as functions of regularization parameter $\beta_c$ in \eqref{eq:reg_c}
    evaluated for optimal solutions $\hat \sigma$ at the coarse scale compared to the results obtained with no
    regularization ($\beta_c = 0$) applied (dashed lines). In both plots, hexagons depict the results for
    $\beta_c = 10^{-4.5}$.}
  \label{fig:m7_reg_comp}
\end{figure}

To improve these results and the overall performance of our computational framework, we apply the regularization at the
coarse scale, as discussed in Section~\ref{sec:penal}, in the assumption of existing knowledge of cancer present in
the tissue, i.e., two constant values $\bar \sigma_l = 0.2$ and $\bar \sigma_h = 0.4$ are given. We ran multiscale optimization
with various values of the weighting parameter $\beta_c$ in \eqref{eq:reg_c} ranging from $10^{-10}$ to $10^0$.
Figure~\ref{fig:m7_reg_comp} compares the obtained results in terms of the objectives (only core part $\mJ$ in
\eqref{eq:obj_reg_full}, not including the regularization component $\mJ_c$) and solution errors both evaluated for the
coarse-scale optimal solutions. As the target values for controls $\sigma_{low}$ and $\sigma_{high}$ ($\bar \sigma_l$ and
$\bar \sigma_h$, respectively) are inputted through the regularization part $\mJ_c$ of the objective $\bar \mJ$, we expect
an improved performance by obtaining images with more accurate shapes. Figure~\ref{fig:m7_reg_comp}(a) suggests this
improvement for $\beta_c$ smaller than $10^{-4}$ (best results are assumed between $10^{-6}$ and $10^{-4}$) when the majority
of the optimization runs finalize with the results (blue dots) better than when no optimization is applied (dashed line).
However, the posterior assessments in Figure~\ref{fig:m7_reg_comp}(b), i.e., solution errors, give a different range
for the best results, from $10^{-6}$ to $10^0$.

\begin{figure}[!htb]
  \begin{center}
  \mbox{
  \subfigure[$\hat \sigma$: $\beta_c = 10^{-4.5}$]
    {\includegraphics[width=0.33\textwidth]{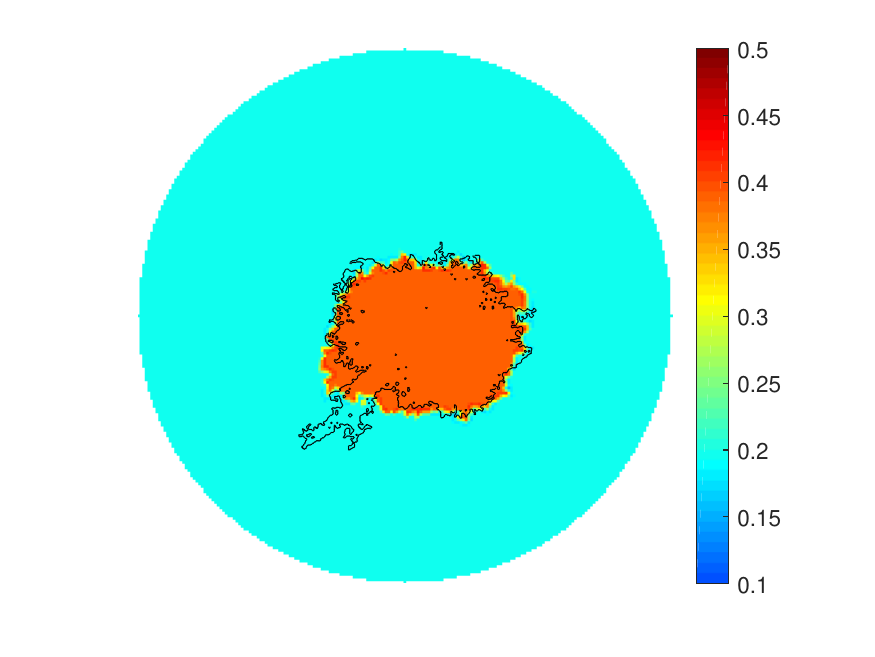}}
  \subfigure[objectives]
    {\includegraphics[width=0.33\textwidth]{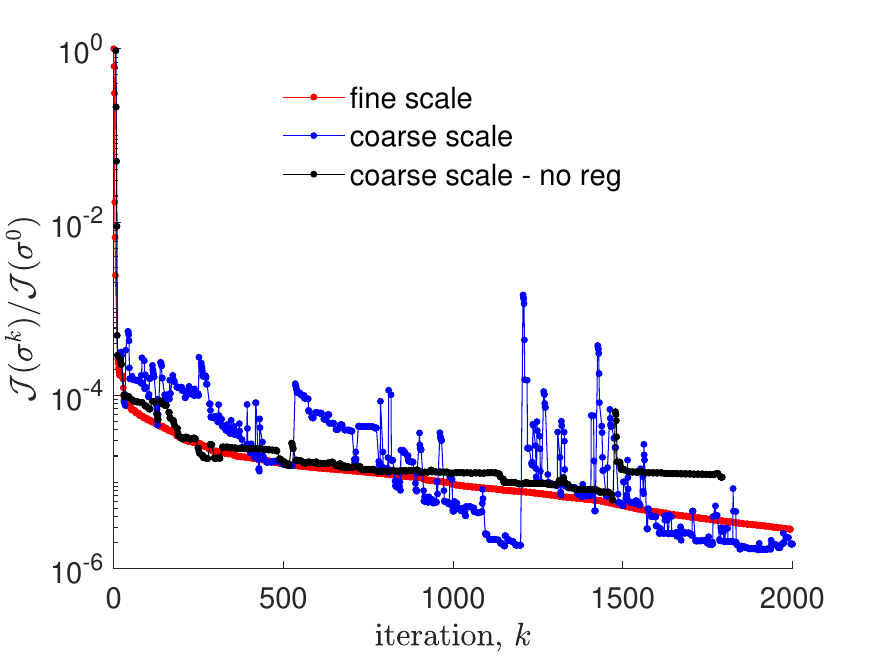}}
  \subfigure[solution error]
    {\includegraphics[width=0.33\textwidth]{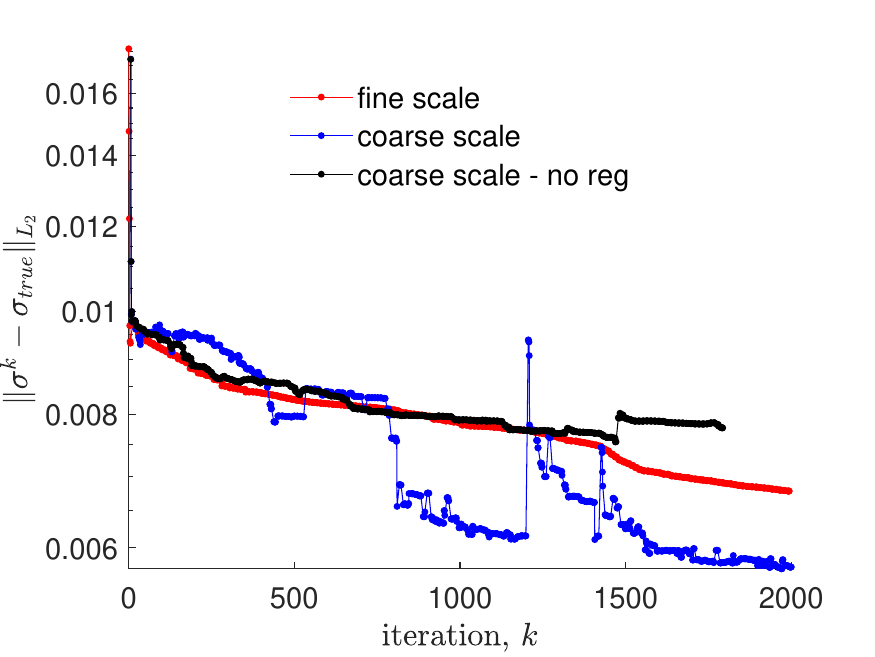}}}
  \end{center}
  \caption{(a)~The image of model~\#2 obtained after optimizing at multiple scales with the enhanced scale switching
    by the tuned PCA when $N_{\max} = 5$ and applied regularization at the coarse scale with weight $\beta_c = 10^{-4.5}$
    in \eqref{eq:reg_c}. The added black curve represents the boundary of the cancer-affected region taken from known
    $\sigma_{true}(x)$ in Figure~\ref{fig:m7_result}(a). (b)~Normalized objective functions $\mJ(\sigma^k)/\mJ(\sigma^0)$
    and (c)~solution errors $\| \sigma^k - \sigma_{true} \|_{L_2}$ as functions of iteration count $k$ evaluated at
    fine (in red) and coarse (in blue) scales compared to the results in the coarse-scale images obtained with no
    regularization ($\beta_c = 0$) applied (black curves).}
  \label{fig:m7_reg_best}
\end{figure}

To check the quality of the solutions obtained in the overlapped interval of $\beta_c$ (namely, between $10^{-6}$
and $10^{-4}$), we choose one value with the corresponding outcomes shown in both plots of Figure~\ref{fig:m7_reg_comp}
as hexagons. The results of the multiscale optimization with added coarse-scale regularization for $\beta_c = 10^{-4.5}$
are provided in Figure~\ref{fig:m7_reg_best}. Here, the coarse-scale image in (a) has significantly improved shape. As
seen in Figure~\ref{fig:m7_reg_best}(b), the gap between the objectives evaluated at the fine and coarse scales is minimal
(not considering the oscillation peaks), reflecting preserved proper communication between solutions at both scales.
Finally, Figure~\ref{fig:m7_reg_best}(c) confirms the notably higher quality of the coarse-scale solution compared to the
one obtained at the fine scale. Here, we conclude that the suitably chosen and properly tuned regularization shows potential
for further improvements in the applications of the proposed approach to real models despite their known complexity.

\subsection{Model~\#3: More Complicated Case of Breast Cancer}
\label{sec:model_3_new}

\begin{figure}[!htb]
  \begin{center}
  \mbox{
  \subfigure[MRI image]
    {\includegraphics[width=0.33\textwidth]{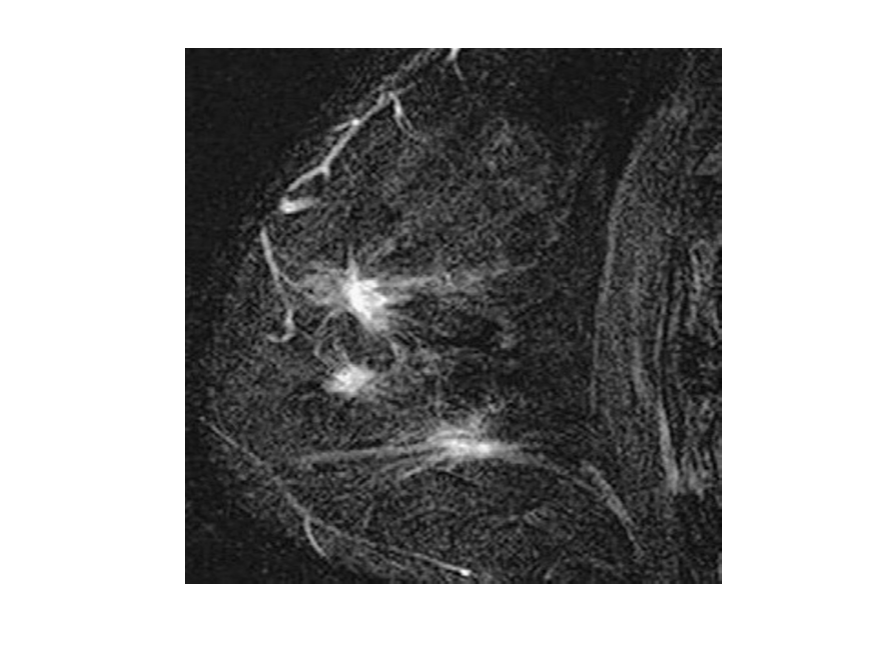}}
  \subfigure[model~\#3]
    {\includegraphics[width=0.33\textwidth]{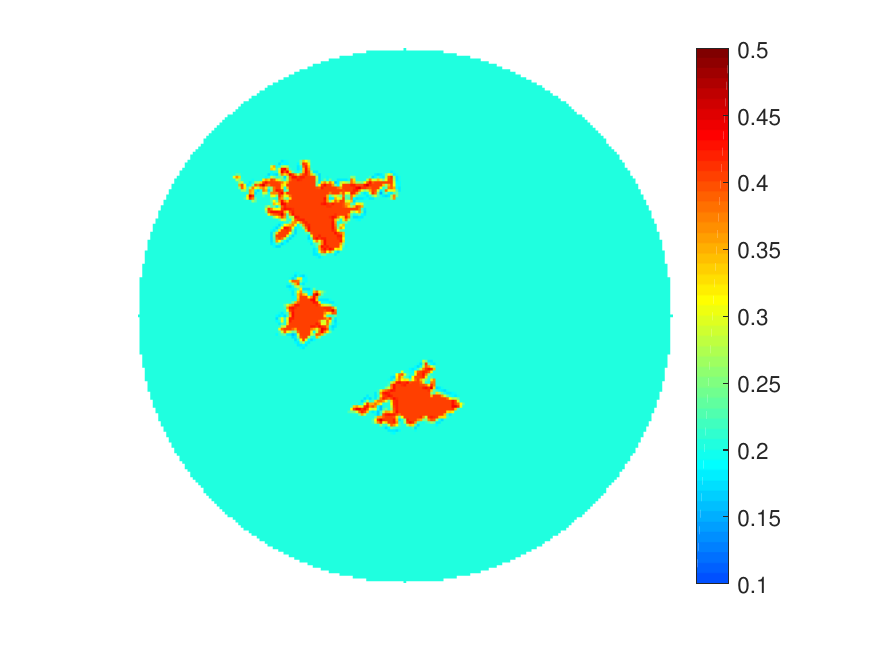}}
  \subfigure[$\hat \sigma$: fine scale only]
    {\includegraphics[width=0.33\textwidth]{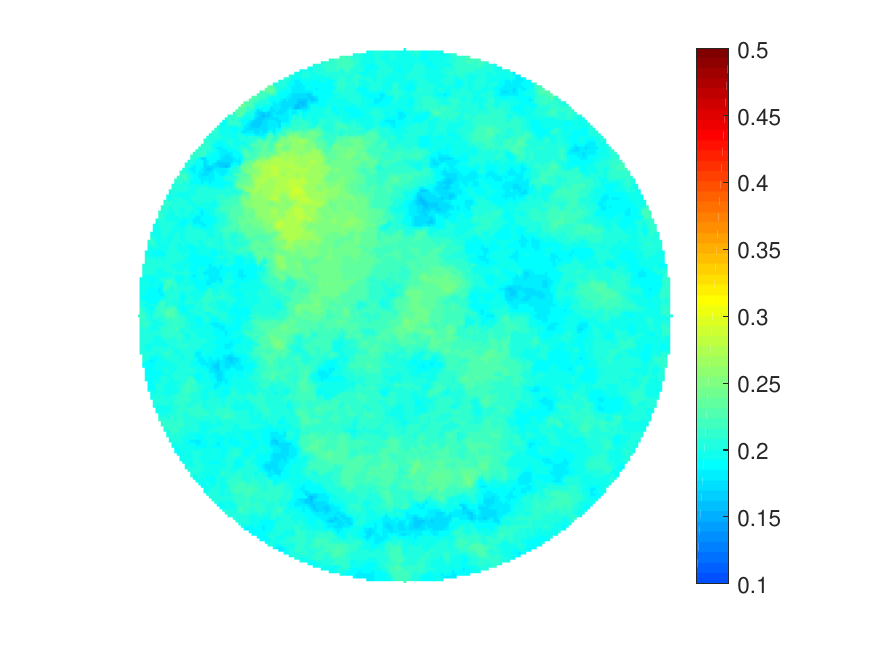}}}
  \mbox{
  \subfigure[$\hat \sigma$: $N_{\max} = 5$ \& PCA]
    {\includegraphics[width=0.33\textwidth]{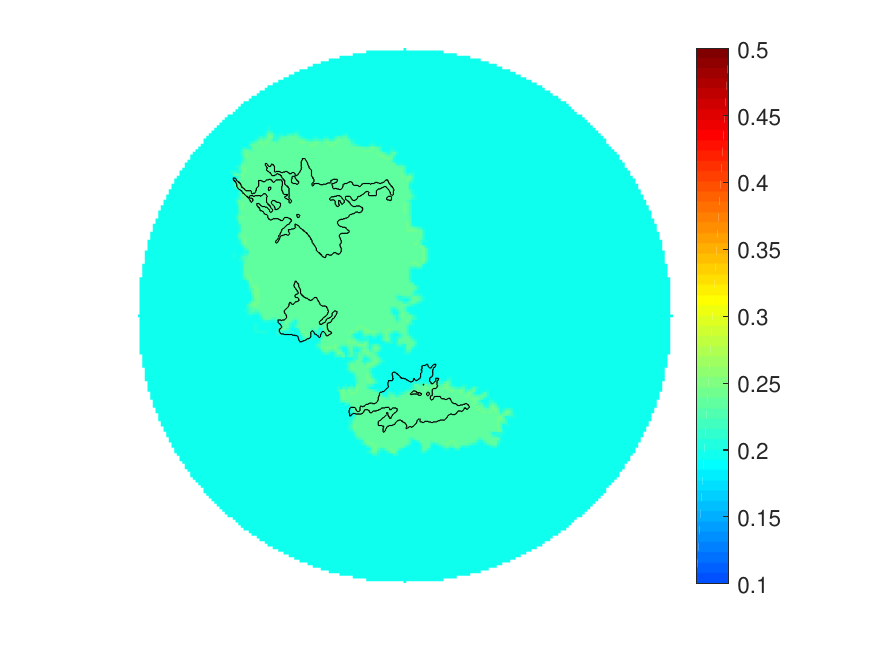}}
  \subfigure[objectives]
    {\includegraphics[width=0.33\textwidth]{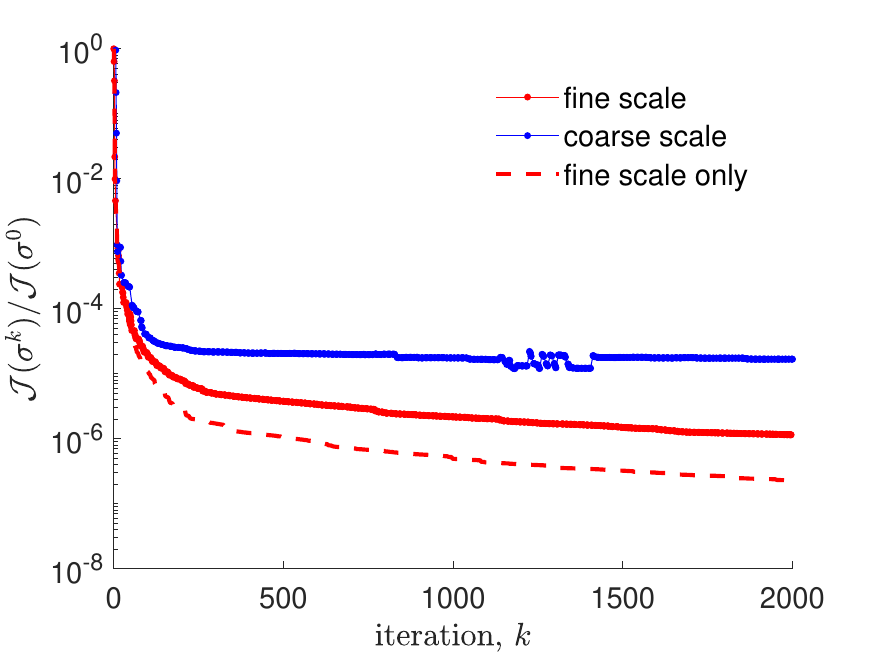}}
  \subfigure[solution error]
    {\includegraphics[width=0.33\textwidth]{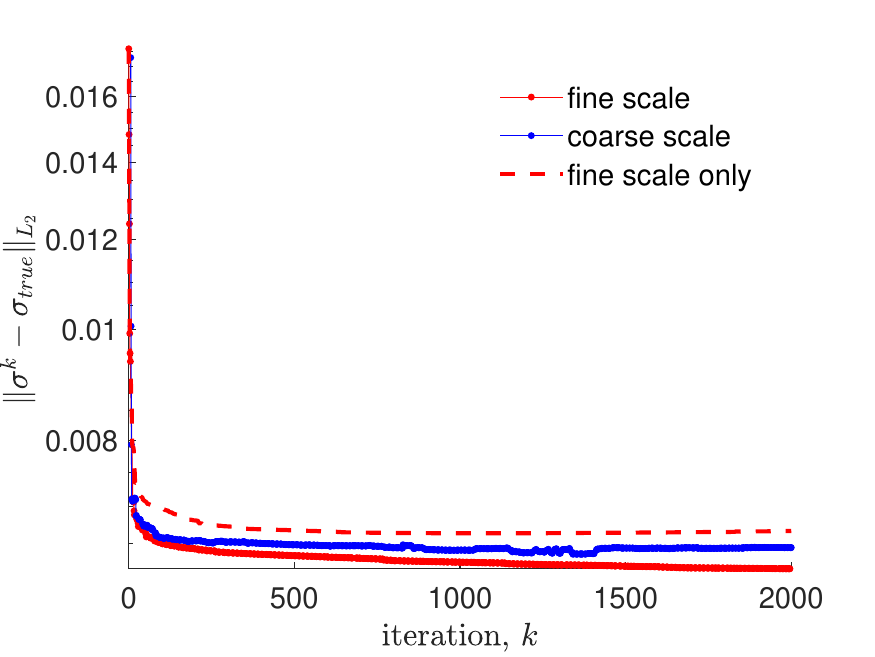}}}
  \end{center}
  \caption{(a)~The real breast cancer image for invasive ductal carcinoma (source: Weinstein, 2009, \cite{Weinstein2009}).
    (b)~EIT model~\#3: true electrical conductivity $\sigma_{true}(x)$. (c,d)~The images of model~\#3 obtained
    after optimizing (c)~only at the fine scale and (d)~at multiple scales with the enhanced scale switching by the
    tuned PCA when $N_{\max} = 5$. In (d), the added black curves represent the boundaries of the cancer-affected
    regions taken from known $\sigma_{true}(x)$ in (b). (e)~Normalized objective functions $\mJ(\sigma^k)/\mJ(\sigma^0)$
    and (f)~solution errors $\| \sigma^k - \sigma_{true} \|_{L_2}$ as functions of iteration count $k$ evaluated at
    fine (in red) and coarse (in blue) scales compared to the results in the fine-scale-only images (dashed lines).}
  \label{fig:m16_result}
\end{figure}

In the final set of our numerical experiments, our focus is on even more complicated cases seen in the medical practice
during cancer-related screening procedures when multiple regions suspicious of cancer are present and characterized
by different sizes and nontrivial shapes. We created our last model (\#3) based on an MRI image of another real breast
cancer case available in \cite{Weinstein2009}; refer to Figure~\ref{fig:m16_result}(a). This model shows multiple
(at least three) spots identified as invasive ductal carcinoma with irregular shapes and spiculated margins.
Similar to our model~\#2, we converted the MRI image to its binary version to obtain synthetic data in place of the
real measurements following the same {\tt MATLAB}-assisted filtering methodology discussed in detail in
Section~\ref{sec:model_2_new}. Figure~\ref{fig:m16_result}(b) shows the ``true'' image of model~\#3, where colors
represent the binary distribution of the electrical conductivity $\sigma_{true}$ we aim to reconstruct.

As before, we ran our computational framework with the complete suite of the proposed multiscale optimization functionalities. Figures~\ref{fig:m16_result}(c,d) demonstrate the images obtained while running optimization at the fine scale only and
in the multiscale mode with PCA tuning activated during the scale switching procedures and the maximum number of principal
components $N_{\max} = 5$, respectively. Here, we must admit that the solution to the inverse EIT problem of model~\#3 is
even more challenging due to the presence of multiple spots, their small sizes, and nontrivial shapes. Unlike in the case
of model~\#2, the fine-scale-only image in Figure~\ref{fig:m16_result}(c) provides almost no information to help in identifying
at least the approximate locations of the cancerous spots. However, similar to model~\#2, this solution fits data better than
the one obtained at multiple scales; see Figure~\ref{fig:m16_result}(e). Contrary to that, the binary image in
Figure~\ref{fig:m16_result}(d) is undoubtedly more informative by providing some (rough) approximation to the spot locations
despite quite a large ``communication gap'' between objectives evaluated at the fine and coarse scales and its overall weakness
in fitting data. The analysis of the solution error in Figure~\ref{fig:m16_result}(f) also suggests better performance at both
scales.

\begin{figure}[!htb]
  \begin{center}
  \mbox{
  \subfigure[objectives]
    {\includegraphics[width=0.33\textwidth]{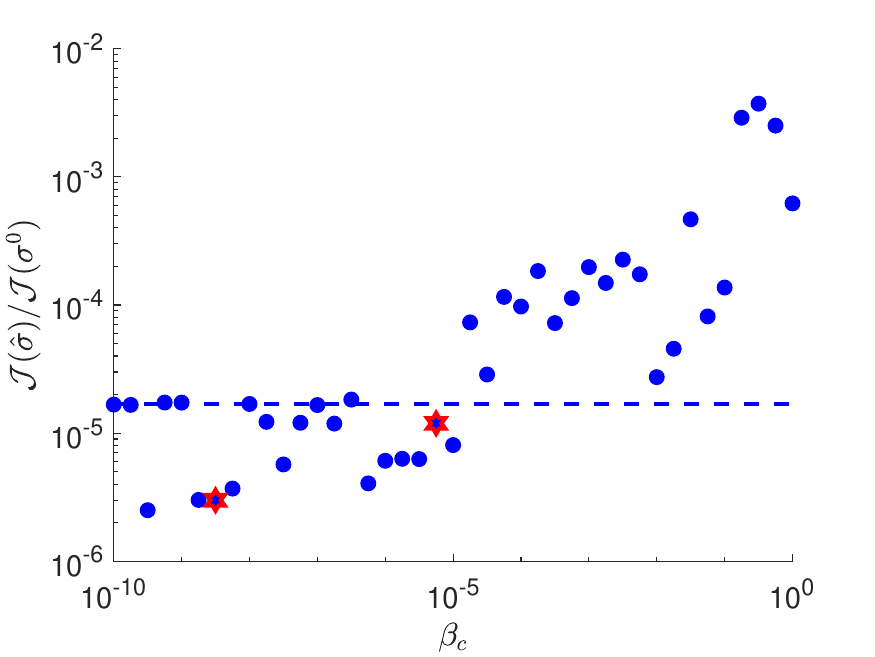}}
  \subfigure[solution error]
    {\includegraphics[width=0.33\textwidth]{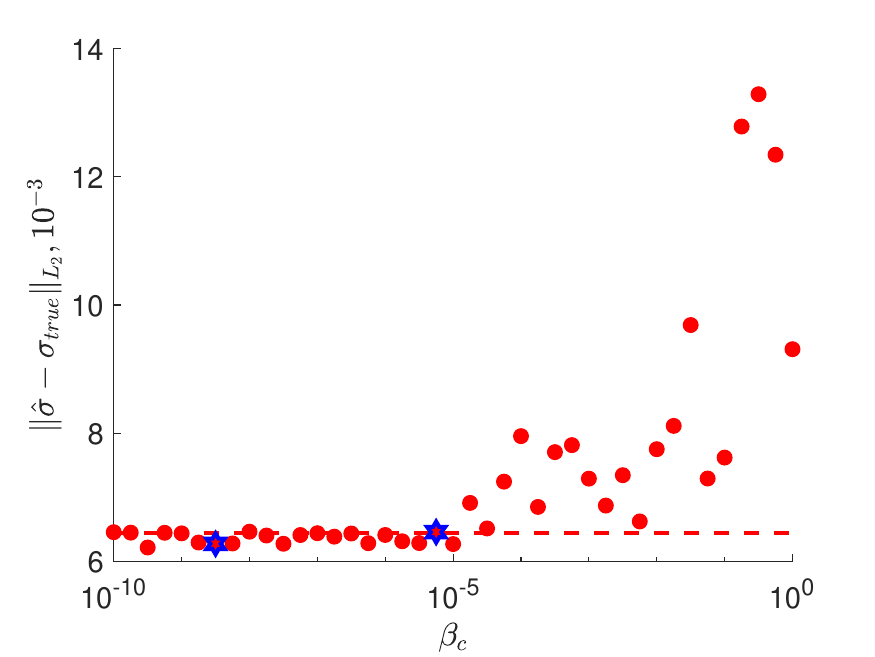}}
  \subfigure[solution error (close look)]
    {\includegraphics[width=0.33\textwidth]{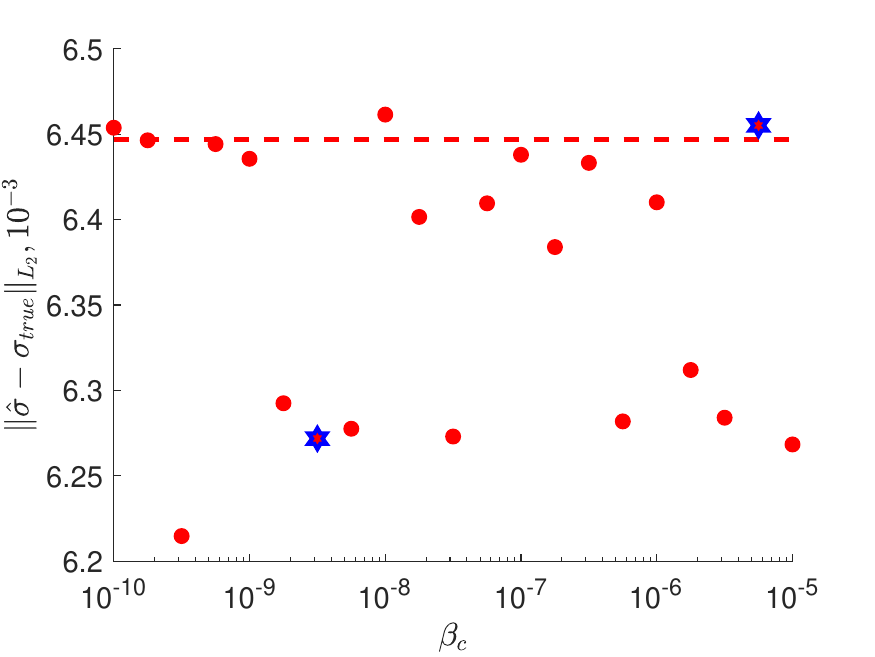}}}
  \end{center}
  \caption{(a)~Normalized objective functions $\mJ(\hat \sigma)/\mJ(\sigma^0)$ and (b,c)~solution errors
    $\| \hat \sigma - \sigma_{true} \|_{L_2}$ as functions of regularization parameter $\beta_c$ in \eqref{eq:reg_c}
    evaluated for optimal solutions $\hat \sigma$ at the coarse scale compared to the results obtained with no
    regularization ($\beta_c = 0$) applied (dashed lines). (c)~A close look at the results in (b) obtained
    with $\beta_c$ between $10^{-10}$ and $10^{-5}$. In all plots, hexagons depict the results for $\beta_c = 10^{-8.5}$
    and $\beta_c = 10^{-5.25}$.}
  \label{fig:m16_reg_comp}
\end{figure}

As we see a noticeable improvement in the performance of our computational framework in application to model~\#2,
see Section~\ref{sec:model_2_new} for details, now we apply the same type of regularization to the solutions
at the coarse scale, as discussed in Section~\ref{sec:penal}. Similarly, we assume the existence of some knowledge
of cancer present in the tissue by providing two constant values, $\bar \sigma_l = 0.2$ and $\bar \sigma_h = 0.4$,
and run multiscale optimization with various values of the weighting parameter $\beta_c$ in \eqref{eq:reg_c} ranging
from $10^{-10}$ to $10^0$. Figure~\ref{fig:m16_reg_comp} compares the obtained results in terms of the objectives
(only core part $\mJ$ in \eqref{eq:obj_reg_full}, not including the regularization component $\mJ_c$) and solution
errors, both evaluated for the coarse-scale optimal solutions. As before, we expect an improved performance by obtaining
images with more accurate shapes to enable a conclusion on cancer present at various locations. Figure~\ref{fig:m16_reg_comp}(a)
suggests this improvement for $\beta_c$ smaller than $10^{-5}$ when the majority of the optimization runs finalize with
the results (blue dots) better than when no optimization is applied (dashed line). The posterior assessments in
Figures~\ref{fig:m16_reg_comp}(b,c), i.e., solution errors, give a similar range.

\begin{figure}[!htb]
  \begin{center}
  \mbox{
  \subfigure[$\hat \sigma$: $\beta_c = 10^{-8.5}$]
    {\includegraphics[width=0.33\textwidth]{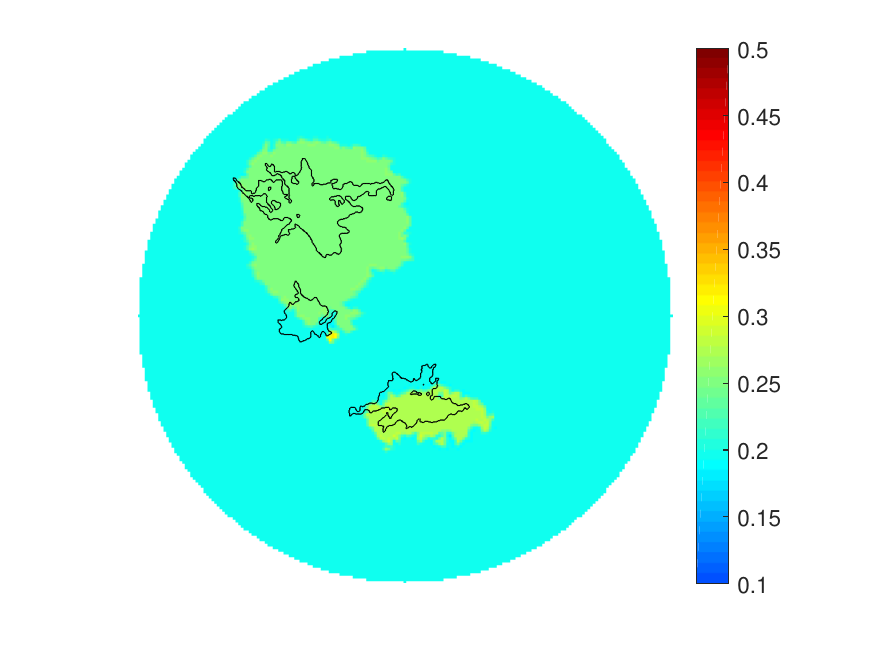}}
  \subfigure[objectives: $\beta_c = 10^{-8.5}$]
    {\includegraphics[width=0.33\textwidth]{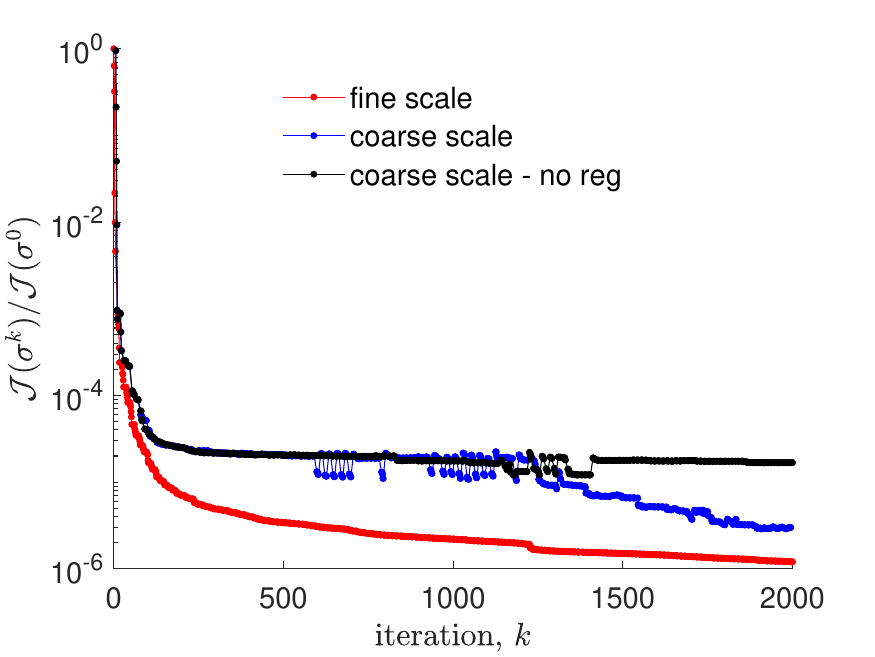}}
  \subfigure[solution error: $\beta_c = 10^{-8.5}$]
    {\includegraphics[width=0.33\textwidth]{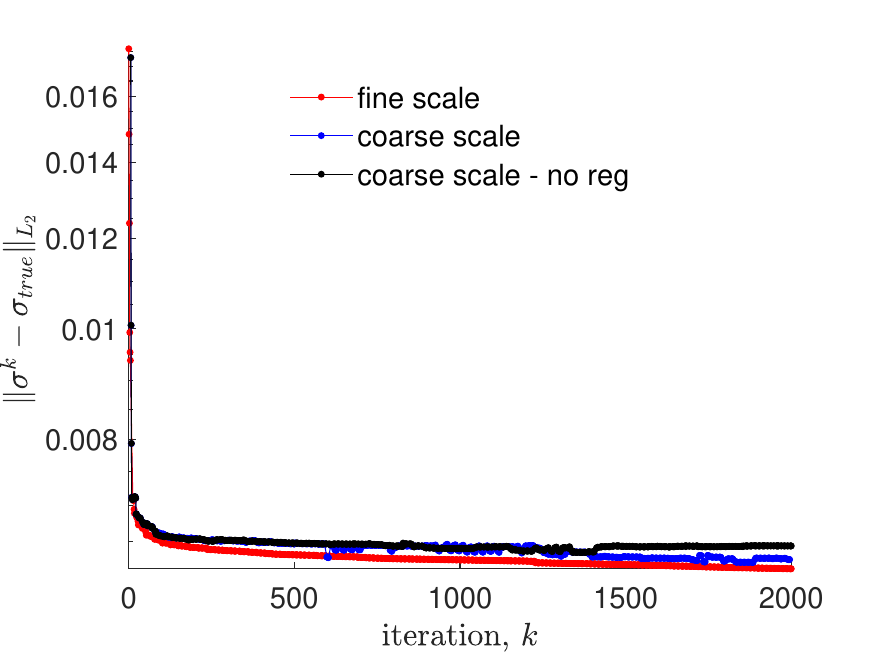}}}
  \mbox{
  \subfigure[$\hat \sigma$: $\beta_c = 10^{-5.25}$]
    {\includegraphics[width=0.33\textwidth]{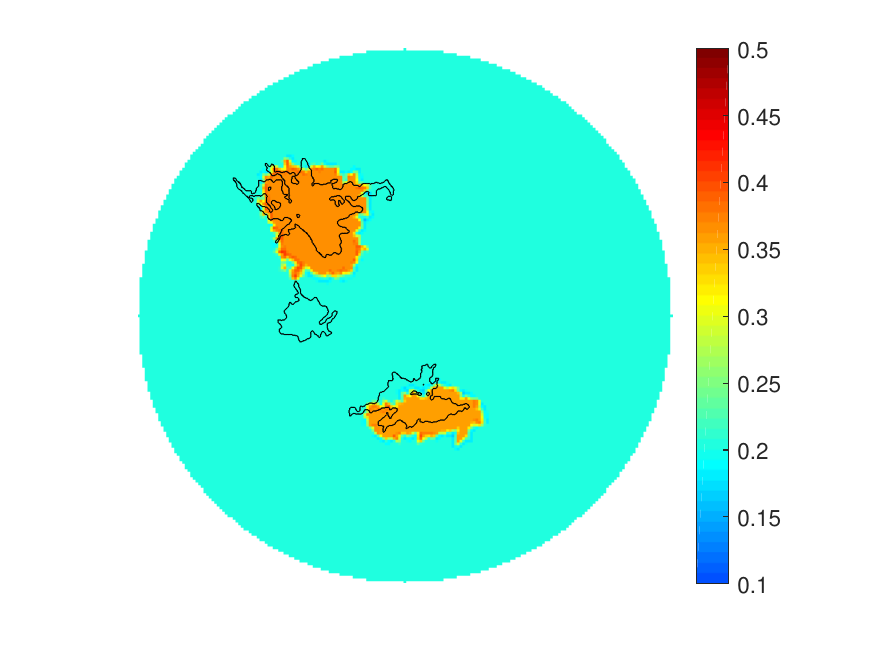}}
  \subfigure[objectives: $\beta_c = 10^{-5.25}$]
    {\includegraphics[width=0.33\textwidth]{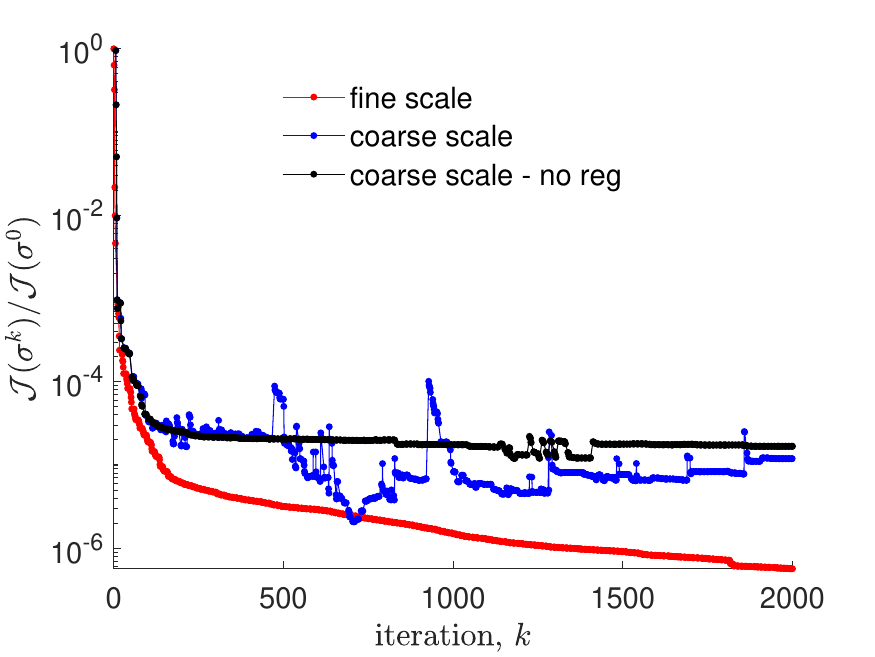}}
  \subfigure[solution error: $\beta_c = 10^{-5.25}$]
    {\includegraphics[width=0.33\textwidth]{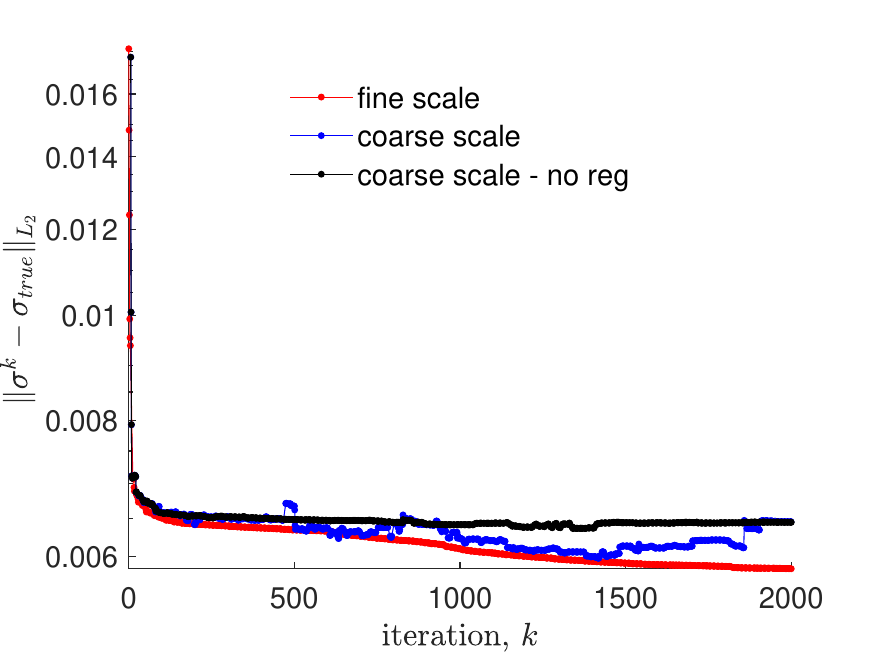}}}
  \end{center}
  \caption{Optimization results obtained with applied regularization at the coarse scale with weights
    (a-c)~$\beta_c = 10^{-8.5}$ and (d-f)~$\beta_c = 10^{-5.25}$ in \eqref{eq:reg_c}. (a,d)~The images of
    model~\#3 obtained after optimizing at multiple scales with the enhanced scale switching by the tuned PCA when
    $N_{\max} = 5$. The added black curves represent the boundary of the cancer-affected regions taken from known
    $\sigma_{true}(x)$ in Figure~\ref{fig:m16_result}(b). (b,e)~Normalized objective functions $\mJ(\sigma^k)/\mJ(\sigma^0)$
    and (c,f)~solution errors $\| \sigma^k - \sigma_{true} \|_{L_2}$ as functions of iteration count $k$ evaluated
    at fine (in red) and coarse (in blue) scales compared to the results in the coarse-scale images obtained with
    no regularization ($\beta_c = 0$) applied (black curves).}
  \label{fig:m16_reg_best}
\end{figure}

In the same fashion, as done before for model~\#2, we check the quality of the solutions obtained for $\beta_c \in
[10^{-10}, 10^{-5}]$ by choosing one value of $\beta_c = 10^{-8.5}$, with the corresponding outcomes shown in all plots
of Figure~\ref{fig:m16_reg_comp} as hexagons (left ones). The results of this multiscale optimization are provided
in Figures~\ref{fig:m16_reg_best}(a-c), showing the coarse-scale image in (a) with improved shapes. This improvement
reflects the presence of multiple (at least two) separate cancerous spots with a better resolution regarding the
estimated boundaries. Figures~\ref{fig:m16_reg_best}(b,c) also demonstrate some improvement in data fitting
and solution error compared to the image obtained without coarse-scale regularization. Finally, Figures~\ref{fig:m16_reg_best}(d-f)
represent another solution and associated measures used to quantify its quality obtained with $\beta_c = 10^{-5.25}$
(right hexagons on all plots of Figure~\ref{fig:m16_reg_comp}). Both measures, namely, objectives in (e) and solution
error in (f), characterize the choice of this $\beta_c$ as less favorable. However, the image in (d) shows better results for
improved shapes. This fact confirms the intricate complexity of model~\#3 and an evident necessity for further development
of the proposed methodology to improve its performance. We leave it as an open problem and discuss the directions for
making new steps in Section~\ref{sec:remarks}.

\section{Concluding Remarks}
\label{sec:remarks}

In this work, we proposed and validated an efficient computational framework for multiscale optimization supplied
with enhanced PCA-based multilevel parameterization. This framework is suitable for the optimal reconstruction of
physical properties (e.g., electrical conductivity in EIT imaging) of various media characterized by distributions
close to binary. For instance, we see this approach as useful in many applications in biomedical sciences to operate
with physical models supplied with some, possibly noisy, measurements. In particular, we explore the possibility of
applying the proposed solution methodology to the IPCD problems to detect defective (cancerous) regions surrounded
by healthy tissues for early cancer detection or easy control of the dynamics of cancer development or treatment
progress.

The core part of our computational framework is the gradient-based multiscale optimization supplied with multilevel
control space reduction. We propose flexible mechanisms for that reduction used interchangeably at fine and coarse
scales to enhance proper communication between solutions obtained at these scales to assure computational efficiency
and the superior quality of obtained results. The current state of this framework enables efficient solutions to
identify multiple regions as cancerous spots accurately characterized by their locations and the shapes of their
boundaries. The impact of the noise in measurements and the employment of regularization techniques are also
systematically analyzed in applications to the synthetic models and models based on real breast cancer images.
The proposed optimization algorithm has an easy-to-follow design tuned by a nominal number of parameters to govern
the entire suite of the computational facilities. In general, we see a high potential of the proposed computational
framework in minimizing possibilities for false positive and false negative screening and improving the overall
quality of the EIT-based procedures.

Despite the superior performance of the proposed framework, there are many ways this multiscale optimization algorithm
can be tested and further extended. We expect an even better performance by applying advanced minimization techniques
to perform local and global searches while optimizing at both fine and coarse scales, using adaptive schemes for flexible
switching between scales, implementing new procedures for finding optimal values of adjustable parameters, experimenting
with various scenarios for initial guesses, and establishing efficient termination criteria. It will be of interest to
involve a further analysis of the measurement structure, e.g., reviewing the value of information depending on the complexity
of models, boosting the efficiency of our rotation scheme to generate sufficient data, considering a 32-electrode scheme, and
improving sensitivity by optimizing the configuration of available data. Also, as many modern EIT systems feature pair-wise
voltage patterns, we will be interested in testing the performance of our new method in applications to such systems. Also
of interest is the extension of our multiscale optimization approach by including various sample structures while constructing
PCA and applying it to models characterized by bimodal and fully anisotropic distributions with irregular shapes and spiculated
margins. Finally, we believe this methodology has future potential in applications to a vast array of problems seen, e.g., in
biomedical sciences, physics, geology, and chemistry.

\subsubsection*{Acknowledgements}
We wish to thank the anonymous reviewer for their valuable comments and suggestions to improve the clarity of the presented
approach and the overall readability of this paper.

\subsubsection*{Funding information}
This study received financial support from Florida Tech's College of Engineering and Science 2021 Institutional Research
Incentive Grant.

\subsubsection*{Nomenclature}
\begin{table}[!h]
  \begin{tabular}{c l}
    $\Omega$ & physical domain \\
    $n$ & space dimension \\
    $x$ & spatial variable (independent) \\
    $\sigma(x)$ & electrical conductivity \\
    $U$ & electrical potential at electrode $E$\\
    $\ell$ & electrode number \\
    $m$ & total number of electrodes \\
    $E$ & boundary electrode \\
    $Z$ & contact impedance \\
    $I$ & electrical current \\
    $\partial \Omega$ & boundary of domain $\Omega$ \\
    $u(x)$ & electrical potential \\
    $I^*$ & measurement of electrical current $I$ \\
    $\mJ(\sigma)$ & objective function \\
    $\beta$ & scalarization weights in objective $\mJ(\sigma)$ \\
    $K$ & number of permutations of potentials within set $U$ \\
    $\hat \sigma(x)$ & optimal solution for $\sigma(x)$ \\
    $\psi(x)$ & adjoint PDE solution \\
    $\xi$ & fine scale reduced dimensional control vector \\
    $\hat \xi$ & optimal solution for control $\xi$ \\
    $\Phi$ & PCA linear transformation matrix \\
    $\bar \sigma(x)$ & prior mean for PCA transform \\
    $\sigma^*(x)$ & sample solution (realization) \\
    $N_r$ & number of realizations to construct matrix $\Phi$ \\
    $N_{\xi}$ & maximum number of principal components in truncated PCA \\
    $N_t$ & current number of principal components in truncated PCA
 \end{tabular}
\end{table}
\begin{table}[!t]
  \begin{tabular}{c l}
    $\hat N_t$ & optimal number of principal components in truncated PCA \\
    $\zeta$ & coarse scale reduced dimensional control vector \\
    $\hat \zeta$ & optimal solution for control $\zeta$ \\
    $N_{\zeta}$ & size of control $\zeta$ \\
    $\Delta$ & area of fine mesh elements \\
    $N$ & number of fine mesh elements \\
    $C$ & partitioning subsets \\
    $\mathcal{M}$ & partitioning map \\
    $P$ & partitioning (indicator) function \\
    $n_s$ & number of iterations between switching scales \\
    $k_s$ & count for switching cycles \\
    $\chi_c$ & coarse-scale indicator function \\
    $\epsilon_f$ & termination tolerance for fine scale \\
    $\epsilon_c$ & termination tolerance for coarse scale \\
    $\alpha_{c \rightarrow f}$ & coarse-to-fine projection parameter \\
    $\hat \alpha_{c \rightarrow f}$ & optimal value of parameter $\alpha_{c \rightarrow f}$ \\
    $\sigma_{PCA}$ & PCA equivalent for projected solution $\sigma(x)$ \\
    $N_{\max}$ & maximum number of expected cancer-affected regions \\
    $\delta_{\zeta}, \delta_{\xi}$ & finite difference scheme (or perturbation) parameters \\
    $r_{\Omega}$ & radius \\
    $w$ & electrode half-width \\
    $\sigma_{true}(x)$ & true electrical conductivity \\
    $\Omega_c$ & cancer affected region of $\Omega$ \\
    $\sigma_c$ & true electrical conductivity for $\Omega_c$ \\
    $\Omega_h$ & healthy tissue part of $\Omega$ \\
    $\sigma_h$ & true electrical conductivity for $\Omega_h$ \\
    $\mJ_c$ & penalization term for coarse scale optimization \\
    $\beta_c$ & coefficient for $\mJ_c$ \\
    $\kappa (\epsilon)$ & $\kappa$-test function \\
    $i,j,n$ & summation indices \\
    $k$ & current major iteration number
 \end{tabular}
\end{table}

\clearpage
\appendix
\section*{Appendix A}
\label{sec:app_1}

\begin{algorithm}[ht!]
\begin{algorithmic}
  \STATE $k \leftarrow 0$
  \STATE $\chi_c \leftarrow 0$
  \STATE $\sigma^0 \leftarrow $ initial guess $\sigma_0(x)$
  \STATE compute $\xi^0$ using $\sigma^0$ by \eqref{eq:map_pca_2}
  \REPEAT
  \STATE compute $u^k$ using $\sigma^k$ by solving \eqref{eq:forward}
  \STATE compute $\psi^k$ using $u^k$ and $\sigma^k$ by solving \eqref{eq:adjoint}
  \STATE compute $\bnabla_{\sigma} \mJ(\sigma^k)$ using $u^k$ and $\psi^k$
    by \eqref{eq:grad_sigma}
  \IF{$\chi_c = 1$}
    \STATE compute $\sigma(\xi^k)$ using $\xi_k$ by \eqref{eq:map_pca_1}
    \STATE compute $\bnabla_{\zeta} \mJ(\zeta^k)$ using $\zeta^k,
      \sigma(\xi^k),$ and $\bnabla_{\sigma} \mJ(\sigma^k)$ by
      \eqref{eq:upscale_map}--\eqref{eq:grad_zeta},
      \eqref{eq:sigma_coarse}, and  \eqref{eq:grad_zeta_a}--\eqref{eq:grad_zeta_b}
    \STATE update control $\zeta$ by computing $\zeta^{k+1}$ using descent directions based
      on $\bnabla_{\zeta} \mJ(\zeta^k)$
    \STATE update control $\sigma$ by computing $\sigma^{k+1}$ using $\zeta^{k+1}$, $\sigma(\xi^{k+1})$,
      and \eqref{eq:sigma_coarse}
  \ELSE
    \STATE compute $\bnabla_{\xi} \mJ(\xi^k)$ using $\bnabla_{\sigma}
      \mJ(\sigma^k)$ by \eqref{eq:grad_xi}
    \STATE update control $\xi$ by computing $\xi^{k+1}$ using descent directions based on
       $\bnabla_{\xi} \mJ(\xi^k)$
    \STATE update control $\sigma$ by computing $\sigma^{k+1}$ using $\xi^{k+1}$ and
      \eqref{eq:map_pca_1}
  \ENDIF
  \STATE $k \leftarrow k + 1$
  \STATE update $\chi_c$ using $k$ by \eqref{eq:scale_ind}
  \IF{$\chi_c(k) \neq \chi_c(k-1)$}
    \IF{$\chi_c = 1$}
      \STATE tune PCA with $N_{\xi}^k$ by solving \eqref{eq:pca_FtoC}
      \STATE update $\sigma^k$ by \eqref{eq:sigma_coarse} using $\xi^k$ updated by \eqref{eq:xi_trunc},
        $\zeta^k$, and $\sigma(\xi^k)$
    \ELSE
      \STATE tune PCA with $N_{\xi}^k$ by solving \eqref{eq:pca_CtoF}
      \STATE update $\xi^k$ using $N_{\xi}^k$, $\sigma^k$, and $\sigma(\xi^k)$ by
        \eqref{eq:xi_trunc}, \eqref{eq:proj_CtoF}, and solving \eqref{eq:proj_CtoF_rlx}
      \STATE update $\sigma^k$ using $\xi^k$ by \eqref{eq:map_pca_1}
    \ENDIF
  \ENDIF
  \STATE evaluate objective $\mJ(\sigma^k)$ by \eqref{eq:obj_fn_final}
  \UNTIL termination criterion \eqref{eq:termination} is satisfied
    to given tolerances $\epsilon_f$ and $\epsilon_c$
\end{algorithmic}
\caption{Workflow for multiscale optimization with multilevel parameterization}
\label{alg:main_opt}
\end{algorithm}

\newpage
\bibliographystyle{spmpsci}
\bibliography{biblio_Bukshtynov,biblio_EIT,biblio_OPT,biblio_Protas}

\end{document}